\documentclass{article}\begin{document}
\centerline{\bf On the Construction of Relativistic Quantum Wave Equation}
\centerline{\bf and }
\centerline{\bf General Solution of the Second Order Differential Equation}
\vskip .4in
\centerline{Nikolaos D. Bagis}
\centerline{Aristotele University of Thessaloniki AUTH}
\centerline{Thessaloniki, Greece}
\[
\]
\centerline{\bf Abstract}
Using the elementary axioms of special relativity and quantum mechanics we construct a wave equation which generalizes the Schrodinger equation. We also solve the general second order differential equation ($y''(x)=V(x)y(x)$).

\section{Introduction}

Einstein belived that the complete equivalence of mass and energy is given from the equation
\begin{equation}
E=mc^2,
\end{equation}
where $m$ is the mass of a moving particle and $E$ the total energy ($c$ is the velocity of light in vacum). 
This formula is not completely right. It has to do with error terms of series representing the mass of a moving particle (see relation (5) below) and restrictions using Newton formula to write 
\begin{equation}
m=m_0\gamma(v)=m_0+\frac{1}{c^2}\left(\frac{1}{2}m_0v^2\right)+\ldots
\end{equation} 
Taking the first two terms of (2) the energy of a moved particle $E=mc^2$ can divided into $E_0=m_0c^2$ (which $m_0$ is the rest mass) and $\frac{1}{2}m_0v^2$ is the kinetic energy.\\
\\
In 1927 Klein and Gordon manage to constuct a relativistic quantum mechanical equation for moving particles using the following equation
\begin{equation}
E^2=p^2c^2+m_0^2c^4.
\end{equation}
Here $p$ is the momentum of a moving particle and $m_0$ is the rest mass. By attaching where $E$ and $p$ the quantum operators $+i\hbar\partial_t$ and $-i\hbar\partial_x$ respectively (energy and momentum), they arived to an equation which bear their names (Klein-Gordon equation):
\begin{equation}
\partial^2_{xx}Y-\frac{1}{c^2}\partial^2_{tt}Y-\frac{m_0^2c^2}{\hbar^2}Y=0.
\end{equation}
In some cases the above equation holds good. But in general is not satisfactory.\\
However the folowing formula for the mass of a traveling particle is reliable
\begin{equation}
m=m_0\gamma(v)=\frac{m_0}{\sqrt{1-\left(\frac{v}{c}\right)^2}}
\end{equation}         
\\
In this article, under the same axioms used by Klein and Gordon and not using energy cuts and avoiding the square in energy we arive to a generalized equation similar to that of Schrodinger. Below we state the two axioms we need:\\
\\
\textbf{Axiom 1.}\\
Relation (5) holds for every moving particle.\\
\\
\textbf{Axiom 2.}\\
i) In the energy $E$ corespond the quantum operator $E=i\hbar\partial_t$\\
ii) In the momentum $p$ corspond the quantum operator $p=-i\hbar\partial_x$.

\section{Main Results}

The form of the equation we looking for must read as
\begin{equation}
EY=\left(f\left(p\right)+E_0 I\right)Y,
\end{equation}
where $f$ is a single valued function. $E_0=m_0c^2$, is the energy of the particle and the quantities $m_0,c$ are the rest mass and the velocity of light in vacum respectively.\\
We shall work with simple single-valued calculus and operators, to arive to an equation formed by (5) (Axiom 1), then we will use the quantum corespodence (Axiom 2) to produce the desired diferential equation.\\
Starting right away we have:
$$
E=\frac{p^2}{2m}+m_0c^2=\frac{m^2v^2}{2m}+m_0c^2=\frac{1}{2}mv^2+m_0c^2=\frac{1}{2}m_0\gamma(v)v^2+m_0c^2=
$$
\begin{equation}
=\frac{m_0c^2}{2\sqrt{1-\left(\frac{v}{c}\right)^2}}\left(\frac{v}{c}\right)^2+m_0c^2.
\end{equation} 
But
\begin{equation}
\frac{v}{c}=\frac{p}{mc}=\frac{p}{m_0c}\sqrt{1-\left(\frac{v}{c}\right)^2}
\end{equation}
We interested to expand $E$ in a form of $E=f(m_0,p)$ and then use Axiom (2). As someone can see this is quite dificult because of the nature of the relation $v=p/m$ in which $m$ depends by $v$ (relation (5)) and $v$ by $p/m$ again. In this point we mention that we can solve (8) with respect to $\frac{v}{c}$ to get  
\begin{equation}
\xi=\frac{v}{c}=\frac{p}{m_0c}\sqrt{1-\xi^2}.
\end{equation}       
Hence
\begin{equation}
\xi=\frac{v}{c}=\pm\frac{p}{m_0c\sqrt{1+\frac{p^2}{m_0^2c^2}}}.
\end{equation}
Seting this in (7) we get 
\begin{equation}
E=\frac{m_0c^2}{2}\frac{\xi^2}{\sqrt{1-\xi^2}}+m_0c^2=\frac{E_0}{2}\frac{p^2}{m_0^2c^2\sqrt{1+\frac{p^2}{m_0^2c^2}}}+E_0.
\end{equation}
Hence from the expansion
$$
\frac{x^2}{\sqrt{1+x^2}}=\sum^{\infty}_{n=0}\small\left(
\begin{array}{cc}
-1/2\\
n
\end{array}\right)\normalsize x^{2n+2}
$$
we get
$$
E=E_0+\frac{E_0}{2}\sum^{\infty}_{n=0}\small\left(
\begin{array}{cc}
-1/2\\
n
\end{array}\right)\normalsize \frac{p^{2n+2}}{(m_0c)^{2n+2}}.
$$
Replacing the arithmetical quntities of $p$ and $E$ with operators (Axiom 2), we get\\
\\
\textbf{Theorem 1.}
\small
\begin{equation}
i\hbar\frac{\partial Y}{\partial t}=E_0Y+\frac{E_0}{2}\sum^{\infty}_{n=0}(-1)^{n+1}\small\left(
\begin{array}{cc}
-1/2\\
n
\end{array}\right)\normalsize \left(\frac{\hbar}{m_0c}\right)^{2n+2}\frac{\partial^{2n+2}Y}{\partial x^{2n+2}}.
\end{equation}   
\normalsize
Provited that for all $(x,t)\in \textbf{R}\times\textbf{R}_{+}$, we have
$$
\lim_{n\rightarrow\infty}\left|\frac{\partial^{2n+2}_xY(x,t)}{\partial^{2n}_xY(x,t)}\right|=l^2\in\textbf{R}\eqno{(12.1)}
$$
and $l<1/l_0$, $l_0=\frac{\hbar}{m_0c}$.\\
\\
\textbf{Proof.}\\
Observe that (from the ratio test)
$$
r=l_0^2\frac{2n+1}{2n+2}\left|\frac{\partial^{2n+4}_{x}Y(x,t)}{\partial^{2n+2}_x Y(x,t)}\right|<1
$$

For solving (12), we use the folowing Lemma's found in [3]:\\
\\
\textbf{Lemma 1.}\\
If 
$$
f(x)=\sum^{\infty}_{n=1}A_ne^{-n x},
$$
then 
$$
f(x)=\sum^{\infty}_{n=1}\frac{1}{e^{nx}-1}\sum_{d|n}A_d\mu(n/d)
$$ 
and the oposite.\\
\\
\textbf{Lemma 2.}\\
If $x>0$
$$
f(x)=\sum^{\infty}_{n=1}\frac{1}{e^{nx}-1}\sum_{d|n}A_d\mu(n/d),
$$ 
then 
$$
f^{(\nu)}(x)=\sum^{\infty}_{n=1}\frac{1}{e^{nx}-1}\sum_{d|n}A_d(-d)^{\nu}\mu(n/d).
$$
\\

Assume now that a solution have the form 
$$
Y(x,t)=\sum^{\infty}_{k=1}\frac{\sum_{d|k}\chi_{d<1/l_0}A_d(t)\mu(k/d)}{e^{kx}-1},\eqno{(12.2)} 
$$
where $l_0<1$. Then
$$
i\hbar\sum^{\infty}_{k=1}\frac{1}{e^{kx}-1}\sum_{d|k}\chi_{d<1/l_0}A'_d(t)\mu(k/d)
=E_0\sum^{\infty}_{k=1}\frac{\sum_{d|k}\chi_{d<1/l_0}A_d(t)\mu(k/d)}{e^{kx}-1}+
$$
\begin{equation}
+\frac{E_0}{2}\sum^{\infty}_{n=0}(-1)^{n+1}C_{-1/2,n}l_0^{2n+2}\sum^{\infty}_{k=1}\frac{1}{e^{kx}-1}\sum_{d|k}\chi_{d<1/l_0}A_d(t)d^{2n+2}\mu(k/d).
\end{equation}
Rearaging the double series in the above equation and summing with respect to $n$ we get 
$$
i\hbar\sum^{\infty}_{k=1}\frac{1}{e^{kx}-1}\sum_{d|k}\chi_{d<1/l_0}A'_d(t)\mu(k/d)=
$$
\begin{equation}
=E_0\sum^{\infty}_{k=1}\frac{1}{e^{kx}-1}\sum_{d|k}\chi_{d<1/l_0}A_d(t)\mu(k/d)-
$$
$$
-\frac{E_0}{2}\sum^{\infty}_{k=1}\frac{1}{e^{k x}-1}\sum_{d|k}\chi_{d<1/l_0}A_d(t)\mu(k/d)\frac{l_0^2d^2}{\sqrt{1-l_0^2d^2}}.
\end{equation}
Hence it must be 
$$
\frac{i\hbar}{m_0c^2} A_d'(t)=\left(1-\frac{l_0^2d^2}{2\sqrt{1-l_0^2d^2}}\right) A_d(t)\textrm{, }d=1,2,\ldots,1/l_0
$$
or equivalently
$$
A_n(t)=f_n\exp\left(-\frac{iE_0t}{\hbar}+\frac{iE_0l_0^2n^2t}{2\hbar\sqrt{1-l_0^2n^2}}\right)\textrm{, for }n=1,2,\ldots,1/l_0.
$$
Hence
$$
Y(x,t)=\sum^{\infty}_{n=1}\frac{1}{e^{nx}-1}\sum_{d|n}f_d\chi_{d<1/l_0}\exp\left(-\frac{iE_0t}{\hbar}+\frac{iE_0l_0^2d^2t}{2\hbar\sqrt{1-l_0^2d^2}}\right)\mu(n/d).
$$
Returning to Taylor series, we arive to\\
\\
\textbf{Theorem 2.}
\begin{equation}
Y(x,t)=\exp\left(-\frac{iE_0t}{\hbar}\right)\sum^{1/l_0}_{n=1}f_n\exp\left(\frac{iE_0l_0^2n^2 t}{2\hbar\sqrt{1-l_0^2n^2}}\right)e^{-nx},
\end{equation}
where it is $Y(x,0)=\sum^{1/l_0}_{n=1}f_ne^{-nx}=f(x)$ and $l_0=\frac{\hbar}{m_0c}<1$ is the Compton's wavelength of the particle.\\
\\
\textbf{Notes.}\\
\textbf{i)} The Compton-Plank lenght $l_0$ along with (12.1), has the following meaning: If $l_0\geq 1$, then we can not have a solution of the form (12.2), as also the series becomes no convergent and equation (12) is no loger valid. For particles such as the neutrino's we have Compton length $l_{0}=\frac{\hbar}{m_{\nu}c}\approx 0.000164\textrm{ cm }$. Electron have $l_0=3,86\times 10^{-13}\textrm{ cm }$. This in literature is called Compton wavelength of electron.\\ 
\textbf{ii)} Taking the first term of the sum in (12)  we arive easily to Schrodinger's equation
\begin{equation}
i\hbar \frac{\partial Y}{\partial t}=E_0Y-\frac{\hbar^2}{2m_0}\frac{\partial^2 Y}{\partial x^2}.
\end{equation}
\\

Set now
\begin{equation}
C_{m,n}:=\small\left(
\begin{array}{cc}
m\\
n
\end{array}\right)\normalsize.
\end{equation}

$$
Y(t;x,y,z)=\sum^{\infty}_{k,l,m=1}\frac{\sum_{d_1|k}\sum_{d_2|l}\sum_{d_3|m}A_{d_1d_2d_3}(t)\mu\left(k/d_1\right)\mu(l/d_2)\mu(m/d_3)}{(e^{kx}-1)(e^{ly}-1)(e^{mz}-1)}.
$$
Hence
$$
\partial^{2\nu+2}_xY(t;x,y,z)=
$$
$$
\sum^{\infty}_{k,l,m=1}\frac{\sum_{d_1|k}\sum_{d_2|l}\sum_{d_3|m}A_{d_1d_2d_3}(t)d_1^{2\nu+2}\mu\left(k/d_1\right)\mu(l/d_2)\mu(m/d_3)}{(e^{kx}-1)(e^{ly}-1)(e^{mz}-1)},
$$
$$
\partial^{2\kappa+2}_yY(t;x,y,z)=
$$
$$
\sum^{\infty}_{k,l,m=1}\frac{\sum_{d_1|k}\sum_{d_2|l}\sum_{d_3|m}A_{d_1d_2d_3}(t)d_2^{2\kappa+2}\mu\left(k/d_1\right)\mu(l/d_2)\mu(m/d_3)}{(e^{kx}-1)(e^{ly}-1)(e^{mz}-1)},
$$
$$
\partial^{2\rho+2}_zY(t;x,y,z)=
$$
$$
\sum^{\infty}_{k,l,m=1}\frac{\sum_{d_1|k}\sum_{d_2|l}\sum_{d_3|m}A_{d_1d_2d_3}(t)d_3^{2\rho+2}\mu\left(k/d_1\right)\mu(l/d_2)\mu(m/d_3)}{(e^{kx}-1)(e^{ly}-1)(e^{mz}-1)}.
$$
Hence as in the one dimensional case we have
$$
i\hbar \sum^{\infty}_{k,l,m=1}\frac{1}{(e^{k x}-1)(e^{ly}-1)(e^{m z}-1)}\sum_{d_1|k,d_2|l,d_3|m}A_{d_1d_2d_3}'(t)\mu(k/d_1)\mu(l/d_2)\mu(m/d_3)=
$$
$$
=E_0\sum^{\infty}_{k,l,m=1}\frac{1}{(e^{kx}-1)(e^{l y}-1)(e^{m z}-1)}\sum_{d_1|k,d_2|l,d_3|m}A_{d_1d_2d_3}(t)\mu(k/d_1)\mu(l/d_2)\mu(m/d_3)\times
$$
$$
\times\left(1-\frac{l_0^2d_1^2}{2\sqrt{1-l_0^2 d_1^2}}-\frac{l_0^2d_2^2}{2\sqrt{1-l_0^2 d_2^2}}-\frac{l_0^2d_3^2}{2\sqrt{1-l_0^2 d_3^2}}\right).
$$
Hence it must be
$$
\frac{i\hbar}{E_0} A_{d_1d_2d_3}'(t)=\left(1-\frac{l_0^2d_1^2}{2\sqrt{1-l_0^2 d_1^2}}-\frac{l_0^2d_2^2}{2\sqrt{1-l_0^2 d_2^2}}-\frac{l_0^2d_3^2}{2\sqrt{1-l_0^2 d_3^2}}\right)A_{d_1d_2d_3}(t),
$$
or
$$
A_{d_1d_2d_3}(t)=
$$
$$
f_{d_1d_2d_3}\exp\left[\left(-i\frac{E_0 t}{\hbar}+i\frac{E_0l_0^2d_1^2t}{2\hbar\sqrt{1-l_0^2 d_1^2}}+i\frac{E_0l_0^2d_2^2t}{2\hbar\sqrt{1-l_0^2 d_2^2}}+i\frac{E_0l_0^2d_3^2t}{2\hbar\sqrt{1-l_0^2 d_3^2}}\right)\right].
$$
Hence the equation 
$$
i\hbar\frac{\partial Y}{\partial t}=E_0Y+\frac{E_0}{2}\sum^{\infty}_{n=0}(-1)^{n+1}C_{-1/2,n}\left(\frac{\hbar}{m_0c}\right)^{2n+2}\left(\frac{\partial^{2n+2}Y}{\partial x^{2n+2}}+\frac{\partial^{2n+2}Y}{\partial y^{2n+2}}+\frac{\partial^{2n+2}Y}{\partial z^{2n+2}}\right),
$$
admits solution
$$
Y(t;x,y,z)=\exp\left(-\frac{iE_0t}{\hbar}\right)\times
$$
$$
\times\sum^{1/l_0}_{k,l,m=1}f_{klm}\exp\left[\frac{iE_0l_0^2t}{2\hbar}\left(\frac{k^2}{\sqrt{\left|l_0^2k^2-1\right|}}+\frac{l^2}{\sqrt{\left|l_0^2l^2-1\right|}}+\frac{m^2}{\sqrt{\left|l_0^2m^2-1\right|}}\right)\right]e^{-kx-ly-mz}.
$$
QED\\
\\
\textbf{Theorem 3.}\\
Assume the equation
$$
i\hbar\frac{\partial Y}{\partial t}=\frac{E_0}{2}\sum^{\infty}_{n=0}(-1)^{n+1}C_{-1/2,n}\left(\frac{\hbar}{m_0c}\right)^{2n+2}\frac{\partial^{2n+2}Y}{\partial x^{2n+2}}+E_0Y+V(x)Y.\eqno{(17.1)}
$$
If the potential $V(x)$ is zero, then
\begin{equation}
Y(x,t)=\frac{\exp\left(-\frac{iE_0 t}{\hbar}\right)}{2\pi}\int^{1/l_0}_{-1/l_0}\widehat{Y}(\gamma,0)\exp\left(-\frac{iE_0 t}{2\hbar}\frac{l_0^2\gamma^2}{\sqrt{1+l_0^2\gamma^2}}\right)e^{i\gamma x}d\gamma.
\end{equation}
If the potential $V(x)$ is a real function of $x$ ($2\pi-$periodic) i.e. 
$$
V(x)=\sum^{\infty}_{k=0}a_kx^k=\sum^{\infty}_{k=-\infty}c_ke^{-ikx},
$$
then
\begin{equation}
i\hbar \frac{\partial}{\partial t}\widehat{Y}(\gamma,t)=\left(E_0+\frac{E_0}{2}\frac{l_0^2\gamma^2}{\sqrt{1+l_0^2\gamma^2}}+D\right)\widehat{Y}(\gamma,t),
\end{equation}
where
\begin{equation}
D(.)=\sum^{\infty}_{k=0}a_k i^k \left(\frac{\partial}{\partial \gamma}\right)^k(.)=V\left(i\frac{\partial}{\partial \gamma}\right)(.).
\end{equation}
Hence if we assume the inner product 
$$
\left(\Phi,\Psi\right)=\int^{\infty}_{-\infty}\overline{\Phi(\gamma,t)}\Psi(\gamma,t)d\gamma ,\eqno{(20.00)}
$$
the Hamiltonian is
$$
H=E_0+\frac{E_0}{2}\frac{l_0^2\gamma^2}{\sqrt{1+l_0^2\gamma^2}}+V\left(i\frac{\partial}{\partial \gamma}\right)\eqno{(20.01)}
$$
and is self adjoint in all Hilbert space $L_2(R)$. If $BL(R)$ denote the space of all Band Limited functions with Fourier transform zero outside of an interval $[-\sigma,\sigma]$, $\sigma>0$, then $BL(R)$ is Hilbert subspace of $L_2(R)$ and also dense in $L_2(R)$. Hence when the potential $V(x)$ is real and $2\pi-$periodic function (not limited to $2\pi$) we can write if $Y((.),t)\in BL(R)$:
$$
i\hbar\frac{\partial}{\partial t}\widehat{Y}(\gamma,t)=H\widehat{Y}(\gamma,t)\textrm{, }\gamma\in[-1/l_0,1/l_0]\eqno{(20.02)}
$$
and more generaly (20.02) hold in all $\gamma\in R$. The solution of (20.02) is
$$
\widehat{Y}(\gamma,t)=\exp\left(-\frac{i}{\hbar}tH\right)\widehat{Y}(\gamma,0).\eqno{(20.03)}
$$
\\
\textbf{Proof.}\\
Assuming the equation
\begin{equation}
i\hbar\frac{\partial Y}{\partial t}=\frac{E_0}{2}\sum^{\infty}_{n=0}(-1)^{n+1}C_{-1/2,n}\left(\frac{\hbar}{m_0c}\right)^{2n+2}\frac{\partial^{2n+2}Y}{\partial x^{2n+2}}+E_0Y+V(x)Y,
\end{equation}
where $Y=Y(x,t)$ and $V(x)=\sum^{\infty}_{k=0}a_kx^k$. Then taking the Fourier transform in both sides of (21), we get
$$
i\hbar \frac{\partial \widehat{Y}(\gamma,t)}{\partial t}=\frac{E_0}{2}\sum^{\infty}_{n=0}(-1)^{n+1}C_{-1/2,n}l_0^{2n+2}(i\gamma)^{2n+2}\widehat{Y}(\gamma,t)+E_0\widehat{Y}(\gamma,t)+
$$
$$
+\int^{\infty}_{-\infty}Y(x,t)V(x)e^{-ix\gamma}dx\textrm{, }|\gamma|\leq 1/l_0\eqno{(i)}
$$
and $\widehat{Y}(\gamma,t)=\int^{+\infty}_{-\infty}Y(x,t)e^{-ix \gamma}dx$. The only way to hold $(i)$ is to have $\widehat{Y}(\gamma,t)=0$, when $\left|\gamma\right|>1/l_0$ i.e. $Y(x,t)$ must be a band limited function. Hence we write equivalent
$$
i\hbar\frac{\partial \widehat{Y}(\gamma,t)}{\partial t}=\left(E_0+\frac{E_0}{2}\frac{l_0^2\gamma^2}{\sqrt{1+l_0^2\gamma^2}}\right)\widehat{Y}(\gamma,t)+\sum^{\infty}_{n=0}a_ni^n\frac{\partial^n \widehat{Y}(\gamma,t)}{\partial \gamma^n}.\eqno{(21.1)}
$$
Hence if we denote the differential operator
\begin{equation}
D(.)=\sum^{\infty}_{n=0}a_n i^n \left(\frac{\partial}{\partial \gamma}\right)^n(.)=V\left(i\frac{\partial}{\partial \gamma}\right)(.),
\end{equation}
we get (19). Where we have make use of the periodicity of $V(x)$, ($V(x)$ is bounded function of $x$, since it is periodic function in all $R$) and the fact that $Y((.),t)$ is in $L_2(R)$:
$$
\int^{\infty}_{-\infty}Y(x,t)V(x)e^{-i x \gamma}dx=\int^{\infty}_{-\infty}Y(x,t)\sum^{\infty}_{n=-\infty}c_ne^{-in x}e^{-ix \gamma}dx=
$$
$$
=\sum^{\infty}_{n=-\infty}c_n\int^{\infty}_{-\infty}Y(x,t)e^{-i x(\gamma+n)}dx=\sum^{\infty}_{n=-\infty}c_n\widehat{Y}(\gamma+n,t)=
$$
$$
=\sum^{\infty}_{n=-\infty}c_n\exp\left(n\frac{\partial}{\partial \gamma}\right)\widehat{Y}(\gamma,t)=V\left(i\frac{\partial}{\partial \gamma}\right)\widehat{Y}(\gamma,t).
$$
For to solve (19), we write $\widehat{Y}(\gamma,t)=A_n(\gamma)B_n(t)$. Then
$$
i\hbar A_n(\gamma)B_n'(t)=\left(E_0+\frac{E_0}{2}\frac{l_0^2\gamma^2}{\sqrt{1+l_0^2\gamma^2}}\right)A_n(\gamma)B_n(t)+B_n(t)DA_n(\gamma).
$$
Hence
$$
i\hbar \frac{B_n'(t)}{B_n(t)}=E^{*}_n\textrm{ and }DA_n(\gamma)=\left(E^{*}_n-E_0-\frac{E_0}{2}\frac{l_0^2\gamma^2}{\sqrt{1+l_0^2\gamma^2}}\right)A_n(\gamma).
$$
Hence
$$
B_n(t)=\Lambda_n\exp\left(\frac{-iE^{*}_n}{\hbar} t\right)\eqno{(22.1)}
$$
and
$$
\sum^{\infty}_{k=0}a_ki^k\left(\frac{\partial}{\partial \gamma}\right)^kA_n(\gamma)=\left(E^{*}_n-E_0-\frac{E_0}{2}\frac{l_0^2\gamma^2}{\sqrt{1+l_0^2\gamma^2}}\right)A_n(\gamma)\Leftrightarrow
$$
$$
V\left(i\frac{\partial}{\partial \gamma}\right)A_n(\gamma)=A_n(\gamma)g_n(\gamma),\eqno{(eq)}
$$
where 
$$
g_{n}(\gamma)=E^{*}_n-E_0-\frac{E_0}{2}\frac{l_0^2\gamma^2}{\sqrt{1+l_0^2\gamma^2}}.\eqno{(22.2)}
$$
This equation $(eq)$ is hold for general potential $V(x)$. Assuming that
$$
V(x)=\sum^{\infty}_{k=-\infty}c_ke^{-ikx},
$$
we have
$$
V\left(i\frac{\partial}{\partial \gamma}\right)A_n(\gamma)=\sum^{\infty}_{k=-\infty}c_k\exp\left(k\frac{\partial}{\partial \gamma}\right)A_n(\gamma)=\sum^{\infty}_{k=-\infty}c_kA_n(\gamma+k).
$$
Hence it is interesting to solve the equation
\begin{equation}
\sum^{\infty}_{k=-\infty}c_kA_m(\gamma+k)=g_m(\gamma)A_m(\gamma)\textrm{, }m\in\textbf{Z}.
\end{equation}
We solve the equation:
\begin{equation}
\sum^{\infty}_{k=-\infty}c_kA_m(\gamma+k)=g_m(\gamma)A_m(\gamma).
\end{equation}
We have
$$
\sum^{\infty}_{k=-\infty}c_kA_m(\gamma+k)=\sum^{\infty}_{k=-\infty}\delta_{k,0}g_m(\gamma)A_m(\gamma+k)\Leftrightarrow
$$
$$
\sum^{\infty}_{k=-\infty}A_m(\gamma+k)\left(c_k-\delta_{k,0}g_m(\gamma)\right)=0.
$$
If
$$
A_{m}(\gamma)=\sum_{l\in Z^{*}}b_{m,l}(\gamma)
$$
and
$$
b_{m,l}(\gamma+k)=b^{*}_{m,l}(\gamma)G_m(kl,\gamma),
$$
we get
$$
\sum^{\infty}_{k=-\infty}c_kA_m(\gamma+k)=c_0A_{m}(\gamma)+\sum_{k\in Z^{*}}c_kA_m(\gamma+k)=
$$
$$
=c_0\sum_{l\in Z^{*}}b_{m,l}(\gamma)+\sum_{(k,l)\in Z^{*}\times Z^{*}}c_kb_{m,l}\left(\gamma+k\right)=
$$
$$
=c_0G_m(0,\gamma)\sum_{l\in Z^{*}}b^{*}_{m,l}(\gamma)+\sum_{(k,l)\in Z^{*}\times Z^{*}}c_kb^{*}_{m,l}(\gamma)G(kl,\gamma)=
$$
$$
=c_0G_m(0,\gamma)\sum_{l\in Z^{*}}b^{*}_{m,l}(\gamma)+\sum_{n\in Z^{*}}\left(\sum_{kl=n}c_kb^{*}_{m,l}(\gamma)\right)G_m(n,\gamma).
$$
Hence
$$
c_0G_m(0,\gamma)\sum_{l\in Z^{*}}b^{*}_{m,l}(\gamma)+\sum_{n\in Z^{*}}G_m(n,\gamma)\left(\sum_{kl=n}c_kb^{*}_{m,l}(\gamma)\right)=
$$
$$
=\sum_{n\in Z}G_m(n,\gamma)\left(\sum_{\scriptsize\begin{array}{cc}kl=n\\ l\neq 0\end{array}\normalsize}c_kb^{*}_{m,l}(\gamma)\right)
=\sum_{n\in Z}\sum_{l\in Z^{*}}\delta_{n,0}g_m(\gamma)b_{m,l}(\gamma+n)\Leftrightarrow
$$
$$
\sum^{\infty}_{n=-\infty}G_m(n,\gamma)\left(\sum_{\scriptsize\begin{array}{cc}kl=n\\ l\neq 0\end{array}\normalsize}c_kb^{*}_{m,l}(\gamma)\right)=\sum_{n\in Z}\sum_{l\in Z^{*}}\delta_{n,0}g_m(\gamma)b^{*}_{m,l}(\gamma)G_m(nl,\gamma)=
$$
$$
=\sum^{\infty}_{n=-\infty}G_m(n,\gamma)\sum_{\scriptsize\begin{array}{cc}kl=n\\l\neq 0\end{array}\normalsize}b^{*}_{m,l}(\gamma)\delta_{k,0}g_m(\gamma)\Leftrightarrow
$$
$$
\sum^{\infty}_{n=-\infty}G_m(n,\gamma)\left(\sum_{\scriptsize\begin{array}{cc}kl=n\\ l\neq 0\end{array}\normalsize}c_kb^{*}_{m,l}(\gamma)\right)=
$$
$$
=\sum^{\infty}_{n=-\infty}G_m(n,\gamma)\left(\sum_{\scriptsize\begin{array}{cc}kl=n\\ l\neq 0\end{array}\normalsize}b^{*}_{m,l}(\gamma)\delta_{k,0}g_m(\gamma)\right)\Leftrightarrow
$$
$$
\sum^{\infty}_{n=-\infty}G_m(n,\gamma)\left(\sum_{\scriptsize\begin{array}{cc}kl=n\\ l\neq 0\end{array}\normalsize}b^{*}_{m,l}(\gamma)\left(c_k-\delta_{k,0}g_m(\gamma)\right)\right)=0.
$$
Hence
$$
\sum_{\scriptsize\begin{array}{cc}kl=n\\ l\neq 0\end{array}\normalsize}b^{*}_{m,l}(\gamma)\left(c_k-g_m(\gamma)\delta_{k,0}\right)=0.\eqno{(24.01)}
$$
We set
$$
a_{m,k}(\gamma)=c_k-g_m(\gamma)\delta_{k,0}.
$$
Then if $a_{m,k}(\gamma)$ are known, the arithmetic inverse of $a_{m,k}(\gamma)$, is $b^{*}_{m,k}(\gamma)$ in the sense
$$
\sum_{\scriptsize\begin{array}{cc}kl=n\\l\neq 0\end{array}\normalsize}b^{*}_{m,l}(\gamma)a_{m,k}(\gamma)=0\Leftrightarrow b^{*}_{m,n}(\gamma)=\{a_{m,n}(\gamma)\}^{(-1)}\eqno{(24.1)}
$$
Equivalent we can find $b^{*}_{m,l}(\gamma)$ from the system: 
$$
\sum_{l\in Z^{*}}b^{*}_{m,l}(\gamma)=0\eqno{(co1)}
$$
and
$$
\sum_{\scriptsize\begin{array}{cc}abs(d)|n\\n\neq 0\end{array}\normalsize}b^{*}_{m,d}(\gamma)c_{n/d}=0\eqno{(co2)}
$$
Thus
\begin{equation}
A_m(\gamma)=G_m(0,\gamma)\sum^{\infty}_{n=-\infty}\{c_n-g_m(\gamma)\delta_{n,0}\}^{(-1)}.
\end{equation}
Then the solution $\widehat{Y}(\gamma,t)$ of (21.1) is
\begin{equation}
\widehat{Y}(\gamma,t)=\sum^{\infty}_{m=-\infty}\Lambda_{m}G_m(0,\gamma)\exp\left(\frac{-iE_m t}{\hbar}\right)\sum^{\infty}_{n=-\infty}\{c_n-g_m(\gamma)\delta_{n,0}\}^{(-1)}.
\end{equation}
Now we show that $H$ is a self adjoint operator in all Hilbert space
$$
\left(\Phi,H\Psi\right)=\int^{\infty}_{-\infty}\overline{\Phi(\gamma,t)}\left(E_0+\frac{E_0}{2}\frac{l_0^2\gamma^2}{\sqrt{1+l_0^2\gamma^2}}+V\left(i\frac{\partial}{\partial\gamma}\right)\right)\Psi(\gamma,t)d\gamma=
$$
$$
=\int^{\infty}_{-\infty}\overline{\Phi(\gamma,t)}\left(E_0+\frac{E_0}{2}\frac{l_0^2\gamma^2}{\sqrt{1+l_0^2\gamma^2}}\right)\Psi(\gamma,t)d\gamma+
$$
$$
+\int^{\infty}_{-\infty}\overline{\Phi(\gamma,t)}V\left(i\frac{\partial}{\partial\gamma}\right)\Psi(\gamma,t)d\gamma=
$$
$$
=\int^{\infty}_{-\infty}\overline{\left(E_0+\frac{E_0}{2}\frac{l_0^2\gamma^2}{\sqrt{1+l_0^2\gamma^2}}\right)\Phi(\gamma,t)}\Psi(\gamma,t)d\gamma+
$$
$$
+\int^{\infty}_{-\infty}\overline{\Phi(\gamma,t)}\sum^{\infty}_{k=-\infty}c_k\exp\left(k\frac{\partial}{\partial\gamma}\right)\Psi(\gamma,t)d\gamma=
$$
$$
=\int^{\infty}_{-\infty}\overline{\left(E_0+\frac{E_0}{2}\frac{l_0^2\gamma^2}{\sqrt{1+l_0^2\gamma^2}}\right)\Phi(\gamma,t)}\Psi(\gamma,t)d\gamma+
$$
$$
+\sum^{\infty}_{k=-\infty}c_k\int^{\infty}_{-\infty}\overline{\Phi(\gamma,t)}\Psi(\gamma+k,t)d\gamma=
$$
$$
=\int^{\infty}_{-\infty}\overline{\left(E_0+\frac{E_0}{2}\frac{l_0^2\gamma^2}{\sqrt{1+l_0^2\gamma^2}}\right)\Phi(\gamma,t)}\Psi(\gamma,t)d\gamma+
$$
$$
+\sum^{\infty}_{k=-\infty}c_k\int^{\infty}_{-\infty}\overline{\Phi(\gamma-k,t)}\Psi(\gamma,t)d\gamma=
$$
$$
=\int^{\infty}_{-\infty}\overline{\left(E_0+\frac{E_0}{2}\frac{l_0^2\gamma^2}{\sqrt{1+l_0^2\gamma^2}}\right)\Phi(\gamma,t)}\Psi(\gamma,t)d\gamma+
$$
$$
+\int^{\infty}_{-\infty}\overline{\sum^{\infty}_{k=-\infty}\overline{c_k}\exp\left(-k\frac{\partial}{\partial \gamma}\right)\Phi(\gamma,t)}\Psi(\gamma,t)d\gamma=
$$
$$
=\int^{\infty}_{-\infty}\overline{\left(E_0+\frac{E_0}{2}\frac{l_0^2\gamma^2}{\sqrt{1+l_0^2\gamma^2}}\right)\Phi(\gamma,t)}\Psi(\gamma,t)d\gamma+
$$
$$
+\int^{\infty}_{-\infty}\overline{\sum^{\infty}_{k=-\infty}\overline{c_{-k}}\exp\left(k\frac{\partial}{\partial \gamma}\right)\Phi(\gamma,t)}\Psi(\gamma,t)d\gamma.
$$
Hence assuming that $\overline{c_{-k}}=c_k$, we get
$$
\left(\Phi,H\Psi\right)=
$$
$$
=\int^{\infty}_{-\infty}\overline{\left(E_0+\frac{E_0}{2}\frac{l_0^2\gamma^2}{\sqrt{1+l_0^2\gamma^2}}+\sum^{\infty}_{k=-\infty}c_k\exp\left(k\frac{\partial}{\partial\gamma}\right)\right)\Phi(\gamma,t)}\Psi(\gamma,t)d\gamma=
$$
$$
=\int^{\infty}_{-\infty}\overline{\left(E_0+\frac{E_0}{2}\frac{l_0^2\gamma^2}{\sqrt{1+l_0^2\gamma^2}}+V\left(i\frac{\partial}{\partial\gamma}\right)\right)\Phi(\gamma,t)}\Psi(\gamma,t)d\gamma=\left(H\Phi,\Psi\right).
$$
But we have asumed that $\overline{c_{-k}}=c_k$. Hence equivalent 
$$
\sum^{\infty}_{k=-\infty}\overline{c_{-k}}e^{-ikx}=\sum^{\infty}_{k=-\infty}c_ke^{-ikx}\Leftrightarrow
$$
$$
\overline{\sum^{\infty}_{k=-\infty}c_{-k}e^{ikx}}=\sum^{\infty}_{k=-\infty}c_ke^{-ikx}\Leftrightarrow
$$
$$
\overline{\sum^{\infty}_{k=-\infty}c_{k}e^{-ikx}}=\sum^{\infty}_{k=-\infty}c_ke^{-ikx}\Leftrightarrow\overline{V(x)}=V(x).
$$
Hence $V(x)$ must be real for $x$ real. QED\\
\\
\textbf{Notes I.}\\
In the most general case we have
$$
i\hbar\frac{\partial Y(x,t)}{\partial t}=E_0Y(x,t)+\frac{E_0}{2}\sum^{\infty}_{n=0}f_n i^n\frac{\partial^{n}Y(x,t)}{\partial x^{n}}+V(x)Y(x,t).\eqno{(26.01)}
$$
We assume in (26.01) that convergence of the sum is uniform and the function $f(x)=\sum^{\infty}_{n=0}f_n x^n$, is such that $f:R \rightarrow R$, analytic in $R$. Then  
$$
H=\frac{E_0}{2}f\left(i\frac{\partial}{\partial x}\right)+E_0 I+V(x) I.\eqno{(26.02)}
$$
Also we set $V(x)$ to be any real and smooth periodic function. That is because we know operator $V\left(i\frac{\partial}{\partial x}\right)$, is self adjoint on any function of $L_2(R)$, (see the proof of Theorem 3 above). Hence the equation 
$$
i\hbar\frac{\partial}{\partial t}Y(x,t)=H Y(x,t),\eqno{(26.03)}
$$
is the semigroup equation of our system and the operator $H$ is independed of $t$. We also have
$$
Y(x,t)=e^{-i/\hbar t H}Y(x,0).
$$
However the operator $H$ is self adjoint.\\
\\
\textbf{Conjecture.}\\
For general functions $f(x),g(x)$, we have the following continuation 
$$
\int^{\infty}_{-\infty}f(x)g(x)e^{-ix \gamma}dx=f\left(i\frac{\partial}{\partial \gamma}\right)\widehat{g}(\gamma)=g\left(i\frac{\partial}{\partial \gamma}\right)\widehat{f}(\gamma).
$$
\\

Taking the Fourier transform in both sides of $(26.01)$ we have (using and the above notes)
$$
i\hbar\frac{\partial \widehat{Y}(\gamma,t)}{\partial t}=E_0 \widehat{Y}(\gamma,t)+\frac{E_0}{2}\sum^{\infty}_{n=0}f_n i^n \int^{\infty}_{-\infty}\frac{\partial^{n}Y(x,t)}{\partial x^{n}}e^{-i x \gamma}dx+
$$
$$
+\int^{\infty}_{-\infty}V(x)Y(x,t)e^{-ix\gamma}dx\Leftrightarrow
$$
$$
i\hbar\frac{\partial \widehat{Y}(\gamma,t)}{\partial t}=E_0 \widehat{Y}(\gamma,t)+\frac{E_0}{2}\sum^{\infty}_{n=0}f_n i^n \sum^{n-1}_{j=0}\left[\left(\frac{\partial^{n-1-j}}{\partial x^{n-1-j}}Y(x,t)\right)e^{-ix\gamma}\right]^{x=+\infty}_{x=-\infty}(i\gamma)^j+
$$
$$
+\frac{E_0}{2}\sum^{\infty}_{n=0}f_n (-\gamma)^n\widehat{Y}(\gamma,t)+V\left(i\frac{\partial}{\partial \gamma}\right)\widehat{Y}(\gamma,t)\Leftrightarrow
$$
$$
E\widehat{Y}(\gamma,t)=\left(E_0+\frac{E_0}{2} f(-\gamma)+V\left(i\frac{\partial}{\partial \gamma}\right)\right)\widehat{Y}(\gamma,t)+
$$
$$
+\frac{E_0}{2}\sum^{\infty}_{n=0}f_n i^n\sum^{n-1}_{j=0}\left[\frac{\partial^{n-1-j}}{\partial x^{n-1-j}}Y(x,t)e^{-ix\gamma}\right]^{x=+\infty}_{x=-\infty}(i\gamma)^j\Leftrightarrow
$$
$$
E\widehat{Y}(\gamma,t)=\left(E_0+\frac{E_0}{2} f(-\gamma)+V\left(i\frac{\partial}{\partial \gamma}\right)\right)\widehat{Y}(\gamma,t)+
$$
$$
+\frac{E_0}{2}\sum^{\infty}_{n=0}f_n i^n\sum^{n-1}_{j=0}\left[Y^{(n-1-j,0)}(h,t)e^{-ih\gamma}\right]^{h=+\infty}_{h=-\infty}(i\gamma)^j\Leftrightarrow
$$
$$
E\widehat{Y}(\gamma,t)=\left(E_0+\frac{E_0}{2} f(-\gamma)+V\left(i\frac{\partial}{\partial \gamma}\right)\right)\widehat{Y}(\gamma,t)+
$$
$$
+\frac{E_0}{2}\left[e^{-ih\gamma}\sum^{\infty}_{n=0}f_n i^n \sum^{n-1}_{j=0}Y^{(j,0)}(h,t)(i\gamma)^{n-1-j}\right]^{h=+\infty}_{h=-\infty}.\eqno{(26.04)}
$$
However it holds
$$
\sum^{n-1}_{j=0}Y^{(j,0)}(h,t)(i\gamma)^{(n-1)-j}=(i\gamma)^{-1}\sum^{n}_{j=0}Y^{(j,0)}(h,t)(i\gamma)^{n-j}-
$$
$$
-(i\gamma)^{-1}Y^{(n,0)}(h,t).
$$
Hence
$$
\left[e^{-ih\gamma}\sum^{\infty}_{n=0}f_n i^n \sum^{n-1}_{j=0}Y^{(j,0)}(h,t)(i\gamma)^{n-1-j}\right]^{h=+\infty}_{h=-\infty}=
$$
$$
=(i\gamma)^{-1}\sum^{\infty}_{n=0}f_n i^n\sum^{n}_{j=0}\left[Y^{(j,0)}(h,t)e^{-ih\gamma}\right]^{h=+\infty}_{h=-\infty}(i\gamma)^{n-j}-
$$
$$
-(i\gamma)^{-1}\sum^{\infty}_{n=0}f_n i^n \left[Y^{(n,0)}(h,t)e^{-ih\gamma}\right]^{h=+\infty}_{h=-\infty}=
$$
$$
=(i\gamma)^{-1}\left[e^{-ih\gamma}\sum^{\infty}_{n=0}f_n i^n \sum^{n}_{j=0}\epsilon_{n-j}Y^{(j,0)}(h,t)(i\gamma)^{n-j}\right]^{h=+\infty}_{h=-\infty},\eqno{(26.05)}
$$
where $\epsilon_{n}=1$, if $n=1,2,3,\ldots$ and $\epsilon_0=0$. Hence
$$
\left(\sum^{\infty}_{k=0}Y^{(k,0)}(h,t)x^k\right)\frac{i\gamma x}{1-i\gamma x}=\sum^{\infty}_{n=0}x^n\sum^{n}_{j=0} Y^{(j,0)}(h,t)\epsilon_{n-j}(i\gamma)^{n-j}.
$$
Hence if we define the $(TY)_n(\gamma,t)$ from the relation
$$
\left[\frac{x e^{-ih\gamma}}{1-i\gamma x}\sum^{\infty}_{k=0}Y^{(k,0)}(h,t)x^k\right]^{h=+\infty}_{h=-\infty}=\sum^{\infty}_{n=0}x^n(TY)_n(\gamma,t),\eqno{(26.06)}
$$
then
$$
\left[e^{-ih\gamma}\sum^{\infty}_{n=0}f_n i^n \sum^{n-1}_{j=0}Y^{(j,0)}(h,t)(i\gamma)^{n-1-j}\right]^{h=+\infty}_{h=-\infty}=
$$
$$
=(i\gamma)^{-1}\left[e^{-ih\gamma}\sum^{\infty}_{n=0}f_n i^n \sum^{n}_{j=0}\epsilon_{n-j}Y^{(j,0)}(h,t)(i\gamma)^{n-j}\right]^{h=+\infty}_{h=-\infty}=
$$
$$
=\sum^{\infty}_{n=0}f_n i^n (TY)_n(\gamma,t).\eqno{(26.07)}
$$
Note also that
$$
(TY)_n(\gamma,t)=\left[(i\gamma)^{-1}e^{-ih\gamma}\sum^{n}_{j=0}Y^{(n-j,0)}(h,t)\epsilon_j(i\gamma)^{j}\right]^{h=+\infty}_{h=-\infty}\Leftrightarrow
$$
$$
(TY)_n(\gamma,t)=(i\gamma)^{-1}\sum^{n}_{j=0}\epsilon_jA_{n-j}(t)(i\gamma)^j.\eqno{(26.08)}
$$
\\
\textbf{Theorem 3.01}\\
Assume the Notes I for the equation
$$
i\hbar\frac{\partial Y(x,t)}{\partial t}=E_0Y(x,t)+\frac{E_0}{2}\sum^{\infty}_{n=0}f_n i^n \frac{\partial^{n}Y(x,t)}{\partial x^{n}}+V(x)Y(x,t).\eqno{(26.09)}
$$
Then taking the Fourier transform in both sides of $(26.09)$, we find
$$
E\widehat{Y}(\gamma,t)=\left(E_0+\frac{E_0}{2} f(-\gamma)+V\left(i\frac{\partial}{\partial \gamma}\right)\right)\widehat{Y}(\gamma,t)+\sum^{\infty}_{n=0}f_n i^n (TY)_n(\gamma,t),\eqno{(26.10)}
$$
where the $(TY)_n(\gamma,t)$ are defined from the relation
$$
\left[\frac{xe^{-ih\gamma}}{1-i\gamma x}\sum^{\infty}_{k=0}Y^{(k,0)}(h,t)x^k\right]^{h=+\infty}_{h=-\infty}=\sum^{\infty}_{n=0}x^n(TY)_n(\gamma,t)\eqno{(26.11)}
$$
and it holds
$$
(TY)_n(\gamma,t)=\left[(i\gamma)^{-1} e^{-ih\gamma}\sum^{n}_{j=0}Y^{(n-j,0)}(h,t)\epsilon_j(i\gamma)^{j}\right]^{h=+\infty}_{h=-\infty}.\eqno{(26.12)}
$$
The symbol $\epsilon_n$ is 1 if $n\neq 0$ and 0 otherwise. Also $f(x)=\sum^{\infty}_{n=0}f_n x^n$.\\
\\ 
\textbf{Lemma.}\\
It holds in general (for any operator $U$ and any function $\Psi$):
$$
i\hbar\left(\partial_t U\right)\Psi=\left[E,U\right]\Psi.
$$
\\
\textbf{Theorem 3.1}\\
We assume the equation
$$
i\hbar\frac{\partial}{\partial t}\Psi(x,t)=H\Psi(x,t),\eqno{(eq)}
$$
where the ''Hamiltonian'' is
$$
H=f(x)+V\left(i\frac{\partial}{\partial x}\right).
$$
If $V(x)=\sum_{n}c_ne^{-i n x}$, is real periodic (and smooth) potential (here limited to $[0,2\pi]$ but not have to in generaly), then $H$ is self adjoint in all Hilbert space $L_2(R)$ with inner product
$$
(X,Y)=\int^{\infty}_{-\infty}\overline{X(x,t)}Y(x,t)dx.
$$
The solution of $(eq)$ is given from the following eigenvalue problem (since $H$ is self adjoint):
$$
E\psi_k(x,t)=H\psi_k(x,t)=E_k\psi_k(x,t).
$$
Hence
$$
f(x)\psi_k(x,t)+V\left(i\frac{\partial}{\partial x}\right)\psi_k(x,t)=E_k\psi_k(x,t)
$$
and
$$
\psi_k(x,t)=\exp\left(-iE_kt/\hbar\right)\phi_k(x).
$$
Hence
$$
\Psi(x,t)=\sum^{\infty}_{k=-\infty}C_k\psi_k(x,t)=\sum^{\infty}_{k=-\infty}C_k\exp\left(-\frac{iE_k t}{\hbar}\right)\phi_k(x).
$$ 
Also then
$$
V\left(i\frac{\partial}{\partial x}\right)\phi_k(x)=\sum^{\infty}_{n=-\infty}c_n\phi_k(x+n)=g_k(x)\phi_k(x),
$$
where
$$
g_k(x)=E_k-f(x).
$$
Hence
$$
V\left(i\frac{\partial}{\partial x}\right)\phi_k(x)=\left(E_k-f(x)\right)\phi_k(x).
$$
The functions $\phi_k(x)$ and $\psi_k(x,t)$ are both orthonormal bases:
$$
\left(\phi_k,\phi_l\right)=\delta_{k,l}\textrm{ and }\left(\psi_k,\psi_l\right)=\delta_{k,l}.
$$
We define the sequences of functions $b_{k,l}(x)$, $b^{*}_{k,l}(x)$ and $G_k(l,x)$, such that
$$
b_{k,l}(x+n)=b^{*}_{k,l}(x)G_k(nl,x)
$$
and
$$
\phi_k(x)=\sum_{l\in Z^{*}}b_{k,l}(x).
$$
Also given the arithmetical functions (with respect to $k$) $\alpha_k(x)$, we define the arithmetical functions $\{\alpha_n(x)\}^{(-1)}$ as
$$
\sum_{\scriptsize\begin{array}{cc}kl=n\\ l\neq 0\end{array}\normalsize}\alpha_k(x)\beta_l(x)=0\Leftrightarrow \beta_k(x)=\{\alpha_k(x)\}^{(-1)}.
$$
This is equivalent to say
$$
\sum_{l\in Z^{*}}\beta_l(x)=0
$$
and
$$
\sum_{0<abs(d)|n}\alpha_{n/d}(x)\beta_d(x)=0.
$$
However in our problem we have
$$
\sum_{\scriptsize\begin{array}{cc}kl=n\\ l\neq 0\end{array}\normalsize}b^{*}_{m,l}(x)\left(c_k-g_m(x)\delta_{k,0}\right)=0
$$
or equivalent
$$
 b^{*}_{m,k}(x)=\{c_k-g_m(x)\delta_{k,0}\}^{(-1)}.\eqno{(co)}
$$
With the above definitions the general problem reduces to:
$$
\sum_{l\in Z^{*}}b^{*}_{m,l}(x)=0\eqno{(co1)}
$$
and
$$
\sum_{\scriptsize\begin{array}{cc}0<abs(d)|n\\n\neq 0\end{array}\normalsize}c_{n/d}b^{*}_{m,d}(x)=0.\eqno{(co2)}
$$
Also writing $Y_k(x,t)=A_k(x)B_k(t)$, we get
$$
Y(x,t)=\sum^{\infty}_{k=-\infty}C_k\psi_k(x,t)=\sum^{\infty}_{k=-\infty}\Lambda_kY_k(x,t)=
$$
$$
=\sum^{\infty}_{m=-\infty}\Lambda_{m}\exp\left(\frac{-iE_m t}{\hbar}\right)G_m(0,x)\sum^{\infty}_{n=-\infty}\{c_n-g_m(x)\delta_{n,0}\}^{(-1)}.
$$
\\
\textbf{Theorem 3.2}\\
Assume Notes I and equation $E\Psi=H_0\Psi : (eq)$, where 
$$
H_0=E_0+\frac{E_0}{2}f\left(i\frac{\partial}{\partial x}\right)+V\left(x\right),
$$ 
and $f(x),V(x)$ also as in Notes I. Taking the Fourier transform (definition)
$$
\widehat{Y}(\gamma,t):=\int^{\infty}_{-\infty}Y(x,t)e^{-ix\gamma}dx,
$$
in both sides of $(eq)$ along with condition
$$
\sum^{\infty}_{n=0}f_ni^n\left(TY\right)_n(\gamma,t)=0\textrm{, }\forall \gamma,t\in R,\eqno{(26.13)}
$$ 
where the function $(TY)_n(\gamma,t)$ is as defined in (26.10),(26.11),(26.12) of Theorem 3.01, the $(eq)$ becomes
$$
E\widehat{Y}(\gamma,t)=\left(E_0+\frac{E_0}{2}f(-\gamma)+V\left(i\frac{\partial}{\partial \gamma}\right)\right)\widehat{Y}(\gamma,t).\eqno{(26.14)}
$$
Here equation (26.14) have Hamiltonian 
$$
H=E_0+\frac{E_0}{2}f(-\gamma)+V\left(i\frac{\partial}{\partial \gamma}\right),\eqno{(26.15)}
$$
which is self adjoint. Hence we can solve $(eq)$ using Theorem 3.1.:
$$
\widehat{Y}(\gamma,t)=
\sum^{\infty}_{m=-\infty}\Lambda_{m}\exp\left(\frac{-iE^{*}_m t}{\hbar}\right)G_m(0,\gamma)\sum^{\infty}_{n=-\infty}\{c_n-g_m(\gamma)\delta_{n,0}\}^{(-1)},
$$
where
$$
g_m(\gamma)=E^{*}_m-E_0-\frac{E_0}{2}f(-\gamma).
$$
\\
\textbf{Proof.}\\
Since we can assume any (arbitrary) function we want for $f(x)$, we may only prove that the operator  $V\left(i\partial_x\right)$, is self adjoint. Assume that $\Phi=\Phi(x,t)$ and $\Psi=\Psi(x,t)$ are any functions of $L_2(R)$. Then
$$
\left(\Phi,V\left(i\frac{\partial}{\partial x}\right)\Psi\right)=\int_{R}\overline{\Phi(x,t)}\sum^{\infty}_{n=-\infty}c_n\exp\left(n\frac{\partial}{\partial x}\right)\Psi(x,t) dx=
$$
$$
=\int_{R}\overline{\Phi(x,t)}\sum^{\infty}_{n=-\infty}c_n\Psi\left(x+n,t\right)dx=\sum^{\infty}_{n=-\infty}c_n\int_{R}\overline{\Phi(x,t)}\Psi(x+n,t)dx=
$$
$$
=\sum^{\infty}_{n=-\infty}c_n\int_{R}\overline{\Phi\left(x-n,t\right)}\Psi(x,t)dx=
$$
$$
=\int_{R}\overline{\sum^{\infty}_{n=-\infty}\overline{c_n}\exp\left(-n\frac{\partial}{\partial x}\right)\Phi(x,t)}\Psi(x,t)dx=
$$
$$
=\int_{R}\overline{\sum^{\infty}_{n=-\infty}c_{-n}\exp\left(-n\frac{\partial}{\partial x}\right)\Phi(x,t)}\Psi(x,t)dx
=
$$
$$
=\int_{R}\overline{\sum^{\infty}_{n=-\infty}c_n\exp\left(n\frac{\partial}{\partial x}\right)\Phi(x,t)}\Psi(x,t)dx=
$$
$$
=\int_{R}\overline{V\left(i\frac{\partial}{\partial x}\right)\Phi(x,t)}\Psi(x,t)dx=\left(V\left(i\frac{\partial}{\partial x}\right)\Phi,\Psi\right).
$$
QED.\\
\\
 
We also have the next\\
\\
\textbf{Theorem 4.1}\\
Assuming the equation
\begin{equation}
y''(x)=V(x)y(x),
\end{equation}
where $V(x)=\sum^{\infty}_{k=-\infty}c_ke^{-ik x}$, we can define the functions $f_l(\gamma)$, $f^{*}_l(\gamma)$, $G(k,\gamma)$, such that:
$$
f_{l}(\gamma+k)=f^{*}_{l}(\gamma)G(k l,\gamma)
$$
and
\begin{equation}
\widehat{y}(\gamma)=\sum_{l\in Z^{*}}f_l(\gamma)=G(0,\gamma)\sum_{l\in Z^{*}}f^{*}_l(\gamma),
\end{equation}
where $f^{*}_n(\gamma)=\{c_n+\gamma^2\delta_{n,0}\}^{(-1)}$ is the arithmetic inverse of $a_{n}(\gamma)=c_n+\gamma^2\delta_{n,0}$ in the sense
$$
\sum_{\scriptsize \begin{array}{cc}
lk=n\\
l\neq 0
\end{array}\normalsize}f^{*}_{l}(\gamma)a_{k}(\gamma)=0.
$$
\\
\textbf{Proof.}\\
Assume that $V(x)=\sum^{\infty}_{n=-\infty}c_ne^{-inx}$. We have
$$
\sum^{\infty}_{k=0}F_ky^{(k)}(x)=V(x)y(x)\Leftrightarrow 
$$
$$
\sum^{\infty}_{k=0}F_k\int^{\infty}_{-\infty}y^{(k)}(x)e^{-ix\gamma}dx=\int^{\infty}_{-\infty}y(x)\left(\sum^{\infty}_{n=-\infty}c_ne^{-in x}\right)e^{-ix\gamma}dx\Leftrightarrow
$$
$$
\sum^{\infty}_{k=0}F_k\cdot (i\gamma)^k\int^{\infty}_{-\infty}y(x)e^{-ix\gamma}dx=\sum^{\infty}_{n=-\infty}c_n\int^{\infty}_{-\infty}y(x)e^{-in x}e^{-i x\gamma}dx\Leftrightarrow
$$
$$
F(i\gamma)\widehat{y}(\gamma)=\sum^{\infty}_{n=-\infty}c_n\exp\left(n\frac{d}{d\gamma}\right)\left(\int^{\infty}_{-\infty}y(x)e^{-ix\gamma}dx\right)\Leftrightarrow
$$
$$
F(i\gamma)\widehat{y}(\gamma)=\sum^{\infty}_{n=-\infty}c_n\widehat{y}(\gamma+n).
$$
Equivalent
$$
F(i\gamma)\widehat{y}(\gamma)=V\left(i\frac{d}{d\gamma}\right)\widehat{y}(\gamma).\eqno{(30.01)}
$$
Hence working as above with $F(x)=x^2$, we get
$$
y''(x)=V(x)y(x)\Leftrightarrow V\left(i\frac{d}{d\gamma}\right)\widehat{y}(\gamma)=-\gamma^2 \widehat{y}(\gamma) \Leftrightarrow
$$
$$
-\gamma^2 \widehat{y}(\gamma)=\sum^{\infty}_{n=-\infty}c_n\widehat{y}(\gamma+n).
$$
Hence
$$
\widehat{y}(\gamma)=\sum_{l\in Z^{*}}f_l(\gamma)
$$
and
$$
\widehat{y}(\gamma)=G(0,\gamma)\sum^{\infty}_{n=-\infty}\{c_n+\gamma^2\delta_{n,0}\}^{(-1)}.
$$
QED.\\
\\

The equation
$$
\sum^{\infty}_{n=0}g_n (-i)^n y^{(n)}(x)=V(x) y(x),\eqno{(30.1)}
$$
where $g_n=\frac{g^{(n)}(0)}{n!}$, can be written in the form
$$
\sum^{\infty}_{n=0}g_n \gamma^n\widehat{y}(\gamma)=V\left(i\frac{d}{d\gamma}\right)\widehat{y}(\gamma).
$$ 
Hence
$$
V\left(i\frac{d}{d\gamma}\right)\widehat{y}(\gamma)=g(\gamma) \widehat{y}(\gamma),
$$
which have solution
$$
\sum^{\infty}_{n=-\infty}c_n\widehat{y}(\gamma+n)=g(\gamma)\widehat{y}(\gamma).
$$
Hence one solution is
$$
\widehat{y}(\gamma)=\frac{1}{g(\gamma)}\sum_{n_1}\frac{c_{n_1}}{g(\gamma+n_1)}\sum_{n_2}\frac{c_{n_2}}{g(\gamma+n_1+n_2)}\sum_{n_3}\frac{c_{n_3}}{g(\gamma+n_1+n_2+n_3)}\ldots .
$$
Hence holds the next\\
\\
\textbf{Theorem 4.2}\\
The equation
$$
\sum^{\infty}_{n=0}g_n(-i)^n\frac{d^n}{dx^n}y(x)=V(x)y(x),
$$
where 
$$
V(x)=\sum^{\infty}_{n=-\infty}c_ne^{-2\pi inx/T},
$$
is any $T-$periodic function, have solution $y(x)$ with Fourier transform $\widehat{y}(\gamma)$, such that
$$
\sum^{\infty}_{n=-\infty}c_n\widehat{y}\left(\gamma+\frac{2\pi n}{T}\right)=g(\gamma)\widehat{y}(\gamma).
$$
The function $g$ is $g(z)=\sum^{\infty}_{n=0}g_n z^n$. Also in case $T=2\pi$ (not limited to), we have
$$
\sum^{\infty}_{n=0}b_n\frac{d^n}{dx^n}y(x)=V(x)y(x),
$$
where 
$$
V(x)=\sum^{\infty}_{n=-\infty}c_ne^{-inx}
$$
and
$$
\sum^{\infty}_{n=-\infty}c_n\widehat{y}\left(\gamma+n\right)=B(i\gamma)\widehat{y}(\gamma).
$$
$$
\widehat{y}(\gamma)=\frac{1}{B(i\gamma)}\sum_{n_1}\frac{c_{n_1}}{B(i\gamma+in_1)}\sum_{n_2}\frac{c_{n_2}}{B(i\gamma+i n_1+in_2)}\sum_{n_3}\frac{c_{n_3}}{B(i\gamma+i n_1+i n_2+in_3)}\ldots,
$$
where
$$
B(x)=\sum^{\infty}_{n=0}b_nx^n.
$$
The Fourier transform of $y(x)$ is defined as:
$$
\widehat{y}(\gamma)=\int^{\infty}_{-\infty}y(t)e^{-it\gamma}dt.
$$
\\

Set $A_n(z)$ and $\theta_l(z)$  such that
$$
\theta_l(z)=\sum^{\infty}_{n=-\infty}A_n(z)\theta_l(z+n)\eqno{(eq1)}
$$
and
$$
\sum^{\infty}_{k=-\infty}c_kA_{n-k}(z+k)=\delta_{n,0} g(z),\eqno{(eq2)}
$$
we have
$$
\sum^{\infty}_{k=-\infty}c_k\theta_l(z+k)=\sum^{\infty}_{k=-\infty}c_k\sum^{\infty}_{n=-\infty}A_n(z+k)\theta_l(z+n+k)=
$$
$$
=\sum^{\infty}_{k=-\infty}c_k\sum^{\infty}_{n=-\infty}A_{n-k}(z+k)\theta_l(z+n)=
$$
$$
=\sum^{\infty}_{n=-\infty}\left(\sum^{\infty}_{k=-\infty}c_kA_{n-k}(z+k)\right)\theta_l(z+n)=
$$
$$
=\sum^{\infty}_{n=-\infty}\delta_{n,0}g(z)\theta_l(z+n)=g(z)\theta_l(z).
$$
Hence we search for functions $\theta_l(z)$ and $A_n(z)$ such that $(eq1)$ and $(eq2)$ hold. But equation $(eq2)$ can be written as
\begin{equation}
\sum^{\infty}_{n=-\infty}\sum^{\infty}_{k=-\infty}c_kA_{n-k}(z+k)e^{inx}=g(z)\sum^{\infty}_{n=-\infty}\delta_{n,0}e^{inx}.
\end{equation}
Hence
$$
\sum^{\infty}_{n,k=-\infty}c_{k}A_{n}(z+k)e^{inx}e^{ikx}=g(z).
$$
Taking in both sides the Fourier transform we find
$$
\sum^{\infty}_{n=-\infty}\sum^{\infty}_{k=-\infty}c_k\widehat{A}_{n}(\gamma)e^{i\gamma k}e^{inx}e^{ikx}=\widehat{g}(\gamma)\Leftrightarrow
$$
$$
\overline{V}(x+\gamma)\sum^{\infty}_{n=-\infty}\widehat{A}_n(\gamma)e^{inx}=\widehat{g}(\gamma)\Leftrightarrow
$$
$$
\sum^{\infty}_{n=-\infty}\widehat{A}_n(\gamma)e^{inx}=\frac{\widehat{g}(\gamma)}{\overline{V}(x+\gamma)}.
$$
Hence
$$
\widehat{A}_n(\gamma)=\widehat{g}(\gamma)\frac{1}{2\pi}\int^{2\pi}_{0}\frac{e^{-itn}}{\overline{V}(t+\gamma)}dt\Leftrightarrow
$$
$$
A_n(x)=\frac{1}{4\pi^2}\int^{\infty}_{-\infty}\widehat{g}(\gamma)\left(\int^{2\pi}_{0}\frac{e^{-itn}}{\overline{V}(t+\gamma)}dt\right) e^{ix\gamma}d\gamma.
$$
Hence knowing $A_n(x)$, then from $(eq1)$ we can find $\theta_l(z)$ by iteration i.e.
$$
\theta_l(z)=\sum_{n_1}A_{n_1}(z)\sum_{n_2}A_{n_2}(z+n_1)\sum_{n_3}A_{n_3}(z+n_1+n_2)\ldots
$$
and we have the next\\
\\
\textbf{Theorem 5.}\\
\textbf{I.} Assume that
\begin{equation}
V(x)=\sum^{\infty}_{n=-\infty}c_ne^{-inx}
\end{equation}
and $g(x)$ is a given function. Then equation
\begin{equation}
V\left(i\frac{d}{d w}\right)\widehat{y}(w)=\sum^{\infty}_{n=-\infty}c_n\widehat{y}(w+n)=g(w)\widehat{y}(w),
\end{equation}
have solution 
\begin{equation}
\widehat{y}(w)=\sum_{n_1}A_{n_1}(w)\sum_{n_2}A_{n_2}(w+n_1)\sum_{n_3}A_{n_3}(w+n_1+n_2)\ldots,
\end{equation}
where $A_n(x)$ is such that
\begin{equation}
A_n(x)=\frac{1}{4\pi^2}\int^{\infty}_{-\infty}\widehat{g}(\gamma)\left(\int^{2\pi}_{0}\frac{e^{-itn}}{\overline{V}(t+\gamma)}dt\right) e^{ix\gamma}d\gamma
\end{equation}
and $\overline{V}(x)=V(x)=\sum^{\infty}_{n=-\infty}c_ne^{-in x}$, $\sum^{\infty}_{n--\infty}|c_n|<\infty$.\\
\textbf{II.} When 
\begin{equation}
V(x)=\sum^{\infty}_{n=-\infty}c_ne^{-in x}
\end{equation}
and $g(x)=-x^2$, then equation
\begin{equation}
y''(x)=V(x)y(x),
\end{equation}
have solution $y(x)$ such that its Fourier transform is
\begin{equation}
\widehat{y}(w)=\sum_{n_1}A_{n_1}(w)\sum_{n_2}A_{n_2}(w+n_1)\sum_{n_3}A_{n_3}(w+n_1+n_2)\ldots,
\end{equation}
where
\begin{equation}
A_n(x)=\frac{1}{2\pi}\int^{\infty}_{-\infty}\delta''(\gamma)\left(\int^{2\pi}_{0}\frac{e^{-itn}}{\overline{V}(t+\gamma)}dt\right) e^{ix\gamma}d\gamma.
\end{equation}
\\
\textbf{Remarks.}\\ 
\textbf{1)} The Fourier transform of $f(x)$ is defined as
\begin{equation}
\widehat{f}(x)=\int^{\infty}_{-\infty}f(t)e^{-itx}dt\Leftrightarrow f(x)=\frac{1}{2\pi}\int^{\infty}_{-\infty}\widehat{f}(\gamma)e^{i\gamma x}d\gamma.
\end{equation}
\textbf{2)} The function $\delta(x)$ is the Dirac delta function, while $\delta_{k,\nu}$ is the Kronecker's delta symbol, i.e. $\delta_{k,\nu}=1$ if $k=\nu$ and $\delta_{k,\nu}=0$ if $k\neq \nu$.\\
\textbf{3)} From 
$$
\frac{1}{2\pi}\int^{\infty}_{-\infty}\delta''(w)F(w)e^{iw x}dw=\frac{1}{2\pi}\left(-x^2F_n(0)+2ixF_n'(0)+F_n''(0)\right),
$$
we get
$$
A_n(x)=\frac{1}{2 \pi}\int^{2\pi}_{0}\left(-\frac{x^2}{\overline{V}(t)}-\frac{2ix\overline{V}'(t)}{\overline{V}(t)^2}+\frac{2\overline{V}'(t)^2}{\overline{V}(t)^3}-\frac{\overline{V}''(t)}{\overline{V}(t)^2}\right)e^{-it n}dt.
$$
\\
\textbf{Example 1.}\\
Assume $V(x)=1$, then $c_n=\delta_{n,0}$. Hence from (39) we get 
$$
A_n(x)=\frac{1}{2\pi}\int^{\infty}_{-\infty}\delta''(\gamma)2\pi\delta_{n,0}e^{i\gamma x}d\gamma=-x^2\delta_{n,0}.\eqno{(a)}
$$
Hence $\theta(z)=-z^2\theta(z)$. Hence one leads to $\theta(z)=\delta(z-i)$. Hence $\widehat{y}(z)=\delta(z-i)$, which gives $y(x)=\frac{1}{2\pi}\int^{\infty}_{-\infty}\delta(z-i)e^{-izx}dz=\frac{e^{x}}{2\pi}\textrm{ : }(a)$. Hence indeed one solution of $y''(x)=V(x)y(x)$, with $V(x)=1$ is $y(x)=\frac{e^{x}}{2\pi}$. However one observes that working with functions of $L_2(\textbf{R})$ is better than working with distributions. In this example the Fourier integral $(a)$ is over the real line and the distribution is $\delta(x-i)$. Hence $x-i$ is never $0$ at the real line. However many books use this Fourier pair (actually we need here a generalization of the Fourier transform or the notion of distribution functions).\\
\\
\textbf{Example 2.}\\
Assume $V(x)=e^{-ix}$. Hence $c_n=\delta_{n,1}$. In this case $\overline{V}(x)=e^{ix}$ and
$$
\frac{1}{2\pi}\int^{2\pi}_{0}\frac{e^{-it n}}{e^{i(t+\gamma)}}dt=2\pi e^{-i\gamma}\delta_{n,-1}.
$$
Hence
$$
A_n(x)=\frac{1}{2\pi}\int^{\infty}_{-\infty}\delta''(\gamma)2\pi e^{-i\gamma}\delta_{n,-1}e^{i\gamma x}d\gamma=-(x-1)^2\delta_{n,-1}.
$$
Hence from
$$
\sum^{\infty}_{n=-\infty}c_n\theta(z+n)=-z^2\theta(z)\Leftrightarrow \theta(z+1)=-z^2\theta(z),
$$
we get easily
$$
\theta(z)=e^{-i\pi z}\Gamma(z)^2.
$$
Therefore, there exists a solution of $y''(x)=e^{-ix}y(x)$ such that its Fourier transform is
$$
\widehat{y}(w)=e^{-i\pi w}\Gamma(w)^2.
$$
However this evaluation is not in $L_2(\textbf{R})$.\\
\\
\textbf{Remarks.}\\
If we assume that $y(it)=(Mg)(t)$, where $(Mg)(t)$ is the Mellin transform of $g(t)$, then the Fourier transform of a function is
$$
\widehat{y}(x)=\int^{+\infty}_{-\infty}y(t)e^{-i t x}dt=\int^{+\infty/i}_{-\infty/i}y(it)e^{-i (it)x}idt=
$$
$$
=-i\int^{+i\infty}_{-i\infty}y(it)e^{tx}dt=-i 2 \pi i\frac{1}{2\pi i}\int^{+i\infty}_{-i\infty}Mg(t)\left(e^{-x}\right)^{-t}dt=
$$
$$
=2\pi g\left(e^{-x}\right).
$$
Hence $\widehat{y}\left(x\right)=2\pi g\left(e^{-x}\right)$ iff $y\left(ix\right)=(Mg)(x)$.\\
\\
\textbf{Example 3.}\\
Assume that
$$
g(x)=\frac{1}{x^2+1},
$$
then
$$
(Mg)(s)=\frac{\pi}{2}\csc\left(\frac{\pi s}{2}\right)\textrm{, }0<Re(s)<2.
$$
Thus
$$
\widehat{y}(\gamma)=\frac{2\pi}{1+e^{-2\gamma}}
$$
and
$$
y(x)=\frac{\pi}{2}\csc\left(\frac{-i\pi x}{2}\right).
$$

\section{About the complete solution of the equation with periodic potentials and some notes on special cases and the equation $y''(x)=V(x)y(x)$.}

Assume now that 
\begin{equation}
\frac{1}{V(x)}=\sum^{\infty}_{n=-\infty}c^{*}_ne^{-i n x}.
\end{equation}
Then
$$
\frac{1}{2\pi}\int^{2\pi}_{0}\frac{e^{-i t n}}{\overline{V}(t+\gamma)}dt=\frac{1}{2\pi}\int^{2\pi}_{0}\left(\sum^{\infty}_{k=-\infty}c^{*}_k e^{-it n}e^{i k (t+\gamma)}\right)dt=
$$
$$
=\sum^{\infty}_{k=-\infty}c^{*}_k\left(\frac{1}{2\pi}\int^{2\pi}_{0}e^{it(k-n)}dt\right)e^{ik\gamma}
=\sum^{\infty}_{k=-\infty}c^{*}_k\delta_{k,n}e^{ik\gamma}=c^{*}_ne^{in\gamma}
$$
and
$$
A_n(x)=\int^{\infty}_{-\infty}\delta''(\gamma)c^{*}_ne^{in\gamma}e^{ix\gamma}d\gamma=-(n+x)^2c^{*}_n.
$$
\\
\textbf{Theorem 5.1}\\
If we denote $c_n^{*}$ such that
$$
\frac{1}{V(x)}=\sum^{\infty}_{n=-\infty}c^{*}_ne^{-i n x},\eqno{(41.1)}
$$
then\\
\textbf{I.} If $g(x)=-x^2$, then
$$
A_n(x)=-(n+x)^2c^{*}_{n}.\eqno{(41.2)}
$$
\textbf{II.} If $g(x)$ is any analytic function, then
$$
A_n(x)=g(n+x)c^{*}_{n}\eqno{(41.3)}
$$
and the solution of 
$$
\sum_{n}c_n\widehat{y}(z+n)=g(z)\widehat{y}(z),\eqno{(41.4)}
$$
is
$$
\widehat{y}(z)=\lim_{N\rightarrow\infty}\sum_{n_1}\sum_{n_2}\ldots\sum_{n_N}g(n_1+z)g(n_1+n_2+z)\ldots
$$
$$
\ldots g(n_1+n_2+\ldots+n_{N}+z)c^{*}_{n_1}c^{*}_{n_2}\ldots c^{*}_{n_N}.\eqno{(41.5)}
$$
\\

Using (41.2) we get
$$
\widehat{y}(w)=\lim_{N\rightarrow\infty}\sum_{n_1}\sum_{n_2}\sum_{n_3}\ldots\sum_{n_N}(-1)^N(n_1+w)^2(n_2+n_1+w)^2(n_3+n_2+n_1+w)^2\ldots
$$
$$
\ldots (n_N+n_{N-1}+\ldots+n_2+n_1+w)^2c^{*}_{n_1}c^{*}_{n_2}\ldots c^{*}_{n_N}=
$$
$$
=\lim_{N\rightarrow\infty}\int^{\infty}_{-\infty}\int^{\infty}_{-\infty}\ldots\int^{\infty}_{-\infty}\delta''(t_1)\delta''(t_2)\ldots\delta''(t_N)\times 
$$
$$
\times \frac{1}{V(t_1+t_2+\ldots+t_N)}\frac{1}{V(t_2+t_3+\ldots+t_N)}\ldots\frac{1}{V(t_{N-1}+t_N)}\frac{1}{V(t_N)}\times 
$$
$$
\times e^{-it_1w}e^{-it_2w}\ldots e^{-it_Nw}dt_1dt_2\ldots dt_N\Rightarrow
$$	
$$
\widehat{y}(w)=\lim_{N\rightarrow\infty}\int^{\infty}_{-\infty}\int^{\infty}_{-\infty}\ldots\int^{\infty}_{-\infty}\delta''(t_1)\delta''(t_2)\ldots\delta''(t_N)\times 
$$
$$
\times \frac{1}{V(t_1)}\frac{1}{V(t_1+t_2)}\ldots\frac{1}{V(t_1+t_2+\ldots+t_N)}\times 
$$
$$
\times e^{-it_1w}e^{-it_2w}\ldots e^{-it_Nw}dt_1dt_2\ldots dt_N=
$$
$$
=\lim_{N\rightarrow\infty}\frac{\partial^2}{\partial h_1^2}\frac{\partial^2}{\partial h_2^2}\ldots\frac{\partial^2}{\partial h_N^2}\left(f_{N}(h_1,h_2,\ldots,h_N)e^{-ih_1w-ih_2 w-\ldots-ih_N w}\right)_{h_1=h_2=\ldots=h_N=0},
$$
where
$$
f_N(t_1,t_2,\ldots,t_N)=\frac{1}{V(t_1)V(t_1+t_2)\ldots V(t_1+t_2+\ldots+t_N)}.
$$
Also if we denote
$$
\epsilon_{n_1,n_2,\ldots,n_N}:=2^{n_1 mod 2+n_2 mod 2+\ldots+n_N mod 2}.
$$
and
$$
f_{N}^{(n_1,n_2,\ldots,n_N)}=\left(\frac{\partial^{n_1+n_2+\ldots+n_N}}{\partial t_1^{n_1}\partial t_2^{n_2}\ldots\partial t_N^{n_N}}f_N(t_1,t_2,\ldots,t_N)\right)_{t_1=t_2=\ldots=t_N=0},
$$
then
$$
\widehat{y}(w)=\lim_{N\rightarrow\infty}\sum^{2}_{n_1,n_2,\ldots,n_N=0} (-1)^{n_1+n_2+\ldots+n_N}\epsilon_{n_1,n_2,\ldots,n_N}f_{N}^{(n_1,n_2,\ldots,n_N)}\cdot (iw)^{2N-n_1-n_2-\ldots-n_N}.
$$
Hence
$$
\widehat{y}\left(w\right)=\lim_{N\rightarrow\infty}(iw)^{2N}\sum^{2}_{n_1,n_2,\ldots,n_N=0}\epsilon_{n_1,n_2,\ldots,n_N}f_{N}^{(n_1,n_2,\ldots,n_N)}\cdot (-iw)^{-n_1-n_2-\ldots-n_N}.
$$
$$
\widehat{y}(w)=\lim_{N\rightarrow\infty}\sum^{2N}_{k=0}\left(\sum_{\scriptsize \begin{array}{cc}0\leq n_1,n_2,\ldots,n_{N}\leq 2 \\ n_1+n_2+\ldots+n_N=k\end{array}\normalsize}\epsilon_{n_1,n_2,\ldots,n_N}f_{N}^{(n_1,n_2,\ldots,n_N)}\right)(-iw)^{2N-k}.
$$
\begin{equation}
\widehat{y}(w)=\lim_{N\rightarrow\infty}\sum^{2N}_{k=0}\left(\sum_{\scriptsize \begin{array}{cc}0\leq n_1,n_2,\ldots,n_{N}\leq 2\\ n_1+n_2+\ldots+n_N=2N-k\end{array}\normalsize}\epsilon_{n_1,n_2,\ldots,n_N}f_{N}^{(n_1,n_2,\ldots,n_N)}\right)(-iw)^{k}.
\end{equation}
But
$$
\left(\frac{\partial^{n_1+n_2+\ldots+n_{N}}}{\partial x_1^{n_1}\partial x_2^{n_2}\ldots\partial x_N^{n_N}}f_{N}(P)\right)_{P=(0,0,\ldots,0)}=\frac{1}{(2\pi)^N}\int_{\textbf{\scriptsize R}\normalsize}\int_{\textbf{\scriptsize R\normalsize}}\ldots\int_{\textbf{\scriptsize R\normalsize}}\widehat{f}_{N}(\gamma_1,\gamma_2,\ldots,\gamma_{N})\times
$$
\begin{equation}
\times (i\gamma_1)^{n_1}(i\gamma_2)^{n_2}\ldots (i\gamma_N)^{n_N}d\gamma_1 d\gamma_2\ldots d\gamma_N.
\end{equation}
Hence
$$
\widehat{y}(w)=\lim_{N\rightarrow\infty}\sum^{2N}_{k=0}\left(\sum_{\scriptsize\begin{array}{cc}0\leq n_1,n_2,\ldots,n_N\leq 2\\n_1+n_2+\ldots+n_N=2N-k\end{array}\normalsize}\epsilon_{n_1,n_2,\ldots,n_N}f_{N}^{(n_1,n_2,\ldots,n_N)}\right)(-iw)^k=
$$
$$
=\lim_{N\rightarrow\infty}\sum^{2N}_{k=0}\left(\sum_{\scriptsize\begin{array}{cc}0\leq n_1,n_2,\ldots,n_N\leq 2\\n_1+n_2+\ldots+n_N=k\end{array}\normalsize}\epsilon_{n_1,n_2,\ldots,n_N}f_{N}^{(2-n_1,2-n_2,\ldots,2-n_N)}\right)(-iw)^k=
$$
$$
=\lim_{N\rightarrow \infty}\frac{1}{(2\pi)^N}\int_{R^N}\widehat{f}_{N}(\gamma_1,\gamma_2,\ldots,\gamma_N)\times
$$
$$
\times\sum^{2N}_{k=0}\left(\sum_{\scriptsize\begin{array}{cc}0\leq n_1,n_2,\ldots,n_N\leq 2\\n_1+n_2+\ldots+n_N=2N-k\end{array}\normalsize}\epsilon_{n_1,n_2,\ldots,n_N} (i\gamma_1)^{n_1}(i\gamma_2)^{n_2}\ldots (i\gamma_N)^{n_N}\right)(-iw)^{2N-k}\times 
$$
$$
\times d\gamma_1 d\gamma_2\ldots d\gamma_N.
$$
Hence we get the next\\
\\
\textbf{Theorem 6.}\\
Assume the equation
\begin{equation}
y''(x)=V(x)y(x).
\end{equation}
Then
\begin{equation}
\widehat{y}(w)=\lim_{N\rightarrow\infty}\left(\frac{\partial^2}{\partial h_1^2}\frac{\partial^2}{\partial  h_2^2}\ldots\frac{\partial^2}{\partial h_N^2}\frac{e^{-ih_1w-ih_2 w-\ldots-ih_N w}}{V(h_1)V(h_1+h_2)\ldots V(h_1+h_2+\ldots+h_N)}\right)_{h_1=h_2=\ldots=h_N=0}.
\end{equation}
Also if we set
\begin{equation}
\widehat{y}_{N}(w):=\frac{1}{(2 \pi)^N}\int_{R^N}\widehat{f}_{N}(\gamma_1,\gamma_2,\ldots,\gamma_N)(\gamma_1-w)^2(\gamma_2-w)^2\ldots (\gamma_N-w)^2d\gamma_1d\gamma_2\ldots d\gamma_N ,
\end{equation}
then
\begin{equation}
\widehat{y}(w)=\lim_{N\rightarrow\infty}\widehat{y}_{N}(w)
\end{equation}
and
$$
\widehat{y}(w)=\lim_{N\rightarrow\infty}\sum^{2N}_{k=0}\left(\sum_{\scriptsize\begin{array}{cc}0\leq n_1,n_2,\ldots,n_N\leq 2\\n_1+n_2+\ldots+n_N=k\end{array}\normalsize}\epsilon_{n_1,n_2,\ldots,n_N}f_{N}^{(2-n_1,2-n_2,\ldots,2-n_N)}\right)(-iw)^k.\eqno{(47.1)}
$$
However
$$
\frac{d^{2N}}{dw^{2N}}\widehat{y}_N(w)=\textrm{costant}.
$$
Hence exists costants $C_{N,k}$, such that
$$
\widehat{y}_N(w)=\sum^{2N}_{k=0}C_{N,k}w^k.
$$
\\
 
Set
$$
Y(x):=\frac{1}{V(x)}.
$$
We have
$$\small
\widehat{f}_N(\gamma_1,\gamma_2,\ldots,\gamma_N)=\int_{R^N}Y(t_1)Y(t_1+t_2)\ldots Y(t_1+t_2+\ldots+t_N)e^{-it_1\gamma_1-it_2\gamma_2-it_N\gamma_N}dt_1dt_2\ldots dt_N.\eqno{(47.2)}\normalsize
$$
Also
$$
\widehat{Y}(x)=\int_{R}Y(t)e^{-it x}dt.
$$
Then from (46) we get
$$
\widehat{y}_N(w)=\frac{1}{(2\pi)^N}\int_{\textbf{\scriptsize R}^N}\widehat{Y}(t_1-t_2)\widehat{Y}(t_2-t_3)\ldots \widehat{Y}(t_{N-1}-t_{N}) \widehat{Y}(t_N+w)t_1^2t_2^2\ldots t_N^2dt_1dt_2\ldots dt_N.
$$
Hence
$$
\widehat{y}_N(w)=\frac{1}{(2\pi)^N}\int_{\textbf{\scriptsize R}^N}\widehat{Y}(t_1+w)\widehat{Y}(t_2-t_1)\widehat{Y}(t_3-t_2)\ldots\widehat{Y}(t_{N-1}-t_{N-2})\widehat{Y}(t_{N}-t_{N-1})\times
$$
$$
\times t_1^2t_2^2\ldots t_N^2dt_1dt_2\ldots dt_N=
$$
$$
=\frac{1}{(2\pi)^N}\int_{\textbf{\scriptsize R}^N}\widehat{Y}(t_1+w)\widehat{Y}(t_2-t_1)\ldots\widehat{Y}(t_{N-1}-t_{N-2})\widehat{Y}(t_{N}-t_{N-1}) t_1^2t_2^2\ldots t_N^2dt_1dt_2\ldots dt_N.
$$
Hence
$$
\frac{1}{2\pi}\int_{\textbf{\scriptsize R}}\widehat{y}_N(w)e^{iww'}dw=\frac{1}{(2\pi)^N}\int_{\textbf{\scriptsize R}^N} \widehat{Y}(t_1-t_2)\widehat{Y}(t_2-t_3)\ldots\widehat{Y}(t_{N-1}-t_{N})\times 
$$
$$
\times \left(\frac{1}{2\pi}\int_{\textbf{\scriptsize R}}\widehat{Y}(t_N+w)e^{iww'}dw\right) t_1^2t_2^2\ldots t_N^2 dt_1dt_2\ldots dt_{N}.
$$
Hence
$$
y_N(w')=\frac{1}{(2\pi)^N}\int_{\textbf{\scriptsize R}^N} \widehat{Y}(t_1-t_2)\widehat{Y}(t_2-t_3)\ldots\widehat{Y}(t_{N-1}-t_{N})\times 
$$
$$
\times \left(\frac{1}{2\pi}\int_{\textbf{\scriptsize R}}\widehat{Y}(w)e^{i(w-t_N)w'}dw\right) t_1^2t_2^2\ldots t_N^2 dt_1dt_2\ldots dt_{N}.
$$
Hence
$$
y_N(w')=\frac{1}{(2\pi)^N}\int_{\textbf{\scriptsize R}^{N}} \widehat{Y}(t_1-t_2)\widehat{Y}(t_2-t_3)\ldots\widehat{Y}(t_{N-1}-t_{N})\times
$$
$$
\times Y(w') e^{-it_N w'} t_1^2t_2^2\ldots t_N^2 dt_1dt_2\ldots dt_{N}\Rightarrow
$$
$$
y_N(w')=Y(w')\frac{1}{(2\pi)^N}\int_{\textbf{\scriptsize R}^{N}} \widehat{Y}(t_1-t_2)\widehat{Y}(t_2-t_3)\ldots\widehat{Y}(t_{N-1}-t_{N})e^{-it_N w'} t_1^2t_2^2\ldots t_N^2 dt_1dt_2\ldots dt_{N}.
$$
Hence
$$
\frac{y_N(w)}{Y(w)}=\frac{1}{(2\pi)^N}\int_{\textbf{\scriptsize R}^{N}} \widehat{Y}(t_1-t_2)\widehat{Y}(t_2-t_3)\ldots\widehat{Y}(t_{N-1}-t_{N})e^{-it_N w} t_1^2t_2^2\ldots t_N^2 dt_1dt_2\ldots dt_{N}.
$$
Hence
$$
\frac{y_N(w)}{Y(w)}=\frac{1}{(2\pi)^N}\int_{\textbf{\scriptsize R}^N}\widehat{Y}(t_2-t_1)\widehat{Y}(t_3-t_2)\ldots\widehat{Y}(t_{N}-t_{N-1})e^{-it_1w}t_1^2t_2^2\ldots t_N^2dt_1dt_2\ldots dt_N=
$$
$$
=\frac{1}{(2\pi)^N}\int_{\textbf{\scriptsize R}}\int_{\textbf{\scriptsize R}}\ldots \int_{\textbf{\scriptsize R}}e^{i(t_2-t_1)w}\widehat{Y}(t_2-t_1)t_1^2dt_1e^{i(t_3-t_2)w}\widehat{Y}(t_3-t_2)t_2^2dt_2\ldots\times
$$
$$
\times e^{i(t_N-t_{N-1})w}\widehat{Y}(t_N-t_{N-1})t_{N-1}^2dt_{N-1}e^{-it_Nw}t_N^2dt_{N}.
$$
Hence if we define the operator $F_{n}(x,w)$ such that
$$
F_{n+1}(x,w)=\frac{1}{2\pi}\int_{\textbf{\scriptsize R}}F_n(t,w)e^{(x-t)iw}\widehat{Y}(x-t)t^2dt,\eqno{(47.2)}
$$
then
$$
\widehat{F}_{N+1}(\gamma,w)=-\frac{1}{2\pi}e^{i\gamma w}\widehat{Y}(\gamma)\frac{d^2}{d\gamma^2}\widehat{F}_{N}(\gamma,w)
$$
and
$$
\frac{y_N(w)}{Y(w)}=\left[\int_{\textbf{\scriptsize R}}F_N(t,w)e^{-it\gamma}t^2dt\right]_{\gamma=w}=-\left[\frac{d^2}{d\gamma^2}\widehat{F}_N(\gamma,w)\right]_{\gamma=w}.
$$
Hence
$$
\frac{y_{\infty}(w)}{Y(w)}=\frac{y(w)}{Y(w)}=-\frac{d^2}{d\gamma^2}\left[-\frac{1}{2\pi}e^{i\gamma w}\widehat{Y}(\gamma)\frac{d^2}{d\gamma^2}\left[-\frac{1}{2\pi}e^{i\gamma w}\widehat{Y}(\gamma)\frac{d^2}{d\gamma^2}\left[\ldots\right]\right]\right]_{\gamma=w}
$$
and we get the next\\
\\
\textbf{Theorem 7.}\\
Assume the equation
\begin{equation}
y''(x)=V(x)y(x),
\end{equation}
where the potential $V(x)$ is of the form
\begin{equation}
V(x)=\sum^{\infty}_{n=-\infty}c_ne^{-inx}.
\end{equation}
If 
\begin{equation}
Y(x):=\frac{1}{V(x)},
\end{equation}
then it holds
\begin{equation}
y(x)=-Y(x)\frac{d^2}{d\gamma^2}\left[-\frac{1}{2\pi}e^{i\gamma x}\widehat{Y}(\gamma)\frac{d^2}{d\gamma^2}\left[-\frac{1}{2\pi}e^{i\gamma x}\widehat{Y}(\gamma)\frac{d^2}{d\gamma^2}\left[\ldots\right]\right]\right]_{\gamma=x}.
\end{equation}
Also and obviously one can see that immediately from (48):
\begin{equation}
y(x)=\frac{1}{V(x)}\frac{d^2}{dx^2}\left[\frac{1}{V(x)}\frac{d^2}{dx^2}\left[\frac{1}{V(x)}\frac{d^2}{dx^2}\left[\frac{1}{V(x)}\ldots\right]\right]\right].
\end{equation}
\\

It is of very interest to find the Fourier transform $\widehat{f}_N(\gamma_1,\gamma_2,\ldots,\gamma_N)$. From relation (47.2) and discusion below of it, we find (setting $t_1=t_N'$, $t_2=t_{N-1}'$, $t_k=t_{N-k+1}'$ and then $t_1'=h_1-h_2$, $t_2'=h_2-h_3$, $\ldots$, $t_{N-2}'=h_{N-2}-h_{N-1}$, $t_{N-1}'=h_{N-1}-h_{N}$, $t_N'=h_N$).
$$
\widehat{f}_N(\gamma_1,\gamma_2,\ldots,\gamma_N)=\int_{\textbf{\scriptsize  R\normalsize}^N}Y(t_N')Y(t_N'+t_{N-1}')\ldots\times
$$
$$
\times\ldots Y(t_N'+t_{N-1}'+\ldots+t_{2}')Y(t_N'+t_{N-1}'+\ldots+t_2'+t_1')\times\ldots
$$
$$
\times\exp{\left(-it_N'\gamma_1-it_{N-1}'\gamma_2-\ldots-it_2'\gamma_{N-1}-it_1'\gamma_N\right)}dt_1'dt_2'\ldots dt_N'=
$$
$$
=\int_{\textbf{\scriptsize R}^N}Y(h_N)Y(h_{N-1})\ldots Y(h_2) Y(h_1)\times
$$
$$
\times\exp{\left(-i\gamma_1 h_N-i(h_{N-1}-h_N)\gamma_2-\ldots-i(h_2-h_1)\gamma_{N-1}-i(h_1-h_2)\gamma_N\right)}dh_1dh_2\ldots dh_N=
$$
$$
=\int_{\textbf{\scriptsize R}^N}Y(h_N)Y(h_{N-1})\ldots Y(h_2) Y(h_1)\times
$$
$$
\times\exp[-i h_N(\gamma_1-\gamma_2)-ih_{N-1}(\gamma_2-\gamma_3)-\ldots\times
$$
$$
\times\ldots-ih_3(\gamma_{N-2}-\gamma_{N-1})-ih_2(\gamma_{N-1}-\gamma_{N})-ih_1\gamma_N]dh_1dh_2\ldots dh_N=
$$
$$
=\widehat{Y}(\gamma_1-\gamma_2)\widehat{Y}(\gamma_2-\gamma_3)\ldots \widehat{Y}(\gamma_{N-2}-\gamma_{N-1})\widehat{Y}(\gamma_{N-1}-\gamma_{N})\widehat{Y}(\gamma_N).
$$
\\
Hence we get the next\\
\\
\textbf{Theorem 8.}\\
Assume
$$
f_N(x_1,x_2,\ldots,x_N)=Y(x_1)Y(x_1+x_2)\ldots Y(x_1+x_2+\ldots+x_N).\eqno{(52.1)}
$$
Then
$$
\widehat{f}_N(\gamma_1,\gamma_2,\ldots,\gamma_N)=\widehat{Y}(\gamma_1-\gamma_2)\widehat{Y}(\gamma_2-\gamma_3)\ldots \widehat{Y}(\gamma_{N-2}-\gamma_{N-1})\widehat{Y}(\gamma_{N-1}-\gamma_{N})\widehat{Y}(\gamma_N),\eqno{(52.2)}
$$
where
$$
Y(x)=\frac{1}{V(x)}.
$$
Hence
$$
\widehat{y}_{N}(w):=\frac{1}{(2 \pi)^N}\int_{R^N}\widehat{Y}(\gamma_1-\gamma_2)\widehat{Y}(\gamma_2-\gamma_3)\ldots \widehat{Y}(\gamma_{N-2}-\gamma_{N-1})\widehat{Y}(\gamma_{N-1}-\gamma_{N})\widehat{Y}(\gamma_N)\times
$$
$$
\times(\gamma_1-w)^2(\gamma_2-w)^2\ldots (\gamma_N-w)^2d\gamma_1d\gamma_2\ldots d\gamma_N
$$
and $\lim_{N\rightarrow\infty}\widehat{y}_{N}(w)=\widehat{y}(w)$.\\
\\
\textbf{Theorem 8.1}\\
Assume the real $2\pi-$periodic potential\\
$$
V(x)=\sum^{\infty}_{n=-\infty}c_ne^{-in x}\textrm{ and }\frac{1}{V(x)}=\sum^{\infty}_{n=-\infty}c^{*}_ne^{-in x}\textrm{, }x\in\textbf{R}. 
$$
If the equation is
$$
\sum^{\infty}_{n=0}g_n(-i)^n\frac{d^n}{dx^n}y(x)=V(x)y(x)\eqno{(eq)}
$$
and $y_l(x)$ is a solution of $(eq)$ and $\widehat{y}_l(z)=\theta_l(z)$, then equations $(eq1), (eq2)$ become equivalent with one equation (see Theorem 5.1):
$$
\theta_l(z)=\sum^{\infty}_{n=-\infty}c^{*}_ng(z+n)\theta_l(z+n).\eqno{(eq0)}
$$
\\
\textbf{Remarks.}\\
If we demand $\theta_l(z)$ to be $1-$periodic, we have
$$
\sum^{\infty}_{n=-\infty}c^{*}_ng(z+n)=1.
$$
The oposite states that if $g(z)=\sum_{n=0}^{\infty}g_nz^n$ is such that
$$
\sum^{\infty}_{n=-\infty}c^{*}_ng(z+n)=1,
$$
then $\theta_l(z)=\widehat{y}_l(z)$ and $y_l(x)$ satisfies $(eq)$, if $\theta_l(z)$ satisfies $(eq0)$. But if $\theta_l(z)$ is $1-$periodic then this is the case.\\
\\

Assuming the notation of Theorem 8.1 we have: If 
$$
\theta(z+1)=H(z)\theta(z),
$$
then
$$
\theta(z+n)=\left(\prod_{j=0}^{n-1}H(z+j)\right)\theta(z).
$$
Hence
$$
\theta(z)=\sum^{\infty}_{n=0}c^{*}_ng(z+n)\theta(z+n)\Leftrightarrow \theta(z)=\sum^{\infty}_{n=0}c^{*}_ng(z+n)\left(\prod^{n-1}_{j=0}H(z+j)\right)\theta(z).
$$
Hence
$$
1=\sum_{n=0}^{\infty}c^{*}_nB^{*}_n(z),
$$
where 
$$
B^{*}_n(z)=g(z+n)\prod^{n-1}_{j=0}H(z+j).
$$
Hence
$$
g(z+n)=\frac{B^{*}_n(z)}{\prod^{n-1}_{j=0}H(z+j)}\Leftrightarrow
$$
$$
B^{*}_{n+1}(z)=g(z+n+1)\prod^{n}_{j=0}H(z+j)=g(n+z+1)H(z+n)\prod^{n-1}_{j=0}H(z+j)=
$$
$$
=g(z+n+1)H(z+n)\frac{B^{*}_n(z)}{g(z+n)}.
$$
Hence
$$
B^{*}_{n+1}(z)=\frac{g(z+n+1)}{g(z+n)}H(z+n)B^{*}_n(z).
$$
Also is very easy to show that\\
\\
\textbf{Lemma.}
$$
B^{*}_{n+1}(z)=B^{*}_n(z+1)H(z),
$$ 
is equivalent to
$$
B^{*}_n(z)=g(z+n)\prod^{n-1}_{j=0}H(z+j),
$$
where $B^{*}_0(z)=g(z)$ and $\theta(z+1)=H(z)\theta(z)$.\\
\\
\textbf{Proof.}\\
We have
$$
B^{*}_{n+1}(z)\theta(z)=B^{*}_n(z+1)\theta(z+1).
$$
Hence if we set $\xi_n(z):=B^{*}_n(z)\theta(z)$, then
$$
\xi_{n+1}(z)=\xi_{n}(z+1).
$$
Hence 
$$
\xi_n(z)=\xi_{n-1}(z+1)\Leftrightarrow\xi_n(z)=\xi_0(z+n)\Leftrightarrow 
$$
$$
B^{*}_n(z)\theta(z)=B^{*}_0(z+n)\theta(z+n).
$$ 
Hence $B^{*}_n(z)=g(z+n)\prod^{n-1}_{j=0}H(z+j)$.\\
\\

If the function $\theta(z)$ is defined from the equation
$$
\theta(z+1)=H(z)\theta(z),
$$
then $\theta(z)$ satisfies 
$$
\sum^{\infty}_{n=0}c_n\theta(z+n)=g(z)\theta(z),
$$
iff
$$
\theta(z)=\sum^{\infty}_{n=0}c^{*}_ng(z+n)\theta(z+n),
$$
iff
$$
\sum^{\infty}_{n=0}B_n(z)=1,\eqno{(a)}
$$
where 
$$
B_n(z)=c_n^{*}g(z+n)\prod^{n-1}_{j=0}H(z+j).
$$
But this last equation is equivalent to:
$$
B_{n+1}(z)=\frac{c^{*}_{n+1}}{c^{*}_{n}}H(z)B_n(z+1).\eqno{(b)}
$$
Hence from $(a),(b)$ we get
$$
H(z)=\sum^{\infty}_{n=0}\frac{c^{*}_n}{c^{*}_{n+1}}B_{n+1}(z).\eqno{(c)}
$$
However we need to show 
$$
\frac{c^{*}_n}{c^{*}_{n+1}}\frac{B_{n+1}(z)}{B_{n}(z+1)}=H(z)=\textrm{independent of }n.
$$
For this, always exist function $\lambda_n(z)$ such that
$$
\frac{B_{n+1}(z)}{B_n(z)}=\lambda_n(z)
$$
and also
$$
\frac{c_n^{*}}{c^{*}_{n+1}}\frac{B_{n+1}(z)}{B_{n}(z+1)}=\frac{c_{n+1}^{*}}{c_{n+2}^{*}}\frac{B_{n+2}(z)}{B_{n+1}(z+1)}\Leftrightarrow\frac{c^{*}_{n+1}B_{n+2}(z)}{c_n^{*}B_{n+1}(z)}=\frac{c^{*}_{n+2}B_{n+1}(z+1)}{c^{*}_{n+1}B_n(z+1)}\Leftrightarrow
$$
$$
\frac{\frac{1}{c^{*}_{n+2}}B_{n+2}(z)}{\frac{1}{c^{*}_{n+1}}B_{n+1}(z)}=\frac{\frac{1}{c^{*}_{n+1}}B_{n+1}(z+1)}{\frac{1}{c^{*}_{n}}B_n(z+1)}.
$$
Hence if we set 
$$
B^{*}_n(z)=\frac{1}{c^{*}_n}B_{n}(z)
$$
and
$$
\lambda^{*}_n(z)=\frac{B^{*}_{n+1}(z)}{B^{*}_n(z)},
$$
then it must be
$$
\lambda^{*}_{n+1}(z)=\lambda^{*}_{n}(z+1)\Leftrightarrow \lambda^{*}_{n}(z)=\lambda^{*}_{n-1}(z+1).
$$
Hence iterating we get equivalent
$$
\lambda^{*}_{n}(z)=\lambda^{*}_{0}(z+n)
$$
and
$$
\frac{B^{*}_{n+1}(z)}{B^{*}_n(z)}=\lambda^{*}_0(z+n)\Leftrightarrow B^{*}_n(z)=\left(\prod^{n-1}_{j=0}\lambda^{*}_0(z+j)\right)B^{*}_0(z).
$$
Hence equivalent
$$
H(z)=\frac{B^{*}_{n+1}(z)}{B^{*}_{n}(z+1)}=\frac{B^{*}_0(z)\prod^{n}_{j=0}\lambda^{*}
_0(z+j)}{B^{*}_0(z+1)\prod^{n-1}_{j=0}\lambda^{*}_0(z+1+j)}=\frac{B^{*}_0(z)}{B^{*}_0(z+1)}\lambda^{*}_0(z).
$$
Hence
$$
\lambda^{*}_0(z)=\frac{B^{*}_0(z+1)}{B^{*}_0(z)}H(z).
$$
Hence
$$
\frac{B^{*}_{n}(z)}{B^{*}_0(z)}=\prod^{n-1}_{j=0}\lambda^{*}_0(z+j)=
$$
$$
=\frac{B^{*}_0(z+1)}{B^{*}_0(z)}H(z)\frac{B^{*}_0(z+2)}{B^{*}_0(z+1)}H(z+1)\frac{B^{*}_0(z+3)}{B^{*}_0(z+2)}H(z+2)\times\ldots
$$
$$
\ldots\times\frac{B^{*}_0(z+n)}{B^{*}_0(z+n-1)}H(z+n-1)=\frac{B^{*}_0(z+n)}{B^{*}_0(z)g(z+n)}g(z+n)\prod^{n-1}_{j=0}H(z+j)=
$$
$$
=\frac{B^{*}_0(z+n)}{B^{*}_0(z)g(z+n)}B^{*}_n(z).
$$
Hence it must be 
$$
B^{*}_{0}(z+n)=g(z+n)\textrm{, }\forall z\in \textbf{R},\forall n\in\textbf{N}\Leftrightarrow B^{*}_0(z)=g(z).
$$

For example for any functions $f_1(z),f_2(z)$ we have (this give us $(a)$):
$$
B_n(z)=\frac{1}{n!}\sum^{n}_{k=0}C_{n,k}\frac{f_{1}^{(k)}(z)f_{2}^{(n-k)}(z)}{f_1(z+1)f_2(z+1)},
$$
where
$$
C_{n,k}:=\left(\begin{array}{cc}n\\k\end{array}\right),
$$
is the binomial symbol. 

Moreover equation $(b)$ is equivalent to equation
$$
\frac{B^{*}_{n+1}(z)}{B^{*}_n(z)}=\lambda^{*}_0(z+n),
$$
where
$$
H(z)=\frac{g(z)}{g(z+1)}\lambda^{*}_0(z)
$$ 
and from $(c)$:
$$
H(z)=\frac{g(z)}{g(z+1)}\lambda^{*}_0(z)=\sum^{\infty}_{n=0}c^{*}_nB^{*}_{n+1}(z)=\sum^{\infty}_{n=0}\frac{c^{*}_n}{c^{*}_{n+1}}B_{n+1}(z).
$$
$$
\frac{B^{*}_{n+1}(z)}{B^{*}_n(z)}=\lambda^{*}_n(z)=\lambda_0^{*}(z+n)=\frac{c_{n}^{*}}{c_{n+1}^{*}}\lambda_n(z).
$$
$$
\frac{B_{n+1}(z)}{B_n(z)}=\lambda_n(z)=\frac{c_0^{*}}{c_1^{*}}\frac{c^{*}_{n+1}}{c^{*}_{n}}\lambda_0(z+n).
$$
$$
\lambda_0^{*}(z)=\frac{c_0^{*}}{c_1^{*}}\lambda_0(z).
$$
$$
B_n(z)=\left(\frac{c_0^{*}}{c_1^{*}}\right)^nc_n^{*}g(z)\prod^{n-1}_{j=0}\lambda_0(z+j)
$$
and
$$
B_n(z)=c_n^{*}g(z+n)\prod^{n-1}_{j=0}H(z+j).
$$
$$
\lambda_0^{*}(z)=\frac{g(z+1)}{g(z)}H(z).
$$
\\

Hence we have the next theorem:\\
\\
\textbf{Theorem 8.2}\\
Given the equation
$$
\sum^{\infty}_{n=0}g_n(-i)^n\frac{d^n}{dx^n}y(x)=V(x)y(x),\eqno{(eq)}
$$
with potential $V(x)=\sum_{n}c_ne^{-in x}$ a $2\pi-$periodic function and $1/V(x)=\sum_{n}c_n^{*}e^{-i n x}$ and $g(x)=\sum^{\infty}_{n=0}g_nx^n$, we assume that exists a family of functions $B_n(z)$ such that:\\
\textbf{I.}
$$
\frac{B_{n+1}(z)}{B_n(z)}=\frac{c^{*}_0}{c^{*}_1}\frac{c^{*}_{n+1}}{c^{*}_n}\lambda_0(z+n),\eqno{(eq1)}
$$
where $\lambda_0(z)$ is a simple function. And\\
\textbf{II.}
$$
\sum^{\infty}_{n=0}B_n(z)=1\textrm{, }\forall z\in D,\eqno{(eq2)}
$$
$$
B_0(z)=g(z).\eqno{(eq3)}
$$
Then
$$
\frac{c_0^{*}}{c_1^{*}}\lambda_0(z)=\frac{g(z+1)}{g(z)}H(z).\eqno{(eq4)}
$$
Hence we find $H(z)$ and determine $\theta(z)$, since
$$
\frac{\theta(z+1)}{\theta(z)}=H(z).\eqno{(eq5)}
$$
Hence we get $\widehat{y}(z)=\theta(z)$. Hence we solve $(eq)$.\\
\\

Assume the functions $B^{*}_n(\tau),C_n(\tau),H(\tau)$ are such that 
$$
B^{*}_n(\tau)=\frac{\xi_n(\tau)}{\theta(\tau)},
$$
where $\xi_n(\tau+1)=\xi_{n+1}(\tau)$. Then
$$
B^{*}_{n+1}(\tau)=B^{*}_{n}(\tau+1)H(\tau)\textrm{, }n\in Z\textrm{, }Im(\tau)>0,
$$
$$
B^{*}_{-n}(\tau)=B^{*}_n(\tau)\textrm{, }n\in Z\textrm{, }Im(\tau)>0
$$
and $f(z)=g(z^2)$ and 
$$
\phi(\tau,z)=\sum^{\infty}_{n=-n_0}C_n(\tau)e^{2\pi i n z}=\sum^{\infty}_{n=-\infty}B^{*}_n(\tau)f\left(z^n\right).
$$
Then
$$
\phi(\tau,-1/z)=\sum^{\infty}_{n=-n_0}C_n(\tau)e^{2\pi i (-1/z)n}=\sum^{\infty}_{n=-\infty}B^{*}_n(\tau)f\left((-1)^nz^{-n}\right)=
$$
$$
=\sum^{\infty}_{n=-\infty}B^{*}_{n}(\tau)g\left(\left((-1)^nz^{-n}\right)^2\right)=\sum^{\infty}_{n=-\infty}B^{*}_n(\tau)f\left(z^{-n}\right)=
$$
$$
=\sum^{\infty}_{n=-\infty}B^{*}_{-n}(\tau)f(z^n)=\sum^{\infty}_{n=-\infty}B^{*}_{n}(\tau)f(z^n)=\phi(\tau,z).
$$
Since $\xi_n(z)=\xi_0(z+n)$ and $\xi_{-n}(z)=\xi_0(z-n)$ and $\xi_n(z+1)=\xi_{n+1}(z)$ and
$$
B^{*}_n(z)=\frac{\xi_n(z)}{\theta(z)}=\frac{\xi_0(z+n)}{\theta(z)},
$$
$$
B^{*}_{-n}(z)=\frac{\xi_0(z-n)}{\theta(z)}\textrm{, }B^{*}_n(z)=B^{*}_{-n}(z),
$$
when $\xi_0(z)$ is any $2-$periodic function. i.e. 
$$
\xi_0(z+2)=\xi_0(z)\Leftrightarrow \xi_0(z+2n)=\xi_0(z)\Leftrightarrow 
$$
$$
\xi_0(z+n)=\xi_0(z-n)\textrm{, }\forall z\in R\textrm{ and }\forall n\in Z.
$$
Also
$$
\phi(\tau,z+1)=\sum^{\infty}_{n=-n_0}C_n(\tau)e^{2\pi i n (z+1)}=\sum^{\infty}_{n=-\infty}B^{*}_n(\tau)f\left((z+1)^n\right)=\phi(\tau,z).
$$
\\
\textbf{Theorem 8.3}\\
Set suitable function $f(z)$, such that exists $g(z)$, with $f(z)=g(z^2)$. Also there exists $C_n(\tau)$, such that
$$
\phi(\tau,z):=\sum^{\infty}_{n=-n_0}C_n(\tau)e^{2\pi i n z}=\sum^{\infty}_{n=-\infty}B^{*}_n\left(\tau\right)f\left(z^n\right)\textrm{, }Im(\tau)>0\textrm{, }Im(z)>0 
$$
and
$$
B^{*}_n(\tau)=\frac{\xi_n(\tau)}{\theta(\tau)},
$$
with
$$
\xi_{n+1}(\tau)=\xi_{n}(\tau+1)\textrm{ and }\frac{\theta(\tau+1)}{\theta(\tau)}=H(\tau).
$$
Then we have
$$
\xi_n(\tau)=\xi_0(\tau+n)
$$
and 
$$
B^{*}_n(\tau)=\frac{\xi_0(\tau+n)}{\theta(\tau)}
$$
$$
B^{*}_{-n}(\tau)=\frac{\xi_0(\tau-n)}{\theta(\tau)}
$$
and
$$
B^{*}_{n+1}(\tau)=B^{*}_{n}(\tau+1)H(\tau)\textrm{, }n\in Z\textrm{, }Im(\tau)>0,
$$
$$
B^{*}_{-n}(\tau)=B^{*}_n(\tau)\textrm{, }n\in Z\textrm{, }Im(\tau)>0.
$$
Then 
$$
\phi(\tau,z+1)=\phi(\tau,z)\textrm{ and }\phi(\tau,-1/z)=\phi(\tau,z).
$$

\section{Some interesting solutions}

We know that equation
$$
\sum^{\infty}_{n=0}g_n(-i)^n\frac{d^n}{dx^n}y(x)=V(x) y(x),
$$ 
can be written in the form
$$
\sum^{\infty}_{n=0}c_n\widehat{y}(z+n)=g(z) \widehat{y}(z),
$$
where $g(x)=\sum^{\infty}_{n=0}g_nx^n$, when
$$
V(x)=\sum^{\infty}_{n=0}c_ne^{-in x}.
$$
Also if 
$$
Y(x)=\frac{1}{V(x)}=\sum^{\infty}_{n=-\infty}c^{*}_ne^{-in x},
$$
then
$$
\widehat{y}(z)=\lim_{N\rightarrow\infty}\sum_{n_1\geq 0}\sum_{n_2\geq 0}\sum_{n_3\geq 0}\ldots\times 
$$
$$
\times\ldots\sum_{n_N\geq 0}A_{n_1}(z)A_{n_2}(z+n_1)A_{n_3}(z+n_1+n_2)\ldots A_{n_N}(z+n_1+n_2+\ldots+n_{N}),\eqno{(52.3)}
$$
where
$$
A_n(x)=g(n+x)c_n^{*}.
$$
Also
$$
g(z)=\sum_{k}c_k\prod^{k-1}_{j=0}f_k(z+j)
$$
Hence $\theta(z)$ is solution of (30.1) and consequently
$$
\widehat{y}(z)=\theta(z),
$$
where
$$
\frac{\theta(z+k)}{\theta(z)}=\prod^{k-1}_{j=0}f_k(z+j).
$$
Hence if $f(z)=z$, then we get $\theta(z)=\Gamma(z)$. If $f(z)=z^2$, then $\theta(z)=\Gamma(z)^2$. If $f(z)=\frac{z^2}{z-1}$, then $\theta(z)=(z-1)\Gamma(z)$ . If $f(z)=2z+1$, then $\theta(z)=2^z\Gamma(z+\frac{1}{2})$. If $f(z)=2z^2+z$, then $\theta(z)=2^z\Gamma(z+\frac{1}{2})\Gamma(z)$. If $f(z)=2z^2+z+1$, then $\theta(z)=2^z\Gamma(z-\rho_1)\Gamma(z-\rho_2)$, where $\rho_1=\frac{-1-i\sqrt{7}}{4}$, $\rho_2=\frac{-1+i\sqrt{7}}{4}$. By this way we get the next\\
\\
\textbf{Theorem 9.}\\
Suppose $P(z)=A\prod^{\nu}_{l=1}(z-\rho_{l})$, $Q(z)=B\prod^{\mu}_{l=1}(z-\sigma_{l})$ and $R(z)=\sum^{\infty}_{l=0}a_lz^l$, are two polynomials of degree $\nu$, $\mu$ and an entire function respectively. $P,Q$, have no common roots and $H(z)=C z^m e^{-R(z)}P(z)/Q(z)$. Then if 
\begin{equation}
g(z)=\sum^{\infty}_{k=0}c_k\left(\frac{A}{B}\right)^k\prod^{k-1}_{j=0}H(z+j)=\sum^{\infty}_{k=0}\prod^{k-1}_{j=0}\frac{A}{B}c^{'}_jH(z+j),
\end{equation}
where $\prod^{k-1}_{j=0}c_j^{'}=c_k$, the equation
\begin{equation}
\sum^{\infty}_{n=0}g_n(-i)^ny^{(n)}(x)=V(x)y(x)\textrm{, }V(x)=\sum^{\infty}_{n=0}c_ne^{-inx},
\end{equation}
can be set in the form
\begin{equation}
\sum^{\infty}_{n=0}c_n\widehat{y}(z+n)=g(z)\widehat{y}(z)
\end{equation}
and a solution of (55) is
\begin{equation}
\widehat{y}(z)=\left(\frac{A}{B}\right)^z C^z \Gamma(z)^m e^{-R_1(z)} \frac{\prod^{\nu}_{j=1}\Gamma(z-\rho_j)}{\prod^{\mu}_{j=1}\Gamma(z-\sigma_j)}.
\end{equation}
If $R(z)$ is a polynomial of degree $\tau$, then the function $R_1(z)$ is a polynomial of degree $\tau+1$ (see notes below). If $\tau=\infty$, then again we use the notes below.\\ 
Hence the solution of (54) is
$$
y(x)=\frac{1}{2\pi}\int^{\infty}_{-\infty}\left(\frac{A}{B}\right)^z C^z \Gamma(z)^m e^{-R_1(z)} \left(\frac{\prod^{\nu}_{j=1}\Gamma(z-\rho_j)}{\prod^{\mu}_{j=1}\Gamma(z-\sigma_j)}\right)e^{iz x}dz.\eqno{(eq)}
$$
\\
\textbf{Notes.}\\
\textbf{I.} Actualy if $R(z)=z$, then $R_1(z)=\frac{z^2}{2}-\frac{z}{2}$. If $R(z)=z^2$, then $R_1(z)=\frac{z^3}{3}-\frac{z^2}{2}+\frac{z}{6}$. If $R(z)=z^3$, then $R_1(z)=\frac{z^4}{4}-\frac{z^3}{2}+\frac{z^2}{4}$, $\ldots$ etc. In general if $R(z)=z^{\tau}$ and $p_{\tau}(n)=\sum^{n-1}_{j=0}j^{\tau}=\sum^{\tau+1}_{k=0}h_kn^k$, then $R(z)=\sum^{\tau+1}_{k=0}h_kz^k$. Obviously this procedure is linear i.e. if for example $R(z)=z^2-2z$, then $R_1(z)=\frac{z^3}{3}-\frac{z^2}{2}+\frac{z}{6}-2\left(\frac{z^2}{2}-\frac{z}{2}\right)=\frac{z^3}{3}-\frac{3z^2}{2}+\frac{7z}{6}$.\\
\textbf{II.} If $H(z)$ is meromorphic, with $\rho_j$, $j=1,2,\ldots$ roots and $\sigma_j$, $j=1,2,\ldots$ poles (infinite roots and poles) and exists entire function $R(z)$ such that: 
$$
H(z)=C z^m e^{-R(z)}\prod^{\infty}_{j=1}\left(\frac{z-\rho_j}{z-\sigma_j}\right),\eqno{(56.00)}
$$
provided that $\sum^{\infty}_{k=1}|\rho_k|<\infty$, $\sum^{\infty}_{k=1}|\sigma_k|<\infty$ and 
$$
g(z)=\sum_{k}c_k\prod^{k-1}_{j=0}H(z+j).\eqno{(56.01)}
$$
Then a solution of
$$
\sum_{k}g_k(-i)^ky^{(k)}(x)=V(x)y(x)\textrm{, }V(x)=\sum^{\infty}_{n=0}c_ne^{-in x},
$$
is such that
$$
\widehat{y}(z)=C^z \Gamma(z)^m e^{-R_1(z)} \prod^{\infty}_{j=1}\frac{\Gamma(z-\rho_j)}{\Gamma(z-\sigma_j)}.
$$
Also it holds
$$
\frac{\widehat{y}(z+1)}{\widehat{y}(z)}=H(z)\Leftrightarrow \widehat{y}(z)=\prod^{\infty}_{n=0}\frac{1}{H(z+n)}\eqno{(56.02)}
$$
and
$$
y(x)=\frac{1}{2\pi}\int^{\infty}_{-\infty}C^z \Gamma(z)^m e^{-R_1(z)} \left(\prod^{\infty}_{j=1}\frac{\Gamma(z-\rho_j)}{\Gamma(z-\sigma_j)}\right) e^{iz x}dz.
$$
\textbf{III.} Note that if $g(z)$ is $1-$periodic, then (see Theorem 10 below) $H(z)$ is also $1-$periodic. Hence 
$$
g(z)=\sum_{k}c_k(H(z))^k.\eqno{(56.03)}
$$
Set
$$
V_1(z)=\sum_{k}c_kz^k,\eqno{(56.04)}
$$
then
$$
H(z)=V^{(-1)}_1\left(g(z)\right).\eqno{(56.05)}
$$
Hence given the $c_k$ of the potential $V (z)$ and $g_n$ of $g(z)$ we get, if $g(z)$ is $1-$periodic
$$
\frac{\widehat{y}(z+1)}{\widehat{y}(z)}=V_1^{(-1)}(g(z)).\eqno{(56.06)}
$$
For example if $g(z)=\cos(2\pi z)$ and $V(z)=\sum^{\infty}_{n=1}\frac{e^{-i n z}}{n}$, then 
$$
\cos(2\pi z)=\sum^{\infty}_{n=1}\frac{H(z)^n}{n}=-\log(1-H(z)).
$$ 
Hence $H(z)=1-e^{-\cos(2\pi z)}$ and
$$
\frac{\widehat{y}(z+1)}{\widehat{y}(z)}=1-e^{-\cos(2\pi z)}.\eqno{(56.1)}
$$
However we can not solve with iteration this equation because it is $1-$periodic and
$$
\widehat{y}(z)=\prod^{\infty}_{n=0}\frac{1}{1-e^{-\cos(2\pi (z+n))}}=\prod^{\infty}_{n=0}\frac{1}{1-e^{-\cos(2\pi z)}}.
$$

Moreover if $H(z)$ is semi-periodic i.e. exist $\lambda$ such that $\lambda\neq\pm 1$ and $0<|\lambda|\leq1$, with 
$$
H(z+1)=\lambda H(z)\textrm{, }\forall z\in D.
$$
Then
$$
g(z)=\sum_{k}c_kH(z)^k\lambda^{k(k-1)/2}.
$$
Hence
$$
\frac{\widehat{y}(z+1)}{\widehat{y}(z)}=H(z).
$$
\\

Moreover than this, if we assume that $\lambda(z)$ is semi-periodic and 
$$
\lambda(z+1)=\xi\lambda(z)\textrm{, }0<|\xi|\leq1\textrm{, }\xi\neq\pm 1.
$$ 
Then if 
$$
H(z+1)=\lambda(z)H(z), 
$$
we get easily
$$
H(z+n)=\xi^{n(n-1)/2}\lambda(z)^{n}H(z).
$$
Also if 
$$
\widehat{y}(z+1)=H(z)\widehat{y}(z),\eqno{(e1)}
$$
then
$$
\widehat{y}(z+n)=\xi^{n(n-1)(n-2)/6}\left(\lambda(z)\right)^{n(n-1)/2}H(z)^n\widehat{y}(z).\eqno{(e2)}
$$
But as someone can see if we take $\widehat{y}(z)$ as
$$
\widehat{y}(z)=\sum^{\infty}_{k=-\infty}A_k(z)\xi^{-k(k-1)(k-2)/6}\lambda(z)^{-k(k-1)/2}H(z)^{-k},\eqno{(e3)}
$$
then $(e1)$ holds if
$$
A_{k-1}(z+1)=A_k(z).\eqno{(e4)}
$$
Hence $(e2)$ will holds also, since from $(e3)$ we can easily see that
$$
\widehat{y}(z+1)=H(z)\sum^{\infty}_{k=-\infty}A_{k-1}(z+1)\xi^{-k(k-1)(k-2)/6}\lambda(z)^{-k(k-1)/2}H(z)^{-k}\Leftrightarrow
$$
$$
\frac{\widehat{y}(z+1)}{\widehat{y}(z)}=H(z).
$$
Also we can write from $(e2)$:
$$
\sum^{\infty}_{k=-\infty}c_k\widehat{y}(z+k)=g(z)\widehat{y}(z)
$$
iff
$$
g(z)=\sum^{\infty}_{k=-\infty}c_k\xi^{k(k-1)(k-2)/6}\lambda(z)^{k(k-1)/2}H(z)^{k}.\eqno{(e5)}
$$
\\
\textbf{Theorem  9.0}\\
Assume the function $\widehat{y}(z)$. If exists $A_k(z)$ such that
$$
\widehat{y}(z)=\sum_{k=-\infty}^{\infty}A_k(z)\xi^{-k(k-1)(k-2)/6}\lambda(z)^{-k(k-1)/2}H(z)^{-k}
$$
and
$$
A_k(z)=A_{k-1}(z+1),
$$
for these $\lambda(z),H(z)$, we will have
$$
\frac{\widehat{y}(z+1)}{\widehat{y}(z)}=H(z)
$$
and the equation
$$
\sum^{\infty}_{k=-\infty}c_k\widehat{y}(z+k)=g(z)\widehat{y}(z)
$$
will be equivalent to
$$
g(z)=\sum^{\infty}_{k=-\infty}c_k\xi^{k(k-1)(k-2)/6}\lambda(z)^{k(k-1)/2}H(z)^{k}.
$$
\textbf{Remarks.}\\
\textbf{I.} The functions $\lambda(x)$ and $H(x)$, are such that
$$
\lambda(x+1)=\xi\lambda(x)\textrm{, }0<|\xi|\leq 1\textrm{, }\xi\neq \pm 1,
$$
$$
H(x+1)=\lambda(x)H(x)
$$
and
$$
V(x)=\sum^{\infty}_{n=-\infty}c_ne^{-inx}.
$$
(For example a possible construction of $\lambda(x)$ function, is to define a function $L(y,x)$ such that at a specific point $y=\xi$ (or maybe all $y\in D$), we have $L(\xi,x+1)=L(\xi,x)$, $\forall x$ and $\lambda(x)=\xi^{x}L(\xi,x)$, then $\lambda(x+1)=\xi^{x+1}L(\xi,x+1)=$ $\xi\left(\xi^x L(\xi,x)\right)$ $=\xi \lambda(x)$.)\\ 
\textbf{II. } Hence if the function $g(x)=\sum^{\infty}_{n=0}g_nx^n$, which is related with the problem   
$$
\sum^{\infty}_{n=0}g_n(-i)^ny^{(n)}(x)=V(x)y(x),
$$
(where $V(x)$ is a continuous $2\pi-$periodic potential), is of the form 
$$
g(x)=\sum^{\infty}_{k=-\infty}c_k\xi^{k(k-1)(k-2)/6}\lambda(x)^{k(k-1)/2}H(x)^{k},
$$
then we have solution (of the problem) $y(x)$ such that
$$
\frac{\widehat{y}(z+1)}{\widehat{y}(z)}=H(z)
$$
and
$$
\widehat{y}(z)=\sum_{k=-\infty}^{\infty}A_k(z)\xi^{-k(k-1)(k-2)/6}\lambda(z)^{-k(k-1)/2}H(z)^{-k},
$$
where $A_k(z)=A_{k-1}(z+1)$.\\
\textbf{III.} Hence for every given function $g(x)$ and sequence $c_k$ we have an expansion of the form (here for simplicity we take $V(x)$ to be a trigonometric polynomial of degree $N$):
$$
g(x)=\sum^{N}_{k=-N}c_k\xi^{k(k-1)(k-2)/6}\lambda(x)^{k(k-1)/2}H(x)^{k},
$$ 
where $\lambda(x)$, $H(x)$ are defined as in $\textbf{I}$ of these remarks. Hence 
$$
g(x)=\sum^{N}_{k=-N}c_k\xi^{k(k-1)(k-2)/6}\left(\frac{H(x+1)}{H(x)}\right)^{k(k-1)/2}H(x)^{k}.
$$
Hence we have the polynomial equation
$$
S_N\left(\xi,\frac{H(x+1)}{H(x)},H(x)\right)=g(x),\eqno{(e6)}
$$
where
$$
S_N(\xi,x,y)=\sum^{N}_{k=-N}c_k\xi^{k(k-1)(k-2)/6}x^{k(k-1)/2}y^{k}.\eqno{(e7)}
$$
Hence solving (e6), with respect to $\lambda(x)=H(x+1)/H(x)$, we have
$$
H(x+1)=H(x)R_N(\xi,g(x),H(x)).\eqno{(e8)}
$$
Hence
$$
\lambda(x)=R_{N}(\xi,g(x),H(x)),\eqno{(e9)}
$$
where $R_N(\xi,x,y)$ is an algebraic function solution of (e7), with respect to $x$.\\
\textbf{IV. } Hence assuming that $g(x)$ is of the form:
$$
g(x)=\sum^{\infty}_{k=0}b_k\lambda(x)^k,
$$
we have 
$$
g(x)=\sum^{\infty}_{k=-\infty}c_k\xi^{k(k-1)(k-2)/6}\lambda(x)^{k(k-1)/2}H(x)^{k},
$$
$$
\sum^{\infty}_{k=0}b_kx^k=\sum^{\infty}_{k=-\infty}c_k\xi^{k(k-1)(k-2)/6}x^{k(k-1)/2}H\left(\lambda^{(-1)}(x)\right)^{k}\Leftrightarrow
$$
$$
\sum^{\infty}_{k=0}b_kx^k=\sum^{\infty}_{k=-\infty}c_k\xi^{k(k-1)(k-2)/6}x^{k(k-1)/2}f(x)^{k}.
$$
Hence exists function $f(x)$ (by inversion) such that ($f(x)$ is known function)
$$
H\left(\lambda^{(-1)}(x)\right)=f(x).
$$
Hence
$$
x=g^{(-1)}\left(\sum^{\infty}_{k=-\infty}c_k\xi^{k(k-1)(k-2)/6}\lambda(x)^{k(k-1)/2}f(\lambda(x))^k\right).
$$
Hence
$$
\lambda^{(-1)}(x)=g^{(-1)}\left(\sum^{\infty}_{k=-\infty}c_k\xi^{k(k-1)(k-2)/6}x^{k(k-1)/2}f(x)^k\right).
$$
Hence $\lambda(x)$ follows. Also then
$$
H(x)=f\left(\lambda(x)\right),
$$
and $H(x)$ also follows.\\Note that 
$$
\lambda(z+1)=\xi\lambda(z)\Leftrightarrow \lambda^{(-1)}(\xi z)-\lambda^{(-1)}(z)=1\Leftrightarrow
$$
$$
\lambda^{(-1)}\left(\xi^{z+1}\right)-\lambda^{(-1)}\left(\xi^z\right)=1.\eqno{(a)}
$$
From relation $(a)$ we get that: for the function 
$$
f(z)=H\left(\lambda^{(-1)}(z)\right),\eqno{(b)}
$$
holds 
$$
f\left(\xi^{z+1}\right)=H\left(\lambda^{(-1)}\left(\xi^{z+1}\right)\right)=H\left(1+\lambda^{(-1)}\left(\xi^z\right)\right)=
$$
$$
=\lambda\left(\lambda^{(-1)}\left(\xi^z\right)\right)H\left(\lambda^{(-1)}\left(\xi^{z}\right)\right)=\xi^zf\left(\xi^z\right).
$$
Hence
$$
f(\xi z)=zf(z)\eqno{(c)}
$$
and
$$
\lambda^{(-1)}(z)=g^{(-1)}\left(\sum^{\infty}_{k=-\infty}c_k\xi^{k(k-1)(k-2)/6}z^{k(k-1)/2}f(z)^k\right),\eqno{(d)}
$$
when the $g(z)$ and $c_k$ are given.\\
It holds from $(a),(d)$:
$$
g^{(-1)}\left(\sum^{\infty}_{k=-\infty}c_{k}\xi^{k(k-1)(k-2)/6}(\xi z)^{k(k-1)/2}f(\xi z)^k\right)-
$$
$$
-g^{(-1)}\left(\sum^{\infty}_{k=-\infty}c_{k}\xi^{k(k-1)(k-2)/6}z^{k(k-1)/2}f(z)^k\right)=1\eqno{(e)}
$$
Hence we have from $(c)$ and $(d)$:
$$
g^{(-1)}\left(\sum^{\infty}_{k=-\infty}c_{k}\xi^{k(k-1)(k-2)/6}(\xi z)^{k(k-1)/2}\left(zf(z)\right)^k\right)-
$$
$$
-g^{(-1)}\left(\sum^{\infty}_{k=-\infty}c_{k}\xi^{k(k-1)(k-2)/6}z^{k(k-1)/2}f(z)^k\right)=1.\eqno{(d1)}
$$
Note that we have an identity also:
$$
g^{(-1)}\left(f(z)^{-1}\sum^{\infty}_{k=-\infty}c_{k-1}\xi^{k(k-1)(k-2)/6}z^{k(k-1)/2}f(z)^k\right)=
$$
$$
=g^{(-1)}\left(\sum^{\infty}_{k=-\infty}c_k\xi^{k(k-1)(k-2)/6}\left(\xi z\right)^{k(k-1)/2}\left(zf(z)\right)^k\right).\eqno{(d2)}
$$
Hence
$$
g^{(-1)}\left(f(z)^{-1}\sum^{\infty}_{k=-\infty}c_{k-1}\xi^{k(k-1)(k-2)/6}z^{k(k-1)/2}f(z)^k\right)-
$$
$$
-g^{(-1)}\left(\sum^{\infty}_{k=-\infty}c_{k}\xi^{k(k-1)(k-2)/6}z^{k(k-1)/2}f(z)^k\right)=1.\eqno{(d3)}
$$
Hence from $(d2),(d3)$, we get
$$
g^{(-1)}\left(\sum^{\infty}_{k=-\infty}c_k\xi^{k(k-1)(k-2)/6}\left(\xi z\right)^{k(k-1)/2}\left(zf(z)\right)^k\right)-
$$
$$
-g^{(-1)}\left(\sum^{\infty}_{k=-\infty}c_{k}\xi^{k(k-1)(k-2)/6}z^{k(k-1)/2}f(z)^k\right)=1.\eqno{(d4)}
$$
Hence equation $(d3)$ is of the form
$$
g^{(-1)}\left(f(z)^{-1}F\left(z,f(z)\right)\right)-g^{(-1)}\left(G\left(z,f(z)\right)\right)=1,\eqno{(d5)}
$$
where $F,G,g$ are the known functions:
$$
F(z,w)=\sum^{\infty}_{k=-\infty}c_{k-1}\xi^{k(k-1)(k-2)/6}z^{k(k-1)/2}w^k
$$
and
$$
G(z,w)=\sum^{\infty}_{k=-\infty}c_k\xi^{k(k-1)(k-2)/6}z^{k(k-1)/2}w^k,
$$
(since they have known coeficients $c_k$).\\ 
Hence assuming $\lambda(z+1)=\xi\lambda(z)$, we get from $H(z+1)=\lambda(z)H(z)$ that if $f(z)$ is  $f(z)=H\left(\lambda^{(-1)}(z)\right)$, then $f(\xi z)=zf(z)$. Also $(d5)$ is a (non-differential) transcendental equation with unkown function the $f(z)$.
Hence we have the next\\
\\ 
\textbf{Theorem 9.1}\\
If for the functions $F(z,w), G(z,w)$ defined above, we have
$$
g^{(-1)}\left(H(z)^{-1}F\left(\lambda(z),H(z)\right)\right)-g^{(-1)}\left(G\left(\lambda(z),H(z)\right)\right)=1,
$$
where $\lambda(z+1)=\xi\lambda(z)$ and $H(z+1)=\lambda(z)H(z)$, then it holds
$$
\sum^{\infty}_{k=-\infty}c_k\widehat{y}(z+k)=g(z)\widehat{y}(z),
$$ 
where
$$
\frac{\widehat{y}(z+1)}{\widehat{y}(z)}=H(z).
$$
Then we also have
$$
\widehat{y}(z)=\sum^{\infty}_{k=-\infty}A_k(z)\xi^{-k(k-1)(k-2)/6}\lambda(z)^{-k(k-1)/2}H(z)^{-k},
$$
where
$$
A_{n-1}(z+1)=A_n(z).
$$
\textbf{Remark.}\\If the function $\exp\left(2\pi i g^{(-1)}(z)\right)$ is one to one, then
$$
H(z)^{-1}F\left(\lambda(z),H(z)\right)=G\left(\lambda(z),H(z)\right)
$$
or equivalently
$$
f(z)^{-1}F\left(z,f(z)\right)=G\left(z,f(z)\right)\textrm{ and }f(\xi z)=zf(z).
$$
\\
\textbf{Theorem 9.2}\\
The differential equation
$$
\sum^{\infty}_{n=0}g_n(-i)^n\frac{d^n}{dx^n}y(x)=V(x) y(x),
$$
$$
V(x)=\sum^{\infty}_{n=-\infty}c_ne^{-in x},
$$ 
have solution
$$
y(x)=\frac{1}{2\pi}\int^{\infty}_{-\infty}\widehat{y}(z)e^{izx}dz,
$$
where
$$
\widehat{y}(z)=\sum^{\infty}_{k=-\infty}A_k(z)\xi^{-k(k-1)(k-2)/6}\lambda(z)^{-k(k-1)/2}f\left(\lambda(z)\right)^{-k},
$$
with $A_{n-1}(z+1)=A_{n}(z)$ and $\lambda(z+1)=\xi\lambda(z)$. The equations of finding $f(z)$ are both 
$$
g^{(-1)}\left(f(z)^{-1}\sum^{\infty}_{k=-\infty}c_{k-1}\xi^{k(k-1)(k-2)/6}z^{k(k-1)/2}f(z)^k\right)-
$$
$$
-g^{(-1)}\left(\sum^{\infty}_{k=-\infty}c_{k}\xi^{k(k-1)(k-2)/6}z^{k(k-1)/2}f(z)^k\right)=1\eqno{(d6)}
$$
and 
$$
f(\xi z)=zf(z).\eqno{(d7)}
$$
Note that
$$
\lambda^{(-1)}(z)=g^{(-1)}\left(\sum^{\infty}_{k=-\infty}c_k\xi^{k(k-1)(k-2)/6}z^{k(k-1)/2}f(z)^k\right).
$$
\\
\textbf{Remarks.}\\
Set now
$$
C(z):=\int^{\infty}_{-\infty}H(t)\lambda(t)^zdt\textrm{, }z\in\textbf{C}.
$$
Then
$$
C(z+1)=\int^{\infty}_{-\infty}H(t)\lambda(t)^{z+1}dt=\int^{\infty}_{-\infty}H(t)\lambda(t)\lambda(t)^{z}dt=
$$
$$
=\int^{\infty}_{-\infty}H(t+1)\lambda(t)^zdt=\int^{\infty}_{-\infty}H(t)\lambda(t-1)^zdt=
$$
$$
=\int^{\infty}_{-\infty}H(t)\left(\xi^{-1}\lambda(t)\right)^zdt=\xi^{-z}\int^{\infty}_{-\infty}H(t)\lambda(t)^zdt=\xi^{-z}C(z).
$$
Hence
$$
C(z+1)=\xi^{-z}C(z)
$$
and
$$
C(z)=C_0(\xi_1,z)\xi^{-z(z-1)/2}\textrm{, }z\in\textbf{C},
$$
$$
\lambda(z)=\xi^{z}\lambda_0(\xi_2,z)\textrm{, }\lambda(z+1)=\xi\lambda(z),
$$
where $C_0(\xi_1,z)$ and $\lambda_0(\xi_2,z)$ are $1-$periodic i.e. $C_0(\xi_1,z+1)=C_0(\xi_1,z)$ and $\lambda_0(\xi_2,z+1)=\lambda_0(\xi_2,z)$. Hence if we set
$$
\Theta(z):=\sum_{n\in\textbf{\scriptsize Z\normalsize}}C(n)\lambda(z)^{-n},
$$
then
$$
\Theta(z+1)=\sum_{n\in\textbf{\scriptsize Z\normalsize}}C(n)\lambda(z+1)^{-n}=\sum_{n\in\textbf{\scriptsize Z\normalsize}}C(n)\left(\xi\lambda(z)\right)^{-n}=
$$
$$
=\lambda(z)\sum_{n\in\textbf{\scriptsize Z\normalsize}}C(n+1)\lambda(z)^{-n-1}=\lambda(z)\Theta(z).
$$
Hence there exists constant $c_0$ such that
$$
\Theta(z)=c_0\sum^{\infty}_{n=-\infty}\xi^{-n(n-1)/2}\lambda(z)^{-n}=c_0\sum^{\infty}_{n=-\infty}\xi^{-n(n+1)/2}\lambda(z)^{n}
$$
and 
$$
\Theta(z+1)=\lambda(z)\Theta(z).
$$
Set
$$
G(z):=\int^{\infty}_{-\infty}C(t)\left(\lambda(z)\right)^{-t}dt\textrm{, }z\in\textbf{C}.
$$
Then
$$
G(z+1)=\int^{\infty}_{-\infty}C(t)\left(\lambda(z+1)\right)^{-t}dt=\int^{\infty}_{-\infty}C(t)\left(\xi\lambda(z)\right)^{-t}dt=
$$
$$
=\int^{\infty}_{-\infty}C(t)\xi^{-t}\left(\lambda(z)\right)^{-t}dt=\int^{\infty}_{-\infty}C(t+1)\left(\lambda(z)\right)^{-t}dt=
$$
$$
=\lambda(z)\int^{\infty}_{-\infty}C(t+1)\left(\lambda(z)\right)^{-t-1}dt=\lambda(z)G(z).
$$
Hence again
$$
G(z+1)=\lambda(z)G(z).
$$
\\

Now as one can see $(d7)$ can be writen in the form $f\left(\xi^{z+1}\right)=\xi^zf\left(\xi^{z}\right)$. Hence exists $C(z)$ such that $f\left(\xi^z\right)=C(z)^{-1}$. Hence $(d6)$ becomes
$$
g^{(-1)}\left(f(\xi^z)^{-1}\sum^{\infty}_{k=-\infty}c_{k-1}\xi^{k(k-1)(k-2)/6}\xi^{zk(k-1)/2}f(\xi^z)^k\right)-
$$
$$
-g^{(-1)}\left(\sum^{\infty}_{k=-\infty}c_{k}\xi^{k(k-1)(k-2)/6}\xi^{zk(k-1)/2}f(\xi^z)^k\right)=1.
$$
Hence
$$
g^{(-1)}\left(C(z)\sum^{\infty}_{k=-\infty}c_{k-1}\xi^{k(k-1)(k-2)/6+zk(k-1)/2}C(z)^{-k}\right)-
$$
$$
-g^{(-1)}\left(\sum^{\infty}_{k=-\infty}c_{k}\xi^{k(k-1)(k-2)/6+zk(k-1)/2}C(z)^{-k}\right)=1.
$$
Hence
$$
g^{(-1)}\left(\sum^{\infty}_{k=-\infty}c_{k}\xi^{k(k-1)(k+1)/6+zk(k+1)/2}C(z)^{-k}\right)-
$$
$$
-g^{(-1)}\left(\sum^{\infty}_{k=-\infty}c_{k}\xi^{k(k-1)(k-2)/6+zk(k-1)/2}C(z)^{-k}\right)=1\Leftrightarrow
$$
$$
g^{(-1)}\left(\sum^{\infty}_{k=-\infty}c_{k}\xi^{k(k-1)(k-2)/6}\xi^{k(k-1)/2}\xi^{zk(k+1)/2}C(z)^{-k}\right)-
$$
$$
-g^{(-1)}\left(\sum^{\infty}_{k=-\infty}c_{k}\xi^{k(k-1)(k-2)/6}\xi^{zk(k-1)/2}C(z)^{-k}\right)=1\Leftrightarrow
$$
$$
g^{(-1)}\left(\sum^{\infty}_{k=-\infty}c_{k}\xi^{k(k-1)(k-2)/6}\xi^{k(k-1)/2}\xi^{zk(k+1)/2}\xi^{-kz}C(z+1)^{-k}\right)-
$$
$$
-g^{(-1)}\left(\sum^{\infty}_{k=-\infty}c_{k}\xi^{k(k-1)(k-2)/6}\xi^{zk(k-1)/2}C(z)^{-k}\right)=1\Leftrightarrow
$$
$$
g^{(-1)}\left(\sum^{\infty}_{k=-\infty}c_{k}\xi^{k(k-1)(k-2)/6}\xi^{(z+1)k(k-1)/2}C(z+1)^{-k}\right)-
$$
$$
-g^{(-1)}\left(\sum^{\infty}_{k=-\infty}c_{k}\xi^{k(k-1)(k-2)/6}\xi^{zk(k-1)/2}C(z)^{-k}\right)=1.
$$

\[
\]

Also if in relation $(d1)$ we set $z\rightarrow \xi^z$, then
$$
g^{(-1)}\left(\sum^{\infty}_{k=-\infty}c_k\xi^{k(k-1)(k-2)/6}\left(\xi^{z+1}\right)^{k(k-1)/2}\left(\xi^z f(\xi^z)\right)^k\right)-
$$
$$
-g^{(-1)}\left(\sum^{\infty}_{k=-\infty}c_k\xi^{k(k-1)(k-2)/6}\left(\xi^{z}\right)^{k(k-1)/2}f(\xi^z)^k\right)=1\Leftrightarrow
$$
$$
g^{(-1)}\left(\sum^{\infty}_{k=-\infty}c_k\xi^{k(k-1)(k-2)/6}\left(\xi^{z+1}\right)^{k(k-1)/2}C(z+1)^{-k}\right)-
$$
$$
-g^{(-1)}\left(\sum^{\infty}_{k=-\infty}c_k\xi^{k(k-1)(k-2)/6}\left(\xi^{z}\right)^{k(k-1)/2}C(z)^{-k}\right)=1.
$$
If we set 
$$
A(z):=g^{(-1)}\left(\sum^{\infty}_{k=-\infty}c_k\xi^{k(k-1)(k-2)/6}\left(\xi^{z}\right)^{k(k-1)/2}C(z)^{-k}\right),
$$
then
$$
A(z+1)-A(z)=1
$$
and the function $\xi^{A(z)}$ is semiperiodic function, with period 1 and weight $\xi$. Hence
$$
\xi^{A(z)}=\xi^{z}\mu_0(z),
$$
where $\mu_0(z)$ is $1-$periodic function. Hence equivalent
$$
A(z)=z+\log_{\xi}\left(\mu_0(z)\right).
$$
Hence we get the next\\
\\
\textbf{Theorem 9.3}\\
The two equations $(d6),(d7)$ of Theorem 9.2, can be writen as:
$$
\sum^{\infty}_{k=-\infty}c_k\xi^{k(k-1)(k-2)/6}z^{k(k-1)/2}f\left(z\right)^{k}=g\left(\log_{\xi}(z)+\log_{\xi}\left(\mu_0\left(\log_{\xi}(z)\right)\right)\right),
$$
where $\mu_0(z)$ is $1-$periodic function ($\mu_0(z+1)=\mu_0(z)$).\\ 
And 
$$
f(z)=C\left(\log_{\xi}z\right)^{-1},
$$
with $C(z+1)=\xi^{-z}C(z)$, (hence $C(z)=\xi^{-z(z-1)/2}C_0(\xi_1,z)$, and $C_0(\xi_1,z)$ is $1-$periodic).\\
\\

From Theorem 9.3 and Theorem 9.1, we have\\
\\
\textbf{Theorem 9.4}\\
Assume the problem:
$$
\sum^{\infty}_{n=0}g_n(-i)^n\frac{d^n}{dx^n}y(x)=V(x) y(x)\textrm{, }g(x)=\sum^{\infty}_{n=0}g_nx^n,
$$
with periodic potential
$$
V(x)=\sum^{\infty}_{n=-\infty}c_ne^{-in x}.
$$
Then if 
$$
f(z)=C\left(\log_{\xi}z\right)^{-1},
$$
with $C(z)=\xi^{-z(z-1)/2}C_0(\xi_1,z)$ and $C_0(\xi_1,z)$ is $1-$periodic function such that 
$$
\sum^{\infty}_{k=-\infty}c_k\xi^{k(k-1)(k-2)/6}\xi^{zk(k-1)/2}\xi^{kz(z-1)/2}C_0\left(\xi_1,z\right)^{-k}=g\left(z+\log_{\xi}\left(\mu_0\left(z\right)\right)\right),
$$
where $\mu_0(z)$ is also any $1-$periodic function. Then function (solution) $y(x)$ (of the above problem) is such that
$$
\widehat{y}(z)=\sum^{\infty}_{k=-\infty}A_k(z)\xi^{-k(k-1)(k-2)/6}\lambda(z)^{-k(k-1)/2}f\left(\lambda(z)\right)^{-k},
$$
where $A_{n-1}(z+1)=A_{n}(z)$ and 
$$
\lambda^{(-1)}(z)=g^{(-1)}\left(\sum^{\infty}_{k=-\infty}c_k\xi^{k(k-1)(k-2)/6}z^{k(k-1)/2}f(z)^k\right).
$$
\\
\textbf{Therorem 9.5}\\
The two equations $(d6),(d7)$ of Theorem 9.1, can be writen as:\\ If
$$
f(z)=C\left(\log_{\xi}z\right)^{-1}
$$
and
$$
g^{(-1)}\left(\sum^{\infty}_{k=-\infty}c_{k}\xi^{k(k-1)(k-2)/6}\xi^{(z+1)k(k-1)/2}C(z+1)^{-k}\right)-
$$
$$
-g^{(-1)}\left(\sum^{\infty}_{k=-\infty}c_{k}\xi^{k(k-1)(k-2)/6}\xi^{zk(k-1)/2}C(z)^{-k}\right)=1,
$$
where $C(z)$ is such that $C(z+1)=\xi^{-z}C(z)$, then
$$
\widehat{y}(z)=\sum^{\infty}_{k=-\infty}A_k(z)\xi^{-k(k-1)(k-2)/6}\lambda(z)^{-k(k-1)/2}f\left(\lambda(z)\right)^{-k},
$$
with $A_{n-1}(z+1)=A_{n}(z)$ and $\lambda(z+1)=\xi\lambda(z)$.\\
\\

Set now
$$
F(z,w):=g^{(-1)}\left(\sum^{\infty}_{k=-\infty}c_k\xi^{k(k-1)(k-2)/6}\xi^{zk(k-1)/2}\xi^{kz(z-1)/2}w^{-k}\right),\eqno{(d8)}
$$
where the $c_k$ are given from the potential $V(z)$ and $g(z)$ are given from the DE's of Theorems 9.1-9.5. We can write also $C(z+1)=\xi^{-z}C(z)$, where $C(z)=\xi^{-z(z-1)/2}C_0(\xi_1,z)$ and $C_0(\xi_1,z)$ is an $1-$periodic function. Then we have
$$
F\left(z+1,C_0(\xi_1,z+1)\right)-F\left(z,C_0(\xi_1,z)\right)=1\eqno{(a)}
$$
and
$$
C_0(\xi_1,z+1)=C_0(\xi_1,z).\eqno{(b)}
$$
Set
$$
P(z,w):=\exp\left(2\pi i F(z,w)\right)\eqno{(c)}
$$
and
$$
G(z):=\frac{1}{C(z)}.\eqno{(c0)}
$$
Then
$$
G(z+1)=\frac{1}{C(z+1)}=\xi^zG(z),\eqno{(c1)}
$$
where $G(z):=\xi^{z(z-1)/2}G_0(\xi_1,z)$ and $G_0(\xi_1,z)$ is any $1-$periodic function. Then from $(a)$ we have 
$$
P(z+1,C_0(\xi_1,z+1))=\exp\left(2\pi i F(z+1,C_0(\xi_1,z+1))\right)=
$$
$$
=\exp\left(2\pi i F(z,C_0(\xi_1,z))\right)=P(z,C_0(\xi_1,z)).\eqno{(c2)}
$$
Also from $(d8)$, when $\xi^{\tau}=e^{2\pi i \tau}$, we have
$$
F(z+1,w)=F(z,\xi^{-z}w)\textrm{, }(z,w)\in D\times S.\eqno{(c3)}
$$
Hence
$$
F(z+1,G(z+1))=F(z,\xi^{-z}G(z+1))=F(z,G(z)).\eqno{(c4)}
$$
By this way we have
$$
P(z+1,G(z+1))=P(z,G(z))\eqno{(c5)}
$$
and
$$
P(z+1,C_0(\xi_1,z+1))=P(z,C_0(\xi_1,z)).\eqno{(c6)}
$$
Also
$$
P(z+1,w)=P(z,\xi^{-z}w).\eqno{(c7)}
$$
Hence we get the next:\\
\\
\textbf{Theorem 9.5.1}\\
Assume that equation
$$
P(z+1,Y(z))=P(z,X(z))\textrm{, }z\in D,\eqno{(c8)}
$$
have equivalent solution 
$$
Y(z)=S(z,X(z))\textrm{, }z\in D.\eqno{(c9)}
$$
Then we can write that: Equation
$$
P(z+1,C_0(\xi_1,z+1))=P(z,C_0(\xi_1,z))\textrm{, }z\in D,\eqno{(g1)}
$$ 
have equivalent unique solution
$$
C_0(\xi_1,z+1)=S(z,C_0(\xi_1,z))\textrm{, }z\in D.\eqno{(g2)}
$$
Then all following three results are equivalent to each other:\\
\textbf{i)}
$$
P(z+1,C_0(\xi_1,z+1))=P(z,C_0(\xi_1,z))
$$
and $C_0(\xi_1,z+1)=C_0(\xi_1,z)$.\\
\textbf{ii)} 
$$
C_0(\xi_1,z+1)=S(z,C_0(\xi_1,z))
$$
and $C_0(\xi_1,z+1)=C_0(\xi_1,z)$.\\
\textbf{iii)}
$$
C_0(\xi_1,z)=S(z,C_0(\xi_1,z))
$$
and
$$
S(z+1,C_0(\xi_1,z+1))=S(z,C_0(\xi_1,z+1)).
$$
\textbf{iv)}
$$
C_0(\xi_1,z)=S(z,C_0(\xi_1,z))
$$
and
$$
S(z+1,C_0(\xi_1,z))=S(z,C_0(\xi_1,z)).
$$
\\
\textbf{Proof.}\\
Case $\textbf{(i)}\Leftrightarrow \textbf{(ii)}$ is obvious.\\ 
The case $\textbf{(ii)}\Rightarrow \textbf{(iii)}$ is
$$
C_0(\xi_1,z)=C_0(\xi_1,z+1)=S(z,C_0(\xi_1,z))\Rightarrow
$$
$$
C_0(\xi_1,z+1)=S(z+1,C_0(\xi_1,z+1))=S(z,C_0(\xi_1,z+1)).
$$
The case $\textbf{(iii)}\Rightarrow \textbf{(ii)}$ is
$$
C_0(\xi_1,z)=S(z,C_0(\xi_1,z))\Rightarrow C_0(\xi_1,z+1)=S(z+1,C_0(\xi_1,z+1))=
$$
$$
=S(z,C_0(\xi_1,z+1)).
$$ 
Hence the equation $X(z)=S(z,X(z))$ have two solutions $C_0(\xi_1,z)$ and $C_0(\xi_1,z+1)$. Hence it must be $C_0(\xi_1,z)=C_0(\xi_1,z+1)$.\\
The case $\textbf{(iv)}\Rightarrow \textbf{(iii)}$ is:\\
$$
S(z+1,C_0(\xi_1,z))=S(z,C_0(\xi_1,z))=C_0(\xi_1,z)
$$
Hence
$$
C_0(\xi_1,z+1)=S(z+1,C_0(\xi_1,z+1))
$$
and
$$
S(z+1,C_0(\xi_1,z))=C_0(\xi_1,z).
$$
Hence $C_0(\xi_1,z)=C_0(\xi_1,z+1)$.\\
The case $\textbf{(iii)}\Rightarrow \textbf{(iv)}$ is
$$
C_0(\xi_1,z+1)=S(z+1,C_0(\xi_1,z+1))=S(z,C_0(\xi_1,z+1))
$$
and
$$
C_0(\xi_1,z)=S(z,C_0(\xi_1,z)).
$$
Hence
$$
C_0(\xi_1,z)=C_0(\xi_1,z+1).
$$
QED\\
\\
\textbf{Theorem 9.5.2}\\
The equations $(a),(d8),(b)$, that solve the problem of Theorems 9.1-9.5 are equivalent to $(c2),(c),(b)$ and hence to (from Theorem 9.5.1):\\
\textbf{i)}
$$
C_0(\xi_1,z)=S(z,C_0(\xi_1,z)),\eqno{(g3)}
$$
$$
C_0(\xi_1,z+1)=C_0(\xi_1,z),\eqno{(g4)}
$$
where $S(z,w)$ is defined from $(c8),(c9)$.\\
Moreover then\\
\textbf{ii)}
$$
S(z,w)\textrm{ is equivalent to } P(z,w)=\exp(2\pi i F(z,w)),\eqno{(g5)}
$$
where the function $F(z,w)$ is given from $(d8)$. The equivalence means that function $C_0(\xi_1,z)$ is solution of ''almost'' the same equations i.e.
$$
P(z+1,C_0(\xi_1,z+1))=P(z,C_0(\xi_1,z))
$$ 
and
$$
S(z+1,C_0(\xi_1,z+1))=S(z,C_0(\xi_1,z)).
$$
\textbf{iii)} 
$$
C_0(\xi_1,z)=P(z,C_0(\xi_1,z)),\eqno{(g6)}
$$
$$
C_0(\xi_1,z+1)=C_0(\xi_1,z). \eqno{(g7)}
$$
\\
\textbf{Proof.}\\
Note that from Theorem 9.5.1 we have
$$
P(z+1,C_0(\xi_1,z+1))=P(z,C_0(\xi_1,z))
$$ 
and $C_0(\xi_1,z+1)=C_0(\xi_1,z)$. Also then equivalent
$$
S(z+1,C_0(\xi_1,z+1))=S(z,C_0(\xi_1,z))
$$
and $C_0(\xi_1,z+1)=C_0(\xi_1,z)$. Hence we can say that $P(z,w)$ and $S(z,w)$ have the same ''shape''.\\
\\
\textbf{IV.} If we assume the equation $y''(x)=V(x)y(x)$, where $V(x)=\sum^{N}_{n=0}c_ne^{-in x}$, is a trigonometric polynomial, then under certain conditions we have
$$
\sum^{N}_{n=0}c_n\widehat{y}(z+n)=-z^2\widehat{y}(z).
$$
This equation is a simple special functions functional equation and we can assume is taking known solutions. For example, with $N=2$ the equation
$$
c_0f(x)+c_1 f(x+1)+c_2f(x+2)=-x^2f(x)
$$ 
gives a solution $f(x)$ and hence the equation 
$$
y''(x)=\left(\sum^{2}_{n=0}c_ne^{-in x}\right)y(x),
$$
have solution $y(x)$ such that $\widehat{y}(x)=f(x)\Leftrightarrow y(x)=\frac{1}{2\pi}\int^{\infty}_{-\infty}f(t)e^{it x}dt$.\\

More general if $g_n=g^{(n)}(0)/n!$, the equation
$$
\sum^{\infty}_{n=0}g_n(-i)^n\frac{d^n}{dx^n}y(x)=V(x)y(x),\eqno{(56.0a)}
$$ 
where the potential $V(x)$ is of the form
$$
V(x)=\sum_{n=0}^{N}c_ne^{-inx},
$$
have solution $y(x)$ such that
$$
y(x)=\frac{1}{2\pi}\int^{\infty}_{-\infty}f(t)e^{it x}dt,
$$
with $f(x)$ solution of the following functional equation
$$
c_0f(x)+c_1f(x+1)+c_2f(x+2)+\ldots+c_{N}f(x+N)=g(x)f(x)
$$
and
$$
g(x)=\sum^{\infty}_{n=0}g_nx^n.
$$
\\

Hence if $y''(x)=V(x)y(x)$, $V(x)=\sum^{N}_{n=0}c_ne^{-in x}$ and $f(x)$ is such that
$$
c_0f(x)+c_1f(x+1)+\ldots+c_Nf(x+N)=-x^2f(x),
$$
then
$$
y(x)=\frac{1}{2\pi}\int^{\infty}_{-\infty}f(t)e^{it x}dt.
$$
\\
\textbf{V.} If $g(z)$ is of the form
$$
g(z)=1+\frac{a_0(z)}{1+}\frac{-a_1(z)}{(1+a_1(z))+}\frac{-a_2(z)}{(1+a_2(z))+}\frac{-a_3(z)}{(1+a_3(z))+}\ldots
=
$$
$$
=1+\sum^{\infty}_{k=1}\prod^{k-1}_{j=0}a_j(z),
$$
where $a_j(z)=\frac{A}{B}c_j^{*}H(z+j)$, then
$$
g(z)=1+\frac{c_{0}^{*}H(z)}{1+}\frac{-\frac{A}{B}c_1^{*}H(z+1)}{\left(1+\frac{A}{B}c_1^{*}H(z+1)\right)+}\frac{-\left(\frac{A}{B}\right)^2c_2^{*}H(z+2)}{\left(1+\left(\frac{A}{B}\right)^2c_2^{*}H(z+2)\right)+}\ldots.\eqno{(cf)}
$$
Hence given $g(z)$, we expand it in the form $(cf)$ and find $H(z)$. Then using Weierstrass theorem
$$
H(z)=C z^m e^{-R(z)}\frac{P(z)}{Q(z)},
$$
or (56.00), we find the solution
$$
\widehat{y}(z)=\left(\frac{A}{B}\right)^z C^z \Gamma(z)^m e^{-R_1(z)} \prod^{\infty}_{j=1}\frac{\Gamma(z-\rho_j)}{\Gamma(z-\sigma_j)}.
$$
\\
\textbf{Theorem 9.6}\\
Assume that
$$
f(z,w)=\sum^{\infty}_{n=-\infty}A_n(z)w^{-n}.
$$
\textbf{i)} If $A_n(z)$ satisfies the relation
$$
A_{n-1}(z+1)=A_n(z),\eqno{(f)}
$$
then equivalent we have that
$$
f\left(z+1,w\right)=wf\left(z,w\right).
$$
\textbf{ii)} Also if 
$$
\widehat{y}(z)=\sum_{n}A_{n}(z)\xi^{-n(n-1)(n-2)/6}\lambda(z)^{-n(n-1)/2}H(z)^{-n}
$$
and holds $(f)$ and $\lambda(z+1)=\xi\lambda(z)$, then:
$$
H(z+1)=\lambda(z)H(z)
$$
is equivalent to
$$
\widehat{y}(z+1)=H(z)\widehat{y}(z).
$$
\\
\textbf{Proof.}\\
The relation $A_{n-1}(z+1)=A_n(z)$ is equivalent to
$$
\sum^{\infty}_{n=-\infty}A_{n-1}(z+1)w^{-n}=\sum^{\infty}_{n=-\infty}A_n(z)w^{-n}\Leftrightarrow
$$
$$
\sum^{\infty}_{n=-\infty}A_{n}(z+1)w^{-n}=w\sum^{\infty}_{n=-\infty}A_n(z)w^{-n}\Leftrightarrow
$$
$$
f\left(z+1,w\right)=wf(z,w).
$$
For the other statements, we have
$$
\widehat{y}(z+1)=
$$
$$
=\sum^{\infty}_{n=-\infty}A_n(z+1)\xi^{-n(n-1)(n-2)/6}\xi^{-n(n-1)/2}\lambda(z)^{-n(n-1)/2}\lambda(z)^{-n}H(z+1)^{-n}=
$$
$$
=\sum^{\infty}_{n=-\infty}A_n(z+1)\xi^{-n(n-1)(n-2)/6}\xi^{-n(n-1)/2}\lambda(z)^{-n(n-1)/2}\lambda(z)^{-n}H(z)^{-n}=
$$
$$
=\sum^{\infty}_{n=-\infty}A_{n}(z+1)\xi^{-n(n-1)(n+1)/6}\lambda(z)^{-n(n+1)/2}H(z)^{-n}=
$$
$$
=H(z)\sum^{\infty}_{n=-\infty}A_{n-1}(z+1)\xi^{-n(n-1)(n-2)/6}\lambda(z)^{-n(n-1)/2}H(z)^{-n}=
$$
$$
=H(z)\widehat{y}(z).
$$
Assuming that $\widehat{y}(z+1)=H(z)\widehat{y}(z)$, we have
$$
\sum_{n}A_n(z+1)\xi^{-n(n-1)(n-2)/6}\lambda(z+1)^{-n(n-1)/2}H(z+1)^{-n}=
$$
$$
=H(z)\sum_{n}A_n(z)\xi^{-n(n-1)(n-2)/6}\lambda(z)^{-n(n-1)/2}H(z)^{-n}\Leftrightarrow 
$$
$$
\sum_{n}A_{n}(z+1)\xi^{-n(n-1)(n-2)/6}\xi^{-n(n-1)/2}\lambda(z)^{-n(n-1)/2}H(z+1)^{-n}=
$$
$$
=H(z)\sum_{n}A_n(z)\xi^{-n(n-1)(n-2)/6}\lambda(z)^{-n(n-1)/2}H(z)^{-n}\Leftrightarrow 
$$
$$
\sum_{n}A_n(z+1)\xi^{-n(n-1)(n+1)/2}\lambda(z)^{-n(n+1)/2}\lambda(z)^{n}H(z+1)^{-n}=
$$
$$
=H(z)\sum_{n}A_n(z)\xi^{-n(n-1)(n-2)/6}\lambda(z)^{-n(n-1)/2}H(z)^{-n}\Leftrightarrow
$$
$$
\sum_{n}A_{n-1}(z+1)\xi^{-n(n-1)(n-2)/6}\lambda(z)^{-n(n-1)/2}\lambda(z)^{n-1}H(z+1)^{-n+1}=
$$
$$
=H(z)\sum_{n}A_{n}(z)\xi^{-n(n-1)(n-2)/6}\lambda(z)^{-n(n-1)/2}H(z)^{-n}\Leftrightarrow
$$
$$
\left(\lambda(z)^{-1}H(z+1)\right)\sum_{n}A_{n}(z)\xi^{-n(n-1)(n-2)/6}\lambda(z)^{-n(n-1)/2}\left(\lambda(z)^{-1}H(z+1)\right)^{-n}=
$$
$$
=H(z)\sum_{n}A_{n}(z)\xi^{-n(n-1)(n-2)/6}\lambda(z)^{-n(n-1)/2}H(z)^{-n}.
$$
Hence obviously
$$
H(z+1)=\lambda(z)H(z).
$$
\\
\textbf{Theorem 9.7}\\
With the notation of Theorem 9.0a we have
$$
H(z+j)=\left(\prod^{j-1}_{k=0}\lambda(z+k)\right)H(z)
$$
and
$$
\widehat{y}(z+n)=\left(\prod^{n-1}_{j=0}\prod^{j-1}_{k=0}\lambda(z+k)\right)H(z)^n\widehat{y}(z).
$$
Hence
$$
\sum^{\infty}_{n=0}c_n\widehat{y}(z+n)=g(z)\widehat{y}(z)
$$
and
$$
g(z)=\sum^{\infty}_{n=0}c_n\left(\prod^{n-1}_{j=0}\prod^{j-1}_{k=0}\lambda(z+k)\right)H(z)^n.
$$
Also then
$$
\widehat{y}(z)=\sum^{\infty}_{n=-\infty}A_n(z)H(z)^n,
$$
where 
$$
A_{n+1}(z+1)\lambda(z)^{n+1}=A_n(z).
$$
\\
\textbf{Theorem 9.8}
$$
g(z)=\sum_nc_n\prod^{n-1}_{j=0}H\left(z+\frac{2\pi j}{T}\right),
$$
then if $H(z)=\exp\left(H_1(z)\right)$, we have
$$
g(z)=\sum_nc_n\exp\left(\sum^{n-1}_{j=0}H_1\left(z+\frac{2\pi j}{T}\right)\right).
$$
Hence assuming the equation
$$
\sum^{n-1}_{j=0}H_{1,n}\left(z+\frac{2\pi j}{T}\right)=g(z)H_{1,n}(z),
$$
we have a solution
$$
\sum^{n-1}_{j=0}f_n\left(z+\frac{2\pi j}{T}\right)=g(z)f_n(z)
$$
and $f_n(z)=H_{1,n}(z)$ and if $\lim_{n\rightarrow\infty}H_{1,n}(z)=\lim_{n\rightarrow\infty}f_n(z)=f(z)$. Hence
$$
\sum^{\infty}_{j=0}f\left(z+\frac{2\pi j}{T}\right)=g(z)f(z).
$$
\\
\textbf{Theorem 9.9}\\
Assume
$$
V(z)=\sum^{\infty}_{n=0}c_ne^{-in z}
$$
and
$$
H(z)=\frac{(z+A_1)(z+A_2)\ldots (z+A_{N})}{(z+B_1)(z+B_2)\ldots (z+B_{M})}e^{-a z}.
$$ 
Then
$$
g(z)=\sum^{\infty}_{k=0}\left(\frac{A}{B}\right)^k c_k \frac{(z+A_1)_k(z+A_2)_k\ldots (z+A_{N})_k}{(z+B_1)_k(z+B_2)_k\ldots (z+B_{M})_k}\frac{e^{-a k z}e^{-ak(k+1)/2}}{k!}.
$$
\\
\textbf{Proof.}\\
Easy follows from (53) of Theorem 9.\\
\\

If $f:\textbf{N}\rightarrow\textbf{N}$ is increasing function and $\chi:\textbf{N}\rightarrow \textbf{C}$, then the function
$$
\theta(q)=\sum_{n=1}^{\infty}\chi(n)q^{f(n)},
$$
can be written as Lambert series:
$$
\theta(q)=\sum^{\infty}_{n=1}\frac{B(n)q^n}{1-q^n},
$$
where
$$
B(n)=\sum_{f(d)|n}\chi(d)\mu\left(\frac{n}{f(d)}\right).
$$
Also then
$$
\theta(q)=\sum^{\infty}_{n=1}q^n\sum_{d|n}B(d)=\sum^{\infty}_{n=1}q^nA(n),
$$
where
$$
A(n)=\sum_{d|n}\sum_{f(\delta)|d}\chi(\delta)\mu\left(\frac{d}{f(\delta)}\right).
$$
If $A(n)$ is $T_0-$periodic, $T_0\in\textbf{N}$, $T_0>1$, then
$$
\theta(q)=\left(1-q^{T_0}\right)^{-1}\sum^{T_0}_{n=1}A(n)q^n,
$$
is rational function of $q$.\\
\\
\textbf{Theorem 10.}\\
Given $g(x)=\sum^{\infty}_{n=0}g_nx^n$ and $V(x)=\sum^{\infty}_{n=0}c_ne^{-2\pi i n x/T}$, the equation 
$$
\sum^{\infty}_{n=0}g_n(-i)^ny^{(n)}(x)=V(x)y(x),
$$
is equivalent to
$$
\sum^{\infty}_{n=0}c_n\widehat{y}\left(z+\frac{2\pi}{T}n\right)=g(z)\widehat{y}(z).
$$
If $H(z)$ is defined from
$$
\widehat{y}\left(z+\frac{2\pi}{T}\right)=H(z)\widehat{y}(z),
$$
then it must hold
$$
g(z)=\sum_{n=0}^{\infty}c_n\prod^{n-1}_{j=0}H\left(z+\frac{2\pi}{T}j\right).
$$
Hence assuming that $H(z)$ is of the form
$$
H(z)=Cz^me^{-R(z)}\theta(z),
$$
where $R(z)$ holomorphic function and $\theta(z)$ meromophic function such that 
$$
\theta(z)=\sum_{n=1}^{\infty}\chi(n)z^{f(n)},
$$
with $f:\textbf{N}\rightarrow\textbf{N}$ is increasing function and $\chi:\textbf{N}\rightarrow \textbf{C}$ and the function 
$$
A(n)=\sum_{d|n}\sum_{f(\delta)|d}\chi(\delta)\mu\left(\frac{d}{f(\delta)}\right),
$$
is $T_0-$periodic, $T_0\in\textbf{N}$, $T_0>1$. Then
$$
\theta(z)=\left(1-z^{T_0}\right)^{-1}\sum^{T_0}_{n=1}A(n)z^n,
$$
is rational function of $z$. Hence
$$
y(x)=\frac{1}{2\pi}\int^{+\infty}_{-\infty}C^z\Gamma(z)^me^{-R_1(z)}\left(\frac{\prod^{T_0}_{j=1}\Gamma(z-\rho_j)}{\prod^{T_0}_{j=1}\Gamma(z-\sigma_j)}\right)e^{i z x}dz.
$$
\\
\textbf{Definition 10.1}\\
We define the operator $\textbf{T}$ such that if
$$
g(z)=\sum^{\infty}_{n=0}c_n\prod^{n-1}_{j=0}H(z+j),
$$
then
$$
\left(\textbf{T}g\right)(z)=\sum^{\infty}_{n=0}c_{n-1}\prod^{n-1}_{j=0}H(z+j),
$$
with $c_{-1}=0$ and $\prod^{-1}_{j=0}a_j=1$.\\
\\
\textbf{Theorem 10.2}
$$
H(z)=\frac{\left(\textbf{T}g\right)(z)}{g(z+1)}.
$$
\\ 
\textbf{Proof.}\\
We have easily
$$
g(z+1)=\sum^{\infty}_{n=0}c_n\prod^{n-1}_{j=0}H(z+j+1)=H(z)^{-1}\sum^{\infty}_{n=0}c_n\prod^{n}_{j=0}H(z+j)=
$$
$$
=H(z)^{-1}\sum^{\infty}_{n=0}c_{n-1}\prod^{n-1}_{j=0}H(z+j)=H(z)^{-1}\left(\textbf{T}g\right)(z).
$$
\\
\textbf{Theorem 10.3}\\
Assume that $e_n(z)$ is a base in $\textbf{A}$ subset of $\textbf{R}$ and $g^{*}_n=\left\langle g(z),e_n(z)\right\rangle$, where $g(z)=\sum^{\infty}_{n=0}g^{*}_{n}e_n(x)$. Also 
$$
V(z)=\sum^{\infty}_{n=0}c_ne^{-i n z}
$$
and 
$$
\sum^{\infty}_{n=0}g_n(-i)^n\frac{d^n}{dx^n}y(z)=V(x)y(x).
$$
For to solve the above equation, set $J_{n,k}(a_{(.)})$ to be an operator such that:
$$
\exp\left(\sum^{\infty}_{k=0}a_kB_{k,n}(z)\right)=\sum^{\infty}_{k=0}J_{n,k}(a_{(.)})B_{k,n}(z),
$$
where
$$
B_{k,n}(z)=\sum^{k}_{l=0}C_{k,l}p_{k-l}(n)z^l\textrm{, }p_{\tau}(n)=\sum^{n-1}_{j=0}j^{\tau}.
$$
If $e_n(z)$ is a certain base, we write 
$$
B_{k,n}(z)=\sum^{\infty}_{l=1}S_{k,n,l}e_l(z).
$$
and $S^{*}_{k,n,l}$ is the matrix inverse of $S_{k,n,l}$, in the sense
$$
\sum^{\infty}_{l=1}S_{k,n,l}S^{*}_{k',n,l}=R_n\delta_{k',k},
$$
where $\delta_{k',k}=1$, if $k=k'$ and $\delta_{k',k}=0$, if $k'\neq k$. Also we define $a'_{n,k}$ to be the arithmetic inverse of $\sum_{d|n}c_d\mu(n/d)$, in the sense
$$
\sum_{d_1|A}\left(\sum_{d|d_1}c_d\mu(d_1/d)\right)a'_{A/d_1,B}=\delta_{A,1}\delta_{B,1}.
$$
Then with a suitable choise of $R_n$, we get
$$
y(x)=\frac{1}{2\pi}\int^{\infty}_{-\infty}\prod^{\infty}_{n=0}\exp\left(-\sum^{\infty}_{k=0}J_{n,k}^{(-1)}\left(G^{-1}\sum^{\infty}_{l=1}g^{*}_l\sum_{d_1|n,d_2|l}S^{*}_{k,d_1,d_2}a'_{n/d_1,l/d_2}\right)(z+n)^k\right)e^{izx}dz.
$$
\\
\textbf{Proof.}\\
We write
$$
g(z)=\sum^{\infty}_{n=0}c_n\prod^{n-1}_{j=0}H(z+j),
$$
where 
$$
H(z)=\exp\left(\sum^{\infty}_{k=0}a_kz^k\right)
$$
Hence if $C_{n,k}$ denotes the binomial function i.e. $C_{n,k}=\frac{n!}{(n-k)!k!}$, then
$$
\prod^{n-1}_{j=0}H(z+j)=\exp\left(\sum^{\infty}_{k=0}a_k\sum^{n-1}_{j=0}(z+j)^k\right)=
$$
$$
=\exp\left(\sum^{\infty}_{k=0}a_k\sum^{n-1}_{j=0}\sum^{k}_{l=0}C_{k,l}j^{k-l}z^l\right)=
$$
$$
=\exp\left(\sum^{\infty}_{k=0}a_k\sum^{k}_{l=0}C_{k,l}\sum^{n-1}_{j=0}j^{k-l}z^l\right)=
$$
$$
=\exp\left(\sum^{\infty}_{k=0}a_k\sum^{k}_{l=0}C_{k,l}p_{k-l}(n)z^l\right).
$$
Hence if we denote
$$
B_{k,n}(z):=\sum^{k}_{l=0}C_{k,l}p_{k-l}(n)z^l,
$$
then
$$
g(z)=\sum^{\infty}_{n=0}c_n\exp\left(\sum^{\infty}_{k=0}a_kB_{k,n}(z)\right).
$$
Now we assume the transformation $J$ such that
$$
\exp\left(\sum^{\infty}_{k=0}a_kB_{k,n}(z)\right)=\sum^{\infty}_{k=0}J_{n,k}(a_{(.)})B_{k,n}(z).
$$
This equivalently give (if we denote $h_{n,k}=J_{n,k}(a_{(.)})$):
$$
\exp\left(\sum^{n-1}_{j=0}\sum^{\infty}_{k=0}a_k(z+j)^k\right)=\sum^{n-1}_{j=0}\sum^{\infty}_{k=0}h_{n,k}(z+j)^k\Leftrightarrow
$$
$$
\prod^{n-1}_{j=0}\exp\left(\sum^{\infty}_{k=0}a_k(z+j)^k\right)=\sum^{n-1}_{j=0}\log\left(\exp\left(\sum^{\infty}_{k=0}h_{n,k}(z+j)^k\right)\right)\Leftrightarrow
$$
$$
\prod^{n-1}_{j=0}\exp\left(\sum^{\infty}_{k=0}a_k(z+j)^k\right)=\log\left(\prod^{n-1}_{j=0}\exp\left(\sum^{\infty}_{k=0}h_{n,k}(z+j)^k\right)\right)\Leftrightarrow
$$
$$
\exp\left(\prod^{n-1}_{j=0}\exp\left(\sum^{\infty}_{k=0}a_k(z+j)^k\right)\right)=\prod^{n-1}_{j=0}\exp\left(\sum^{\infty}_{k=0}h_{n,k}(z+j)^k\right).
$$
Hence the problem reduces to find if for any function $X(z)$ exists function $Y_n(z)$ such that 
$$
\prod^{n-1}_{j=0}X(z+j)=\sum^{n-1}_{j=0}Y_n(z+j)\textrm{, }\forall n\in\textbf{N}.
$$
However
$$
g(z)=\sum^{\infty}_{n=1}c_n\sum^{\infty}_{k=1}J_{n,k}(a_{(.)})B_{k,n}(z)=\sum^{\infty}_{n=1}\sum^{\infty}_{k=1}c_nJ_{n,k}(a_{(.)})B_{k,n}(z).\eqno{(a)}
$$
Assume now $e_n(z)$ is any base function in a set $\textbf{A}$ subset of $\textbf{R}$, then
$$
B_{k,n}(z)=\sum^{\infty}_{l=1}S_{k,n,l}e_{l}(z)
$$
and relation $(a)$ becomes
$$
g(z)=\sum^{\infty}_{l=1}g^{*}_le_{l}(z)=\sum^{\infty}_{n=1}\sum^{\infty}_{k=1}c_nJ_{n,k}(a_{(.)})\sum^{\infty}_{l=1}S_{k,n,l}e_l(z)=
$$
$$
=\sum^{\infty}_{l=1}\left(\sum^{\infty}_{n,k=1}c_nJ_{n,k}(a_{(.)})S_{k,n,l}\right)e_l(z).
$$
Hence
$$
g^{*}_l=\sum^{\infty}_{n,k=1}c_nJ_{n,k}(a_{(.)})S_{k,n,l}.
$$
If
$$
\sum^{\infty}_{n=1}I_{n,l}c_n=g_l^{*},
$$
where
$$
I_{n,l}=\sum^{\infty}_{k=1}J_{n,k}(a_{(.)})S_{k,n,l}
$$
and
$$
\sum^{\infty}_{l=1}S_{k,n,l} S^{*}_{k',n,l}=R_n\delta_{k,k'},
$$
then
$$
\sum^{\infty}_{l=1}I_{n,l}S^{*}_{k',n,l}=\sum^{\infty}_{k=1}J_{n,k}(a_{(.)})\sum^{\infty}_{l=1}S_{k,n,l}S^{*}_{k',n,l}=
$$
$$
=\sum^{\infty}_{k=1}J_{n,k}(a_{(.)})R_n\delta_{k,k'}=J_{n,k'}(a_{(.)})R_n.
$$
Assume
$$
\sum^{\infty}_{n_1,l_1=1}I_{nn_1,ll_1}b_{n_1,l_1}=\phi(n,l)
$$
Then
$$
\sum^{\infty}_{n_1,l_1=1}I_{nn_1n_2,ll_1l_2}b_{n_1,l_1}=\phi(nn_2,ll_2)
$$
$$
\sum^{\infty}_{n_1,l_1=1}I_{nn_1n_2,ll_1l_2}b_{n_1,l_1}a'_{n_2,l_2}=\phi(nn_2,ll_2)a'_{n_2,l_2}
$$
$$
\sum^{\infty}_{n_2,l_2=1}\sum^{\infty}_{n_1,l_1=1}I_{nn_1n_2,ll_1l_2}b_{n_1,l_1}a'_{n_2,l_2}=\sum^{\infty}_{n_2,l_2=1}\phi(nn_2,ll_2)a'_{n_2,l_2}
$$
$$
\sum^{\infty}_{A,B=1}I_{nA,lB}\sum_{d_1|A,d_2|B}b_{d_1,d_2}a'_{A/d_1,B/d_2}=\sum^{\infty}_{n_2,l_2=1}\phi(nn_2,ll_2)a'_{n_2,l_2}
$$
Assume
$$
\sum_{d_1|A,d_2|B}b_{d_1,d_2}a'_{A/d_1,B/d_2}=\delta_{A,1}\delta_{B,1}.\eqno{(b)}
$$
Then
$$
I_{n,l}=\sum^{\infty}_{n_2,l_2=1}\phi(nn_2,ll_2)a'_{n_2,l_2}
$$
Assume now that $b_{n_1,l_1}=\delta_{l_1,1}c^{*}_{n_1}$, then
$$
\sum^{\infty}_{n=1}\sum^{\infty}_{n_1,l_1=1}I_{nn_1,ll_1}b_{n_1,l_1}=\sum^{\infty}_{n=1}\phi(n,l)\Rightarrow
$$
$$
\sum^{\infty}_{n,l_1=1}I_{n,ll_1}\sum_{d|n}b_{d,l_1}=\sum^{\infty}_{n=1}\phi(n,l)=g_l^{*}.
$$
Hence if 
$$
\sum_{d|n}b_{d,l_1}=\delta_{l_1,1}c_{n}\Leftrightarrow b_{n,l}=\delta_{l,1}\sum_{d|n}c_{d}\mu(n/d)=\delta_{l,1}c^{*}_n,
$$
then
$$
\sum^{\infty}_{n=1}I_{n,l}c_n=g^{*}_l.
$$
Also if $a'_{n,l}$ is the inverse arithmetic function of $b_{n,l}$ in the sense $(b)$, then
$$
\sum_{d_1|A,d_2|B}\delta_{d_2,1}\left(\sum_{d|d_1}c_d\mu(d_1/d)\right)a'_{A/d_1,B/d_2}=\delta_{A,1}\delta_{B,1}\Leftrightarrow
$$
$$
\sum_{d_1|A}\left(\sum_{d|d_1}c_d\mu(d_1/d)\right)a'_{A/d_1,B}=\delta_{A,1}\delta_{B,1}.\eqno{(c)}
$$
From equation $(c)$ we get the $a'_{n,l}$. Finaly
$$
B_{k,n}(z)=\sum^{\infty}_{l=1}S_{k,n,l}e_l(z)
$$
$$
\sum^{\infty}_{l=1}S_{k,n,l}S^{*}_{k',n,l}=R_n\delta_{k,k'}
$$
$$
J_{n,k}(a_{(.)})=\frac{1}{R_n}\sum^{\infty}_{l=1}I_{n,l}S^{*}_{k,n,l}
$$
and
$$
I_{n,l}=\sum^{\infty}_{n_2,l_2=1}\phi(nn_2,ll_2)a'_{n_2,l_2},
$$
$$
\sum^{\infty}_{n=1}\phi(n,l)=g^{*}_{l}.
$$
Assuming further that $\phi(n,l)=\phi_1(n)\phi_2(l)$, we get
$$
\sum^{\infty}_{n=1}\phi_1(n)\phi_2(l)=g^{*}_{l}.
$$
Hence
$$
\phi_2(l)=\frac{g^{*}_l}{G}\textrm{, }G=\sum^{\infty}_{n=1}\phi_1(n)
$$
Hence there exists constant $G$ such that
$$
\phi(n,l)=\phi_1(n)\frac{g^{*}_l}{G}
$$
and
$$
I_{n,l}=\frac{1}{G}\sum^{\infty}_{n_2,l_2=1}\phi_1(nn_2)g^{*}_{ll_2}a'_{n_2,l_2}
$$
Also then 
$$
J_n(a_k)=\frac{1}{GR_n}\sum^{\infty}_{l,n_2,l_2=1}\phi_1(nn_2)g^{*}_{ll_2}a'_{n_2,l_2}S^{*}_{k,n,l}
$$
$$
g(z)=\frac{1}{G}\sum^{\infty}_{n,k=1}\frac{c_n}{R_n}\sum^{\infty}_{l,n_2,l_2=1}\phi_1(nn_2)g^{*}_{ll_2}a'_{n_2,l_2}S^{*}_{k,n,l}B_{k,n}(z).
$$
$$
g(z)=\frac{1}{G}\sum^{\infty}_{n,k=1}\frac{c_n}{R_n}\sum^{\infty}_{l=1}\left(\sum^{\infty}_{n_2,l_2=1}\phi_1(nn_2)g^{*}_{ll_2}a'_{n_2,l_2}\right)S^{*}_{k,n,l}B_{k,n}(z).
$$
Assume that $\psi(n)$ is such that
$$
\sum^{\infty}_{n=1}J_{n,k}(a_{(.)})\psi(n)=\sum^{\infty}_{n=1}\frac{\psi(n)}{GR_n}\sum^{\infty}_{l,n_2,l_2=1}\phi_1(nn_2)g^{*}_{ll_2}a'_{n_2,l_2}S^{*}_{k,n,l}.
$$
Hence
$$
\sum^{\infty}_{n=1}J_{n,k}(a_{(.)})\psi(n)=\sum^{\infty}_{n=1}\sum^{\infty}_{l,n_2,l_2=1}\frac{\psi(n)}{GR_n}\phi_1(nn_2)g^{*}_{ll_2}a'_{n_2,l_2}S^{*}_{k,n,l}\Rightarrow
$$
$$
\sum^{\infty}_{n=1}J_{n,k}(a_{(.)})\psi(n)=\sum^{\infty}_{n,l=1}\phi_1(n)g^{*}_{l}\sum_{d_1|n,d_2|l}\frac{\psi(d_1)}{GR_{d_1}}S^{*}_{k,d_1,d_2}a'_{n/d_1,l/d_2}=
$$
$$
=\sum^{\infty}_{n=1}\phi_1(n)\sum^{\infty}_{l=1}g^{*}_{l}\sum_{d_1|n,d_2|l}S^{*}_{k,d_1,d_2}a'_{n/d_1,l/d_2}.
$$
If we set now $\phi_1(n)=\psi(n)=R_n$, ($R_n$ is any arithmetical function), then
$$
\sum^{\infty}_{n=1}J_{n,k}(a_{(.)})R_{n}=\sum^{\infty}_{n=1}R_nG^{-1}\sum^{\infty}_{l=1}g^{*}_l\sum_{d_1|n,d_2|l}S^{*}_{k,d_1,d_2}a'_{n/d_1,l/d_2},
$$
where the $a'_{n,m}$ is given from $(c)$. Hence
$$
J_{n,k}\left(a_{(.)}\right)=G^{-1}\sum^{\infty}_{l=1}g^{*}_l\sum_{d_1|n,d_2|l}S^{*}_{k,d_1,b_2}a'_{n/d_1,l/d_2}
$$
Using now
$$
H(z)=\exp\left(\sum^{\infty}_{k=0}a_kz^k\right),
$$
we get
$$
y(x)=\frac{1}{2\pi}\int^{\infty}_{-\infty}\prod^{\infty}_{n=0}\exp\left(-\sum^{\infty}_{k=0}J_{n,k}^{(-1)}\left(G^{-1}\sum^{\infty}_{l=1}g^{*}_l\sum_{d_1|n,d_2|l}S^{*}_{k,d_1,d_2}a'_{n/d_1,l/d_2}\right)(z+n)^k\right)e^{izx}dz.
$$
\\
\textbf{Theorem 10.4}\\Given a potential
$$
V(z;\tau)=\sum^{\infty}_{n=0}c_n(\tau)e^{-i \tau n z}
$$
and
$$
B(z;\tau)=\sum^{\infty}_{n=0}b_n(\tau)z^n=g(-iz;\tau)=\sum^{\infty}_{n=0}(-i)^ng_n(\tau)z^n,
$$
we form the equation
$$
\sum^{\infty}_{n=0}b_n(\tau)\frac{d^n}{dz^n}y(z;\tau)=V(z;\tau)y(z;\tau),
$$
which may be formed as
$$
V\left(i\frac{d}{dz}\right)\widehat{y}(z;\tau)=g(z;\tau)\widehat{y}(z;\tau),
$$
where
$$
\widehat{y}(z;\tau)=\int^{\infty}_{-\infty}y(x;\tau)e^{-ix z}dx.
$$
Also this may be written as
$$
\sum^{\infty}_{n=0}c_n(\tau)\widehat{y}(z+\tau n;\tau)=g(z;\tau)\widehat{y}(z;\tau).\eqno{(56.1a)}
$$
\textbf{I.} If we assume that
$$
\widehat{y}(z;\tau)=\sum^{\infty}_{l=1}C_l(z;\tau),\eqno{(a)}
$$
where
$$
C_l(z+n\tau;\tau)=A_n(z;\tau)C_l(z;\tau)\textrm{, }n\in\textbf{Z},
$$
then it must be
$$
g(z;\tau)=\sum^{\infty}_{n=0}c_n(\tau)A_n(z;\tau).\eqno{(b)}
$$
From where we find $A_n(z;\tau)$. In this case also
$$
\widehat{y}(z+n\tau;\tau)=A_n(z;\tau) \widehat{y}(z;\tau).\eqno{(c)}
$$
\textbf{II.} If $H(z;\tau)$ is such that
$$
\frac{\widehat{y}(z+\tau;\tau)}{\widehat{y}(z;\tau)}=H(z;\tau)\textrm{, }Im(\tau),Im(z)>0\eqno{(d)}
$$
and $H(z+\tau;\tau)=A^{*}(z;\tau) H(z;\tau)$, then
$$
H(z+j\tau;\tau)=\left(\prod^{j-1}_{l=0}A^{*}(z+l\tau;\tau)\right) H(z;\tau)
$$
and from $(d)$:
$$
\widehat{y}(z+n\tau;\tau)=\left(\prod^{n-1}_{j=0}H(z+j\tau;\tau)\right) \widehat{y}(z;\tau).\eqno{(e)}
$$
Hence from $(c)$:
$$
A_n(z;\tau)=\prod^{n-1}_{j=0}H(z+j\tau;\tau).\eqno{(56.2)}
$$ 
Hence\\
\textbf{III.}
$$
A_1(z;\tau)=H(z;\tau)
$$
and from $(b)$ we get
$$
\widehat{y}(z;\tau)=\prod^{\infty}_{n=0}\frac{1}{A_1(z+n\tau;\tau)}.\eqno{(f)}
$$
\\
\textbf{Proof.}\\
For the proof of (II), we have
$$
g(z;\tau)=\sum^{\infty}_{n=0}c_n(\tau)\prod^{n-1}_{j=0}H(z+j\tau;\tau).\eqno{(56.4)}
$$
Hence from (56.2), we get 
$$
g(z;\tau)=\sum^{\infty}_{n=0}c_n(\tau)\left(\prod^{n-1}_{j=0}\prod^{j-1}_{l=0}A^{*}(z+l\tau;\tau)\right)H(z;\tau)^n.
$$
Hence given the $(-i)^ng_n(\tau)=b_n(\tau)$ of $B(z;\tau)$ and $c_n(\tau)$ of the potential, we get $H(z;\tau)$ from (56.2) and the solution of $(eq)$ is 
$$
\widehat{y}(z;\tau)=\prod^{\infty}_{n=0}\frac{1}{H(z+n\tau;\tau)}.
$$
We have
$$
\sum^{\infty}_{k=0}b_k(\tau)y^{(k)}(x;\tau)=V(x;\tau)y(x;\tau)\Leftrightarrow 
$$
$$
\sum^{\infty}_{k=0}b_k(\tau)\int^{\infty}_{-\infty}y^{(k)}(x;\tau)e^{-ix\gamma}dx=\int^{\infty}_{-\infty}y(x;\tau)\left(\sum^{\infty}_{n=0}c_n(\tau)x^n\right)e^{-ix\gamma}dx\Leftrightarrow
$$
$$
\sum^{\infty}_{k=0}b_k(\tau)(i\gamma)^k\int^{\infty}_{-\infty}y(x;\tau)e^{-ix\gamma}dx=\sum^{\infty}_{n=0}c_n(\tau)\int^{\infty}_{-\infty}y(x;\tau)x^ne^{-i x\gamma}dx\Leftrightarrow
$$
$$
B(i\gamma;\tau)\widehat{y}(\gamma)=\sum^{\infty}_{n=0}c_n(\tau)i^n\frac{d^n}{d\gamma^n}\left(\int^{\infty}_{-\infty}y(x;\tau)e^{-ix\gamma}dx\right)\Leftrightarrow
$$
$$
B(i\gamma;\tau)\widehat{y}(\gamma;\tau)=\sum^{\infty}_{n=0}c_n(\tau)i^n\frac{d^n}{d\gamma^n}\widehat{y}(\gamma;\tau)\Leftrightarrow
$$
$$
B(i\gamma;\tau)\widehat{y}(\gamma;\tau)=V\left(i\frac{d}{d\gamma}\right)\widehat{y}(\gamma;\tau).\eqno{(57)}
$$
Assume now that $V(x;\tau)=\sum^{\infty}_{n=0}c_n(\tau)e^{-i n\tau x}$. Hence
$$
B(i\gamma;\tau)\widehat{y}(\gamma;\tau)=\sum^{\infty}_{n=0}c_n(\tau)\exp\left(\tau n\frac{d}{d\gamma}\right)\widehat{y}(\gamma;\tau).\eqno{(58)}
$$
Or equivalent
$$
B(i\gamma;\tau)\widehat{y}(\gamma;\tau)=\sum^{\infty}_{n=-\infty}c_n(\tau)\widehat{y}(\gamma+n\tau;\tau).\eqno{(59)}
$$
Assume that $Im(z)$, $Im(\tau)>0$ and
$$
\widehat{y}(z;\tau)=\sum^{\infty}_{l=1}C_l(z;\tau).\eqno{(60)}
$$
Then if $n$ is integer, we have
$$
C_l(z+n\tau;\tau)=A_{n}(z;\tau) C_l(z;\tau).
$$
Hence for almost all $a_n(\tau)$ we have
$$
\sum^{\infty}_{n=0}a_n(\tau)\widehat{y}(z+n\tau)
=\sum^{\infty}_{n=0}a_n(\tau)\sum^{\infty}_{l=1}C_l(z+n\tau;\tau)=
$$
$$
=\sum^{\infty}_{n=0}\sum^{\infty}_{l=1}a_n(\tau)A_{n}(z;\tau)C_l(z;\tau)=
\sum^{\infty}_{n=0}a_n(\tau)A_{n}(z;\tau)\sum^{\infty}_{l=1}C_l(z;\tau)=
$$
$$
=\left(\sum^{\infty}_{n=0}a_n(\tau)A_{n}(z;\tau)\right)\widehat{y}(z;\tau).
$$
Hence we can write
$$
\sum^{\infty}_{n=0}a_n(\tau)\left(\widehat{y}(z+n\tau;\tau)-A_{n}(z;\tau)\widehat{y}(z;\tau)\right)=0
$$
and it must be
$$
\widehat{y}(z+n\tau;\tau)=A_{n}(z;\tau)\widehat{y}(z;\tau).\eqno{(61)}
$$
Hence if $a_n(\tau)=c_n(\tau)$ and
$$
g(z;\tau)=B(iz;\tau)=\sum^{\infty}_{n=0}b_n(\tau)(iz)^n=\sum^{\infty}_{n=0}c_n(\tau)A_n(z;\tau),
$$
then from (d) $\widehat{y}(z)$ is given from (f) and satisfies (56.1a) and $y(x;\tau)$ is solution of 
$$
\sum^{\infty}_{n=0}b_n(\tau)y^{(n)}(x;\tau)=V(x;\tau)y(x;\tau),\eqno{(62)}
$$
where 
$$
V(x;\tau)=\sum^{\infty}_{n=0}c_n(\tau)e^{-in\tau x}.
$$

\section{A relativistic wave equation}

\[
\]

We have
$$
\left[E,U\right]\Psi=\left(EU\Psi-UE\Psi\right)=i\hbar\partial_t\left(U\Psi\right)-i\hbar U\partial_t\Psi
=i\hbar\left(\partial_t U\right)\Psi.
$$
Hence\\
\\
\textbf{Lemma 1.}\\
It holds in general (for any operator $U$ and any function $\Psi$):
$$
i\hbar\left(\partial_t U\right)\Psi=\left[E,U\right]\Psi.
$$
\\
\textbf{Theorem 11.}\\
Assume the equation
$$
E\Psi=H\Psi\eqno{(eq)}
$$
Assume operator $U$ satisfies the equation
$$
i\hbar\left(\partial_tU\right)=\left[H,U\right],\eqno{(eq1)}
$$
where $H$ is the self adjoint Hamiltonian of $(eq)$. Then\\
\textbf{I.} $U$ sends a solution $\Psi$ of $(eq)$ to $U\Psi$, which is also a solution of $(eq)$.\\
\textbf{II.} $U$ is almost self adjoint i.e. for every $\Psi,\Phi$ solutions of $(eq)$ exists constant $c$ such that:
$$
\left(U\Phi,\Psi\right)-\left(\Phi,U\Psi\right)=c=\textrm{constant}.
$$
\textbf{III.} The operator $U$ is almost self-adjoint operator and also invariant i.e. is a conserved quantity of the system.\\
\textbf{Remarks.} Equation $(eq1)$ is equivalent to $(\textbf{I})$.\\
\\
\textbf{Proof.}\\ 
It holds for every operator $U$ the relation
$$
i\hbar\left(\partial_tU\right)=\left[E,U\right].
$$
Hence from $(eq1)$, we get easily that if $\Psi$ is solution of $E\Psi=H\Psi$, then
$$
\left[U,H\right]\Psi+\left[E,U\right]\Psi=0\Leftrightarrow
$$
$$
UH\Psi-HU\Psi+EU\Psi-UE\Psi=0\Leftrightarrow
$$
$$
EU\Psi=HU\Psi.
$$
Hence $\Psi'=U\Psi$ is clearly solution of $(eq)$. The inverse of this is also easily proved.\\ If $\Phi$ and $\Psi$ are solutions of $(eq)$, then differentiating $(U\Phi,\Psi)-(\Phi,U\Psi)$ with respect to $t$, we arive equivalently to the equation 
$$
\partial_t\left(\left(U\Phi,\Psi\right)-\left(\Phi,U\Psi\right)\right)=0. 
$$
Hence easily we conclude that
$$
\left(U\Phi,\Psi\right)-\left(\Phi,U\Psi\right)=c=\textrm{constant.}
$$
$U$ is almost self adjoint.\\
It is also easy to see someone that $i\hbar\frac{d}{dt}\left\langle U\right\rangle=0$. Hence $U$ is a conserved quantity.\\
\\
\textbf{Examples.}\\
Assume that the Hamiltonian is independent of time $t$ (for example $H=\left(p\right)^2/(2m)+V(x)$). Then from 
$$
i\hbar\partial_tH=0=\left[H,H\right]=0,
$$
we get that if $\Psi$ is solution of $(eq)$, then $H\Psi$ is also a solution of $(eq)$. Also then 
$$
\left\langle E\right\rangle=\left\langle H\right\rangle=\textbf{constant.}
$$ 
Hence the energy is conserved. The same holds if we replace 
$$
U=f(H)=\sum^{\infty}_{k=0}a_k H^k.
$$
Then if $\Psi$ is solution of $(eq)$, the quantity $f(H)\Psi$ is also solution of $(eq)$. Also then  $f(H)$ is almost self adjoint and conserved quantity. Note that
$$
i\hbar\left(\partial_tf(H)\right)=0=\left[H,f(H)\right].
$$
\\

The Lorentz transformation is
$$
x_{k}=\sum^{2}_{\nu=1}\Lambda_{k,\nu} x'_{\nu},\eqno{:(\Lambda)}
$$
where $x_1=x$, $x_2=t$, $x'_1=x'$, $x'_2=t'$ and $\Lambda_{1,1}=\gamma_0$, $\Lambda_{1,2}=-\gamma_0 v$, $\Lambda_{2,1}=-\gamma_0 v/c^2$, $\Lambda_{2,2}=\gamma_0$, $\gamma_0=\left(1-v^2/c^2\right)^{-1/2}$.\\

Assume that we have two equations
\begin{equation}
i\hbar\frac{\partial}{\partial t}\Psi(x,t)=H\Psi(x,t)\textrm{ and }i\hbar\frac{\partial}{\partial x}\Psi(x,t)=W\Psi(x,t).
\end{equation}
Here $x,t$ are the variables and $\Psi=\Psi(x,t)$ is the solution of our system. We now introduce a Lorentz boost along $x$:
\begin{equation}
x=\gamma_0 x'+\gamma_0 v t'\textrm{ and }t=\gamma_0v c^{-2}x'+\gamma_0 t',
\end{equation}
where $\gamma_0=(1-v^2/c^2)^{-1/2}$. The inverse of this transformation is
\begin{equation}
x'=\gamma_0x-\gamma_0vt\textrm{ and }t'=-\gamma_0vc^{-2}x+\gamma_0 t
\end{equation}
 The derivatives of all functions $F(x,t)$ are obeying
\begin{equation}
\frac{\partial}{\partial t}=-\gamma_0 v\frac{\partial}{\partial x'}+\gamma_0\frac{\partial}{\partial t'}
\end{equation}
and
\begin{equation}
\frac{\partial}{\partial x}=\gamma_0\frac{\partial}{\partial x'}-\gamma_0 v c^{-2}\frac{\partial}{\partial t'}
\end{equation}
and the inverse is
\begin{equation}
\frac{\partial}{\partial t'}=\gamma_0 v\frac{\partial}{\partial x}+\gamma_0\frac{\partial}{\partial t}
\end{equation}
and
\begin{equation}
\frac{\partial}{\partial x'}=\gamma_0\frac{\partial}{\partial x}+\gamma_0vc^{-2}\frac{\partial}{\partial t}.
\end{equation}
Before the Lorentz transformation our equations where that of (55). After, they become
\begin{equation}
i\hbar\frac{\partial}{\partial t'}\Psi'(x',t')=H'\Psi'(x',t')\textrm{ and }i\hbar\frac{\partial}{\partial x'}\Psi'(x',t')=W'\Psi'(x',t').
\end{equation}
We define the operator $S$ which sends $\Psi(x,t)$ to $\Psi'(x',t')$, such:
\begin{equation}
\Psi'(x',t')=S \Psi(x,t).
\end{equation}
Also we define the opertors $S^{*}, T$ such that:
\begin{equation}
S\Psi(x,t)=\Psi'(x',t')=S^{*}\Psi(x',t')=S^{*}T\Psi(x,t),
\end{equation}
where
\begin{equation}
S=S^{*}T.
\end{equation}
At this point we also define the operators $\Gamma,G$ such that $\Gamma=i\hbar\frac{\partial}{\partial T}$ and the equation $\Gamma \Psi=G\Psi\Leftrightarrow\Psi\in L$. Then if $T,U$ are any operators, we have
\begin{equation}
i\hbar\left(\frac{\partial U}{\partial T}\right)\Phi=\left[\Gamma,U\right]\Phi\textrm{, }\forall \Phi\in M=L_2(R).
\end{equation}
The new generalized operator (it sends the operators $\Gamma,P$ to an operator $(A^{*}_TP)$):
\begin{equation}
\left(A^{*}_{T}P\right):=i\hbar\left(\frac{\partial P}{\partial T}\right)+\left[P,G\right].
\end{equation}
For example
\begin{equation}
\left(A^{*}_tP\right):=i\hbar\left(\frac{\partial P}{\partial t}\right)+\left[P,H\right]
\end{equation}
and
\begin{equation}
\left(A^{*}_xP\right):=i\hbar\left(\frac{\partial P}{\partial x}\right)+\left[P,W\right],
\end{equation}
for all operators $P$. We define also (from the equation $\Gamma \Psi=G\Psi$) the operators:
\begin{equation}
\left(A_TP\right)\Phi:=i\hbar\left(\frac{\partial P}{\partial T}\right)\Phi+P(G\Phi)
\end{equation}
and
\begin{equation}
\left(A_{T'}P\right)\Phi:=i\hbar\left(\frac{\partial P}{\partial T'}\right)\Phi+P(G'\Phi).
\end{equation}
Hence if $E\Psi=H\Psi$, then
\begin{equation}
\left(A_tP\right)\Phi=i\hbar\left(\frac{\partial P}{\partial t}\right)\Phi+P(H\Phi)
\end{equation}
and
\begin{equation}
\left(A_{t'}P\right)\Phi=i\hbar\left(\frac{\partial P}{\partial t'}\right)\Phi+P(H'\Phi).
\end{equation}
In the same way we have (with $\Psi$ satisfying $-p\Psi=W\Psi$):
\begin{equation}
\left(A_xP\right)\Psi=i\hbar\left(\frac{\partial P}{\partial x}\right)\Psi+P(W\Psi)
\end{equation}
and
\begin{equation}
\left(A_{x'}P\right)\Psi=i\hbar\left(\frac{\partial P}{\partial x'}\right)\Psi+P(W'\Psi).
\end{equation}
Hence
$$
\left(A_tP\right)\Psi=i\hbar\left(\frac{\partial P}{\partial t}\right)\Psi+P(H\Psi)=\left(EP\right)\Psi+P\left(H\Psi\right)=
$$
$$
=E\left(P\Psi\right)-P\left(E\Psi\right)+P\left(H\Psi\right)
=E\left(P\Psi\right).
$$
Also
$$
\left(A_xP\right)\Psi=i\hbar\left(\frac{\partial P}{\partial x}\right)\Psi+P(W\Psi)=-\left(pP\right)\Psi+P\left(W\Psi\right)=
$$
$$
=-p\left(P\Psi\right)+P\left(p\Psi\right)+P\left(W\Psi\right)
=-p\left(P\Psi\right).
$$
Also (see relations (68),(69)):
\begin{equation}
\left(A_{x'}P\right)\Psi=-p'\left(P\Psi\right)\textrm{, }p'=-i\hbar\frac{\partial}{\partial x'}
\end{equation}
and
\begin{equation}
\left(A_{t'}P\right)\Psi=E'\left(P\Psi\right)=i\hbar\frac{\partial}{\partial t'}\left(P\Psi\right).
\end{equation}
Also $A_t$ sends an operator $P$ to operator $A_tP$:
\begin{equation}
A_tP=i\hbar\left(\frac{\partial P}{\partial t}\right)+PH
\end{equation}
and
\begin{equation}
A_xP=i\hbar\left(\frac{\partial P}{\partial x}\right)+PW.
\end{equation}
We can also say that (under equation $E\Psi=H\Psi$ and $E'\Psi'=H'\Psi'$):
\begin{equation}
\left(A_{t'}P\right)\Psi'=i\hbar\left(\frac{\partial P}{\partial t'}\right)\Psi'+P\left(H'\Psi'\right)
\end{equation}
and
\begin{equation}
\left(A_{x'}P\right)\Psi'=i\hbar\left(\frac{\partial P}{\partial x'}\right)\Psi'+P\left(W'\Psi'\right).
\end{equation}

Then we can write
\begin{equation}
i\hbar\frac{\partial }{\partial t}\Psi=H\Psi\Leftrightarrow i\hbar\frac{\partial}{\partial t}\left(S^{-1}\Psi'\right)=H\left(S^{-1}\Psi'\right).
\end{equation}
Hence
$$
I=i\hbar\frac{\partial }{\partial t}\Psi=i\hbar\frac{\partial}{\partial t}\left(S^{-1}\Psi'\right)=i\hbar\left(\frac{\partial}{\partial t}S^{-1}\right)\Psi'+S^{-1}i\hbar\frac{\partial}{\partial t}\Psi'=
$$
$$
=i\hbar \gamma_0\left(\frac{\partial S^{-1}}{\partial t'}\right)\Psi'-\gamma_0 v i\hbar \left(\frac{\partial S^{-1}}{\partial x'}\right)\Psi'+\gamma_0i\hbar S^{-1}\frac{\partial \Psi'}{\partial t'}-i\hbar\gamma_0 v S^{-1}\frac{\partial \Psi'}{\partial x'}=
$$
$$
=i\hbar \gamma_0\left(\frac{\partial S^{-1}}{\partial t'}\right)\Psi'-\gamma_0 v i\hbar \left(\frac{\partial S^{-1}}{\partial x'}\right)\Psi'+\gamma_0 S^{-1}H'\Psi'-\gamma_0 v S^{-1}W'\Psi'=
$$
$$
=\gamma_0\left(i\hbar\left(\frac{\partial S^{-1}}{\partial t'}\right)+S^{-1}H'\right)\Psi'-\gamma_0 v \left(i\hbar \left(\frac{\partial S^{-1}}{\partial x'}\right)+S^{-1}W'\right)\Psi'=
$$
$$
=\gamma_0\left(A_{t'}S^{-1}\right)\Psi'-\gamma_0v\left(A_{x'}S^{-1}\right)\Psi'=
$$
$$
=\gamma_0\left(i\hbar\left(\frac{\partial S^{-1}}{\partial t'}\right)+\left(S^{-1}H'-H'S^{-1}\right)\right)\Psi'+\gamma_0H'S^{-1}\Psi'-
$$
$$
-\gamma_0 v\left(i\hbar\left(\frac{\partial S^{-1}}{\partial x'}\right)+\left(S^{-1}W'-W'S^{-1}\right)\right)\Psi'-\gamma_0 vW'S^{-1}\Psi'=
$$
\begin{equation}
=\gamma_0\left(A^{*}_{t'}S^{-1}\right)\Psi'-\gamma_0 v\left(A^{*}_{x'}S^{-1}\right)\Psi'+\left(\gamma_0H'-\gamma_0v W'\right)S^{-1}\Psi'=H\Psi.
\end{equation}
\\
\textbf{Theorem 12.0}\\
Assume the system of equations
\begin{equation}
E\Psi=H\Psi\textrm{ and }-p\Psi=W\Psi
\end{equation}
and denote $L$ the set of its solutions. Then under a Lorentz transform this becomes
\begin{equation}
E'\Psi'=H'\Psi'\textrm{ and }-p'\Psi'=W'\Psi'.
\end{equation}
We introduce the operators $f_{11},f_{22}$ such that
\begin{equation}
H'=\gamma_0f_{11}+v\gamma_0f_{22}
\end{equation}
and
\begin{equation}
W'=\gamma_0 v c^{-2}f_{11}+\gamma_0f_{22}.
\end{equation}
Then we have
\begin{equation}
f_{11}\Psi=(\gamma_0 H'-\gamma_0 v W')\Psi
\end{equation}
and
\begin{equation}
f_{22}\Psi=(\gamma_0 W'-\gamma_0 v c^{-2}H')\Psi.
\end{equation}
Then from (84),(85), we also have
\begin{equation}
E\Psi'=f_{11}\Psi'\textrm{ and }-p\Psi'=f_{22}\Psi'.
\end{equation}
\\
\textbf{Proof.}\\
From relation $E'\Psi'=H'\Psi'$, we have
$$
E'\Psi'=H'\Psi'\Leftrightarrow i\hbar\frac{\partial \Psi'}{\partial t'}=\left(\gamma_0 f_{11}+v\gamma_0 f_{22}\right)\Psi'\Leftrightarrow
$$
$$
i\hbar v\gamma_0\left(\frac{\partial \Psi'}{\partial x}\right)+i\hbar\gamma_0\left(\frac{\partial \Psi'}{\partial t}\right)=\left(\gamma_0 f_{11}+v\gamma_0 f_{22}\right)\Psi'.
$$
Also
$$
-p'\Psi'=W'\Psi'\Leftrightarrow i\hbar\frac{\partial \Psi'}{\partial x'}=\gamma_0 v c^{-2}f_{11}\Psi'+\gamma_0f_{22}\Psi'\Leftrightarrow
$$
$$
i\hbar\gamma_0\frac{\partial \Psi'}{\partial x}+i\hbar\gamma_0 v c^{-2}\frac{\partial \Psi'}{\partial t}=\gamma_0 v c^{-2}f_{11}\Psi'+\gamma_0f_{22}\Psi'.
$$
\\

For $f_{11},f_{22}$ see also Theorem 12.4\\
\\

Continuing from (83) we get
$$
B\Psi'+f_{11}\Psi=H\Psi\Leftrightarrow
$$
$$
B\Psi'+f_{11}\Psi=i\hbar\left(\frac{\partial S^{-1}}{\partial t}\right)\Psi'+S^{-1}\left(E\Psi'\right)=H\Psi,
$$
where
$$
f_{11}\Psi=\left(\gamma_0H'-\gamma_0 v W'\right)\Psi.
$$
Hence
\begin{equation}
B\Psi'=i\hbar\left(\frac{\partial S^{-1}}{\partial t}\right)\Psi'+S^{-1}\left(E\Psi'\right)-f_{11}\left(S^{-1}\Psi'\right)=\left(H-f_{11}\right)\Psi.
\end{equation}
However if we set
\begin{equation}
B\Psi'=i\hbar\left(\frac{\partial S^{-1}}{\partial t}\right)\Psi'+S^{-1}\left(H\Psi'\right)-H\left(S^{-1}\Psi'\right),
\end{equation}
we actually demand
\begin{equation}
S^{-1}\left(E\Psi'\right)-f_{11}\left(S^{-1}\Psi'\right)=S^{-1}\left(H\Psi'\right)-H\left(S^{-1}\Psi'\right)
\end{equation}
and
\begin{equation}
B\Psi'=H\Psi-f_{11}\Psi.
\end{equation}
Then we have equivalent:
$$
S^{-1}(E\Psi'-H\Psi')=(f_{11}-H)S^{-1}\Psi'\Leftrightarrow
$$
$$
(E-H)S\Psi=S(f_{11}-H)\Psi\Leftrightarrow
$$
$$
E\Psi'-H\Psi'=S(f_{11}\Psi)-S(H\Psi)\Leftrightarrow 
$$
$$
f_{11}\Psi'-H(S\Psi)=S(f_{11}\Psi)-S(H\Psi)\Leftrightarrow
$$
\begin{equation}
\left[H-f_{11},S\right]\Psi=0.
\end{equation}
Hence (under (93),(94)) the operators $B$, $S$, commute each other. Hence if there exists $S^{-1}$, then $B\Psi'=B(S\Psi)$ and
\begin{equation}
B(S\Psi)=B\Psi'=S(B\Psi).
\end{equation}
(proof)\\Set $B_0:=H-f_{11}$, then from above (relation (95)): $\left[B_0,S\right]=0$. Also (from (94)): $B\Psi'=B(S\Psi)=B_0\Psi$. Hence $B=B_0S^{-1}$ and $\left[B,S\right]\Psi=\left[B_0S^{-1},S\right]\Psi=B_0S^{-1}S\Psi-SB_0S^{-1}\Psi=0$, since $\left[B_0,S\right]=0$, in the set of all solutions of $E\Psi=H\Psi : (e1)$ and $-p\Psi=W\Psi : (e2)$. Hence (also in the set of solutions of $(e1),(e2)$), holds
\begin{equation}
B=\left(H-f_{11}\right)S^{-1}.
\end{equation}
$(qed)$.

Also holds from (91),(92),(93),(94):
$$
(A^{*}_tS)\Psi=i\hbar\left(\frac{\partial S}{\partial t}\right)\Psi-\left[H,S\right]\Psi=\left[E,S\right]\Psi-\left[H,S\right]\Psi=
$$
$$
=\left[E-H,S\right]\Psi=\left(E-H\right)S\Psi=(f_{11}-H)\Psi'=
$$
$$
=-S(H-f_{11})\Psi=-S(B\Psi')=-S\left(A^{*}_{t}S^{-1}\right)\Psi',
$$
where
$$
\left(A^{*}_{t}S^{-1}\right)\Psi'=i\hbar\left(\frac{\partial S^{-1}}{\partial t}\right)\Psi'+S^{-1}\left(H\Psi'\right)-H\left(S^{-1}\Psi'\right).
$$
We define
$$
\left(D_tP\right)\Psi:=i\hbar\left(\frac{\partial P}{\partial t}\right)\Psi-\left[H,P\right]\Psi.
$$
Also we have the identity
\begin{equation}
\left(A^{*}_{t}S^{-1}\right)\Psi'=-S^{-1}\left(D_tS\right)\Psi,
\end{equation}
which we can also write it as
\begin{equation}
\left(D_t S^{-1}\right)\Psi'=-S^{-1}\left(D_tS\right)\Psi.
\end{equation}
Note that operator $\left(D_tP\right)\Psi$ is that of the invariant set of the equation $E\Psi=H\Psi$, of Therorem 11, $(eq1)$. Also we have the invariant transformation (under (91)-(94)):
\begin{equation}
\left(A^{*}_tS^{-1}\right)\Psi'=\gamma_0\left(A^{*}_{t'}S^{-1}\right)\Psi'-\gamma_0 v\left(A^{*}_{x'}S^{-1}\right)\Psi'.
\end{equation}
Also we have the identity
\begin{equation}
i\hbar\left(\frac{\partial S}{\partial t}\right)\Psi=\left[H,S\right]\Psi-S\left(A^{*}_{t}S^{-1}\right)S\Psi.
\end{equation}
\\

In the same way as above 
\begin{equation}
-p\Psi=W\Psi,
\end{equation}
we get 
$$
-p\left(S^{-1}\Psi'\right)=W\left(S^{-1}\Psi'\right).
$$
Hence
$$
-p\left(S^{-1}\Psi'\right)=-\left(pS^{-1}\right)\Psi'-S^{-1}\left(p\Psi'\right)=i\hbar\left(\frac{\partial S^{-1}}{\partial x}\right)\Psi'+S^{-1}i\hbar\frac{\partial \Psi'}{\partial x}=
$$
$$
i\hbar\gamma_0\left(\frac{\partial S^{-1}}{\partial x'}\right)\Psi'-i\hbar\gamma_0vc^{-2}\left(\frac{\partial S^{-1}}{\partial t'}\right)\Psi'+S^{-1}\gamma_0 i\hbar\frac{\partial}{\partial x'}\Psi'-\gamma_0 v c^{-2}S^{-1}i\hbar\frac{\partial}{\partial t'}\Psi'=
$$
$$
i\hbar\gamma_0\left(\frac{\partial S^{-1}}{\partial x'}\right)\Psi'-i\hbar\gamma_0vc^{-2}\left(\frac{\partial S^{-1}}{\partial t'}\right)\Psi'+S^{-1}\gamma_0 W'\Psi'-\gamma_0 v c^{-2}S^{-1}H'\Psi'=
$$
$$
=\gamma_0\left(i\hbar\left(\frac{\partial S^{-1}}{\partial x'}\right)+S^{-1}W'\right)\Psi'-\gamma_0vc^{-2}\left(i\hbar\left(\frac{\partial S^{-1}}{\partial t'}\right)+S^{-1}H'\right)\Psi'=
$$
$$
=\gamma_0 \left(A_{x'}S^{-1}\right)\Psi'-\gamma_0v c^{-2}\left(A_{t'}S^{-1}\right)\Psi'=
$$
$$
=\gamma_0 \left(A^{*}_{x'}S^{-1}\right)\Psi'-\gamma_0v c^{-2}\left(A^{*}_{t'}S^{-1}\right)\Psi'+f_{22}\Psi=W\Psi,
$$
where (for $f_{22}$ see Theorem 12.0 and Theorem 12.4 relations (123),(124)):
$$
f_{22}\Psi=\gamma_0W'\Psi-\gamma_0vc^{-2}H'\Psi.
$$
Hence it holds
$$
-p\Psi=\Gamma\Psi'+f_{22}\Psi=W\Psi=i\hbar\left(\frac{\partial S^{-1}}{\partial x}\right)\Psi'+S^{-1}\left(-p\Psi'\right).
$$
Then
$$
\Gamma\Psi'=i\hbar\left(\frac{\partial S^{-1}}{\partial x}\right)\Psi'+S^{-1}(-p\Psi')-f_{22}\left(S^{-1}\Psi'\right).
$$
We want to show that
$$
\Gamma\Psi'=\left(A^{*}_{x}S^{-1}\right)\Psi'=i\hbar\left(\frac{\partial S^{-1}}{\partial x}\right)\Psi'+S^{-1}\left(W\Psi'\right)-W\left(S^{-1}\Psi'\right).
$$
Hence equivalent it must be
\begin{equation}
S^{-1}\left(-p\Psi'\right)-f_{22}\left(S^{-1}\Psi'\right)=S^{-1}\left(W\Psi'\right)-W\left(S^{-1}\Psi'\right)
\end{equation}
and
\begin{equation}
\Gamma \Psi'=W\Psi-f_{22}\Psi
\end{equation}
Hence equivalent
$$
\Gamma\Psi'=W\Psi-f_{22}\Psi=S^{-1}\left(W\Psi'\right)-S^{-1}\left(-p\Psi'\right)\Leftrightarrow
$$
$$
S\left(\Gamma \Psi'\right)=W\Psi'+p\Psi'=S\left(W\Psi-f_{22}\Psi\right)\Leftrightarrow
$$
$$
W(S\Psi)+p(S\Psi)=S(W\Psi)-S(f_{22}\Psi)\Leftrightarrow
$$
$$
-p(S\Psi)-W(S\Psi)=-S(W-f_{22})\Psi\Leftrightarrow
$$
$$
(-p-W)\Psi'=-S(W-f_{22})\Psi\Leftrightarrow
$$
$$
(f_{22}\Psi'-W\Psi')=-S(W-f_{22})\Psi\Leftrightarrow
$$
\begin{equation}
(W-f_{22})S\Psi=S(W-f_{22})\Psi.
\end{equation}
Hence we have equivalently
$$
\left[W-f_{22},S\right]\Psi=0.
$$
Hence one can easily see that $\Gamma$ and $S$ commute (as with $B$ relation (96) above). Also holds the identity (when $S^{-1}$ exists):
$$
\left(D_xS\right)\Psi=-S\Gamma\Psi'\Leftrightarrow \left(D_xS\right)\Psi=-S\left(W-f_{22}\right)\Psi
$$
and (as with $B$ above):
\begin{equation}
\Gamma=\left(W-f_{22}\right)S^{-1}.
\end{equation}
Hence\\
\\
\textbf{Theorem 12.1}\\
From the above notations and evaluations, we have $\Gamma=(W-f_{22})S^{-1}$ and
\begin{equation}
\left[\Gamma,S\right]\Psi=0, 
\end{equation}
iff
\begin{equation}
\Gamma\Psi'=-S^{-1}(D_x S)\Psi=(D_x S^{-1})\Psi',
\end{equation}
iff
\begin{equation}
(A^{*}_{x}S^{-1})\Psi'=\gamma_0(A^{*}_{x'}S^{-1})\Psi'-\gamma_0v c^{-2} (A^{*}_{t'}S^{-1})\Psi'.
\end{equation}
\\
\textbf{Theorem 12.2}\\
Also $B=(H-f_{11})S^{-1}$ and
\begin{equation}
\left[B,S\right]\Psi=0, 
\end{equation}
iff
\begin{equation}
B\Psi'=-S^{-1}(D_t S)\Psi=(D_t S^{-1})\Psi',
\end{equation}
iff
\begin{equation}
(A^{*}_{t}S^{-1})\Psi'=\gamma_0(A^{*}_{t'}S^{-1})\Psi'-\gamma_0 v (A^{*}_{x'}S^{-1})\Psi'.
\end{equation}
\\
\textbf{Theorem 12.3}\\
We define the operators
\begin{equation}
\left(D_tA\right)\Psi:=i\hbar\left(\frac{\partial A}{\partial t}\right)\Psi-\left[H,A\right]\Psi
\end{equation}
and
\begin{equation}
\left(D_xA\right)\Psi:=i\hbar\left(\frac{\partial A}{\partial x}\right)\Psi-\left[W,A\right]\Psi.
\end{equation}
Then $D_t$ and $D_x$ are derivatives, in the set of operators and holds:
\begin{equation}
\left(D_t (A B)\right)\Psi=\left(\left(D_t A\right)B\right)\Psi+\left(A \left(D_t B\right)\right)\Psi
\end{equation}
and
\begin{equation}
\left(D_x (AB)\right)\Psi=\left(\left(D_x A\right)B\right)\Psi+\left(A \left(D_x B\right)\right)\Psi.
\end{equation}
Also it holds under (95),(98):
\begin{equation}
\left(A^{*}_{t}S^{-1}\right)\Psi'=\gamma_0\left(A^{*}_{t'}S^{-1}\right)\Psi'-\gamma_0 v\left(A^{*}_{x'}S^{-1}\right)\Psi',
\end{equation}
\begin{equation}
\left(A^{*}_x S^{-1}\right)\Psi'=\gamma_0\left(A^{*}_{x'}S^{-1}\right)\Psi'-\gamma_0v c^{-2}\left(A^{*}_{t'}S^{-1}\right)\Psi'.
\end{equation}
\textbf{Remarks.}\\ \textbf{I.} Assume $A,S$ are any operators with $S$ invertible in the sense $S(S^{-1}\Psi)=S^{-1}(S\Psi)=I\Psi=\Psi$ and $\Psi$ any function, then
\begin{equation}
(A*(S*S^{-1}))\Psi=0.
\end{equation}
\textbf{II.} For simplicity we write the Lie product between two operators $(A*B)\Psi=\left[A,B\right]\Psi=(AB)\Psi$ and the composition product as $(A\circ B)\Psi=A(B\Psi)$. Hence $(A*B)\Psi=(AB)\Psi=A(B\Psi)-B(A\Psi)=(A\circ B)\Psi-(B\circ A)\Psi=A(B\Psi)-B(A\Psi)$.\\ 
\textbf{III.} Every operator $A$ is a differential operator in the space of operators. That is because
\begin{equation}
A(B\Psi)=\left[A,B\right]\Psi+B(A\Psi).
\end{equation}
\\
\textbf{Proof.}\\
It holds
$$
\left(\frac{\partial}{\partial t}\left(AB\right)\right)\Psi=\left(\left(\frac{\partial A}{\partial t}\right)B\right)\Psi+\left(A\left(\frac{\partial B}{\partial t}\right)\right)\Psi
$$
and
$$
\left(\frac{\partial}{\partial x}(AB)\right)\Psi=\left(\left(\frac{\partial A}{\partial x}\right)B\right)\Psi+\left(A\left(\frac{\partial B}{\partial x}\right)\right)\Psi.
$$
Also holds in general
$$
\left[C,AB\right]\Psi=\left(\left[C,A\right]B\right)\Psi+\left(A\left[C,B\right]\right)\Psi.
$$
Hence we get easily (115),(116).\\
The condition to hold
\begin{equation}
\left(A^{*}_x S^{-1}\right)\Psi'=\gamma_0\left(A^{*}_{x'}S^{-1}\right)\Psi'-\gamma_0v c^{-2}\left(A^{*}_{t'}S^{-1}\right)\Psi',
\end{equation}
is equivalent to:
\begin{equation}
i\hbar\left(\frac{\partial S}{\partial x}\right)\Psi-\left[W,S\right]\Psi=-S\left(A_x^{*}S^{-1}\right)S\Psi=(f_{22}-H)S\Psi.
\end{equation}
\\
\textbf{Remarks.}\\
We have for any operator $P$:
$$
\left(D_t\left(P^{-1}P\right)\right)\Psi=0\Leftrightarrow \left(\left(D_tP^{-1}\right)P\right)\Psi+\left(P^{-1}\left(D_tP\right)\right)\Psi=0\Leftrightarrow
$$
$$
\left[\left(D_tP^{-1}\right),P\right]\Psi=-\left[P^{-1},\left(D_tP\right)\right]\Psi\Leftrightarrow 
$$
$$
\left[\left(D_tP^{-1}\right),P\right]\Psi=\left[\left(D_tP\right),P^{-1}\right]\Psi.
$$
Hence if we set the operators
$$
U_t(P):=\left[\left(D_tP\right),P^{-1}\right]\textrm{, }\forall P\textrm{ operator, }
$$
and
$$
U_x(P):=\left[\left(D_xP\right),P^{-1}\right]\textrm{, }\forall P\textrm{ operator, }
$$
then
$$
U_t\left(P^{-1}\right)\Psi=U_{t}\left(P\right)\Psi,
$$
and
$$
U_x\left(P^{-1}\right)\Psi=U_x\left(P\right)\Psi,
$$
for all invertable operators $P$ and any function $\Psi$.\\
\\
\textbf{Theorem 12.4}\\
Assume a physical system described from the equations $E\Psi=H\Psi$ and $-p\Psi=W\Psi$. If we perform a Lorentz transformation $(\Lambda)$, the new equations become $E'\Psi'=H'\Psi'$ and $-p'\Psi'=W'\Psi'$. The operator $S$, which connects the system before the Lorentz transform and the system after the Lorentz transform, is described as $\Psi'=S\Psi$. Assuming also, that $H,W$ are transforming as $H'=\gamma_0f_{11}+v\gamma_0f_{12}$ and $W'=\gamma_0 v c^{-2}f_{21}+\gamma_0f_{22}$, it is easy to show that most general case happens if $f_{12}=f_{22}$ and $f_{21}=f_{11}$. Hence we have
\begin{equation}
H'\Phi=\gamma_0f_{11}\Phi+v\gamma_0f_{22}\Phi
\end{equation}
and 
\begin{equation}
W'\Phi=\gamma_0 v c^{-2}f_{11}\Phi+\gamma_0 f_{22}\Phi,
\end{equation}
for any function $\Phi$, in all Hilbert space $M$.\\
We define the operators $B,\Gamma$, as
\begin{equation}
B=(H-f_{11})S^{-1}\textrm{ and }\Gamma=(W-f_{22})S^{-1}\textrm{, resp}
\end{equation}
for any function in $S(L)$, (note that $L$ is the space of all solutions $\Psi$ of equations $(E-H)\Psi=0$ and $(-p-W)\Psi=0$). Then we have:\\
\textbf{I.} The invariant transformation
\begin{equation}
\left(A^{*}_tS^{-1}\right)\Psi'=\gamma_0\left(A^{*}_{t'}S^{-1}\right)\Psi'-\gamma_0 v\left(A^{*}_{x'}S^{-1}\right)\Psi'\textrm{, }\Psi'\in S(L),
\end{equation}
is equivalent to equation
\begin{equation}
\left[S,H-f_{11}\right]\Psi=0\textrm{, }\Psi\in L
\end{equation}
and this is equivalent to say that $B,S$ commute.\\
\textbf{II.} The invariant transformation
\begin{equation}
\left(A^{*}_x S^{-1}\right)\Psi'=-\gamma_0v c^{-2}\left(A^{*}_{t'}S^{-1}\right)\Psi'+\gamma_0\left(A^{*}_{x'}S^{-1}\right)\Psi',
\end{equation}
is equivalent to equation
\begin{equation}
\left[S,W-f_{22}\right]\Psi=0\textrm{, }\Psi\in L
\end{equation}
and this is equivalent to say that $\Gamma$, $S$ commute in $L$.\\
\\
\textbf{Proof.}\\
From relation $E'\Psi'=H'\Psi'$, we have
$$
E'\Psi'=H'\Psi'\Leftrightarrow i\hbar\frac{\partial \Psi'}{\partial t'}=\left(\gamma_0 f_{11}+v\gamma_0 f_{22}\right)\Psi'\Leftrightarrow
$$
$$
i\hbar v\gamma_0\left(\frac{\partial \Psi'}{\partial x}\right)+i\hbar\gamma_0\left(\frac{\partial \Psi'}{\partial t}\right)=\left(\gamma_0 f_{11}+v\gamma_0 f_{22}\right)\Psi'.\eqno{(a)}
$$
Hence from $(a)$
$$
-p\Psi'=f_{22}\Psi'.
$$
Continuing our arguments, we have the inverse (see Theoem 12.5 below):
$$
i\hbar\frac{\partial (S\Psi)}{\partial x}=f_{22}(S\Psi)\Leftrightarrow i\hbar\left(\frac{\partial S}{\partial x}\right)\Psi+i\hbar S \left(\frac{\partial \Psi}{\partial x}\right)=f_{22}(S\Psi)\Leftrightarrow
$$
$$
\left[W,S\right]\Psi-S\left(A^{*}_{x}S^{-1}\right)S\Psi+S(W\Psi)=f_{22}(S\Psi)\Leftrightarrow
$$
$$
-S\left(A^{*}_{x}S^{-1}\right)\Psi'+\left(W-f_{22}\right)\Psi'=0\Leftrightarrow
$$
$$
\left(D_xS\right)\Psi+\left(W-f_{22}\right)\Psi'=0\Leftrightarrow
$$
$$
\left(f_{22}-W\right)S\Psi+\left(W-f_{22}\right)\Psi'=0\Leftrightarrow 0=0\textrm{, true. }
$$
Hence indeed the solution of 
$$
i\hbar\frac{\partial}{\partial x}\Psi'=f_{22}\Psi'
$$
is 
$$
\Psi'=S\Psi.
$$
Also
$$
S\left(\Gamma \Psi'\right)=S\left(W-f_{22}\right)\Psi=\left(W-f_{22}\right)\Psi'=\left(W+p\right)\Psi'.
$$
Also
$$
i\hbar\frac{\partial \left(S\Psi\right)}{\partial t}=f_{11}(S\Psi)\Leftrightarrow
$$
$$
i\hbar\left(\frac{\partial S}{\partial t}\right)\Psi+i\hbar S\left(\frac{\partial \Psi}{\partial t}\right)=f_{11}S\Psi\Leftrightarrow
$$
$$
\left[H,S\right]\Psi-S\left(A^{*}_t S^{-1}\right)S\Psi+S(H\Psi)=f_{11}(S\Psi)\Leftrightarrow
$$
$$
H(S\Psi)-S(H\Psi)+(D_{t}S)\Psi+S(H\Psi)=f_{11}(S\Psi)\Leftrightarrow
$$
$$
H(S\Psi)-f_{11}(S\Psi)=(H-f_{11})S\Psi\Leftrightarrow
$$
$$
 (H-f_{11})S\Psi=(H-f_{11})S\Psi.
$$
\\
\textbf{Theorem 12.5}\\
Under (123),(124) and $E\Psi=H\Psi$, $-p\Psi=W\Psi$, it always holds:
\begin{equation}
(D_tS)\Psi=-S(D_tS^{-1})S\Psi\textrm{ and }(D_xS)\Psi=-S(D_xS^{-1})S\Psi.
\end{equation}
Also then we have
\begin{equation}
(D_tS)\Psi=(f_{11}-H)S\Psi=(E-H)S\Psi
\end{equation}
and
\begin{equation}
(D_xS)\Psi=(f_{22}-W)S\Psi=(-p-W)S\Psi.
\end{equation}
\\
\textbf{Proof}\\
We have
$$
E\Psi=H\Psi\textrm{ and }-p\Psi=W\Psi\eqno{(eq)}
$$ 
and its Lorentz 
$$
E'\Psi'=H'\Psi'\textrm{ and }-p'\Psi'=W'\Psi'.\eqno{(eq')}
$$
Assume that $f_{11}$ and $f_{22}$ are such that 
\begin{equation}
H'=\gamma_0f_{11}+v\gamma_0f_{22}
\end{equation}
and 
\begin{equation}
W'=\gamma_0 v c^{-2}f_{11}+\gamma_0 f_{22}.
\end{equation}
Then
\begin{equation}
E\Psi'=f_{11}\Psi'
\end{equation}
and
\begin{equation}
-p\Psi'=f_{22}\Psi',
\end{equation}
where $\Psi'=S\Psi$. 
That is because $(eq')$ can be writen as 
$$
i\hbar\frac{\partial \Psi'}{\partial t'}=(\gamma_0f_{11}+v\gamma_0f_{22})\Psi'
$$
and
$$
i\hbar\frac{\partial \Psi'}{\partial x'}=(\gamma_0v c^{-2}f_{11}+\gamma_0f_{22})\Psi'.
$$
Hence writing the partial derivatives in $t$ and $x$ (see relations (58)-(61)) we get:
$$
i\hbar\left(\gamma_0v\frac{\partial}{\partial x}+\gamma_0\frac{\partial}{\partial t}\right)\Psi'=(\gamma_0f_{11}+v\gamma_0f_{22})\Psi'
$$
and
$$
i\hbar\left(\gamma_0\frac{\partial}{\partial x}+\gamma_0v c^{-2}\frac{\partial}{\partial t}\right)\Psi'=\left(\gamma_0 v c^{-2}f_{11}+\gamma_0 f_{22}\right)\Psi'.
$$
Hence we get (133),(134). Then also must hold 
$$
E(S\Psi)=f_{11}(S\Psi)\Leftrightarrow (ES)\Psi+S(E\Psi)=f_{11}(S\Psi)\Leftrightarrow
$$
$$
i\hbar\left(\frac{\partial S}{\partial t}\right)\Psi+S(H\Psi)=f_{11}(S\Psi)\Leftrightarrow
$$
$$
\left(D_tS\right)\Psi=(f_{11}-H)S\Psi.\eqno{(a)}
$$
Also in the same way we have
$$
(D_xS)\Psi=(f_{22}-W)S\Psi.\eqno{(b)}
$$
Set now $A=-p-W$, then
$$
(D_xS)\Psi=i\hbar\left(\frac{\partial S}{\partial x}\right)\Psi-\left[W,S\right]\Psi=-(pS)\Psi-\left[W,S\right]\Psi=
$$
$$
=-\left[p,S\right]\Psi-\left[W,S\right]\Psi=-\left[p+W,S\right]\Psi=\left[A,S\right]\Psi=
$$
$$
=-(p+W)S\Psi-S(-p-W)\Psi=-(p+W)(S\Psi)=A(S\Psi).
$$
Hence
$$
(D_xS)\Psi=A(S\Psi)\eqno{(c)}
$$
In the same way we have
$$
(D_xS^{-1})\Psi=A(S^{-1}\Psi)
$$
and
$$
(D_xS^{-1})\Psi'=i\hbar\left(\frac{\partial S^{-1}}{\partial x}\right)\Psi'-\left[W,S^{-1}\right]\Psi'=
$$
$$
=-(pS^{-1})\Psi'-\left[W,S^{-1}\right]\Psi'=-\left[p,S^{-1}\right]\Psi'-\left[W,S^{-1}\right]\Psi'=
$$
$$
=-\left[p+W,S^{-1}\right]\Psi'=-(p+W)(S^{-1}\Psi')+S^{-1}(p+W)\Psi'=-S^{-1}(A\Psi').
$$
Hence
$$
(D_xS^{-1})\Psi'=-S^{-1}A\Psi'.\eqno{(d)}
$$
Finaly we have from $(c)$:
$$
-S^{-1}(D_xS)\Psi=-S^{-1}(A(S\Psi))
$$
and from $(d)$:
$$
(D_xS^{-1})S\Psi=-S^{-1}(A\Psi')=-S^{-1}(A(S\Psi)).
$$
Hence
\begin{equation}
(D_xS^{-1})S\Psi=-S^{-1}(D_xS)\Psi.
\end{equation}
In the same way one can obtain
\begin{equation}
(D_tS^{-1})S\Psi=-S^{-1}(D_tS)\Psi.
\end{equation}
QED.\\
\\
\textbf{Theorem 12.6}\\
If we have the system of equations 
$$
E\Psi=H\Psi\textrm{ and }-p\Psi=W\Psi.\eqno{(eq)}
$$
Under a Lorentz transform $x=\sigma(x',t')$, $t=\mu(x',t')$, the $(eq)$ become
\begin{equation}
E'\Psi'=H'\Psi'\textrm{ and }-p'\Psi'=W'\Psi'.
\end{equation}
Assume now that exist operators $f_{11},f_{22}$, such that
\begin{equation}
H'=\gamma_0f_{11}+v\gamma_0f_{22},
\end{equation} 
and
\begin{equation}
W'=\gamma_0 v c^{-2}f_{11}+\gamma_0 f_{22}.
\end{equation}
Then
\begin{equation}
E\Psi'=f_{11}\Psi'\textrm{ and }-p\Psi'=f_{22}\Psi',
\end{equation}
where $\Psi'(x',t')=S\Psi(x,t)$ and the following conditions are equivalent:\\
\textbf{I.}
\begin{equation}
\left(D_t S^{-1}\right)\Psi'=\gamma_0\left(D_{t'}S^{-1}\right)\Psi'-\gamma_0 v\left(D_{x'}S^{-1}\right)\Psi',
\end{equation}
\begin{equation}
\left(D_x S^{-1}\right)\Psi'=-\gamma_0v c^{-2}\left(D_{t'}S^{-1}\right)\Psi'+\gamma_0\left(D_{x'}S^{-1}\right)\Psi'.
\end{equation}
\textbf{II.}
\begin{equation}
\left[H-f_{11},S\right]\Psi=0\textrm{ and }\left[W-f_{22},S\right]\Psi=0.
\end{equation}
\textbf{III.}\\
If we define the operators $B,\Gamma$, as
\begin{equation}
B=(H-f_{11})S^{-1}\textrm{ and }\Gamma=(W-f_{22})S^{-1},
\end{equation}
then 
\begin{equation}
\left[B,S\right]\Psi=0\textrm{ and }\left[\Gamma,S\right]\Psi=0\textrm{, (resp.)}
\end{equation}
\\

If $\Psi^{*}=G\Psi'$, then
\begin{equation}
\left(D_{x'}G^{-1}\right)\Psi^{*}=G^{-1}\left(p'+W'\right)\Psi^{*},
\end{equation}
\begin{equation}
\left(D_{t'}G^{-1}\right)\Psi^{*}=-G^{-1}\left(E'-H'\right)\Psi^{*}.
\end{equation}
After some calculations we find
\begin{equation}
\gamma_0\left(D_{x'}G^{-1}\right)\Psi^{*}-\gamma_0 v c^{-2}\left(D_{t'}G^{-1}\right)\Psi^{*}=G^{-1}\left(p+f_{22}\right)\Psi^{*}
\end{equation}
and
\begin{equation}
\gamma_0\left(D_{t'}G^{-1}\right)\Psi^{*}-\gamma_0 v\left(D_{x'}G^{-1}\right)\Psi^{*}=-G^{-1}\left(E-f_{11}\right)\Psi^{*}.
\end{equation}
Hence
$$
\gamma_0\left(D_{x'}G^{-1}\right)\Psi^{*}-\gamma_0 v c^{-2}\left(D_{t'}G^{-1}\right)\Psi^{*}=
$$
$$
=\gamma_0G^{-1}\left(p'+W'\right)\Psi^{*}+\gamma_0 v c^{-2}G^{-1}\left(E'-H'\right)\Psi^{*}=
$$
\begin{equation}
=G^{-1}\left(p+f_{22}\right)\Psi^{*}.
\end{equation}
Also
$$
\left(D_{x}G^{-1}\right)\Psi^{*}=i\hbar\left(\frac{\partial G^{-1}}{\partial x}\right)\Psi^{*}-\left[W,G^{-1}\right]\Psi^{*}=\left[-p-W,G^{-1}\right]\Psi^{*}=
$$
\begin{equation}
=-\left(p+W\right)\Psi'+G^{-1}\left(p+W\right)\Psi^{*}.
\end{equation}
Hence for to have
\begin{equation}
\gamma_0\left(D_{x'}G^{-1}\right)\Psi^{*}-\gamma_0 v c^{-2}\left(D_{t'}G^{-1}\right)\Psi^{*}=\left(D_xG^{-1}\right)\Psi^{*},
\end{equation}
it must be equivalent
$$
G^{-1}\left(p+f_{22}\right)\Psi^{*}=G^{-1}(p+W)\Psi^{*}-(p+W)G^{-1}\Psi^{*}\Leftrightarrow
$$
$$
G^{-1}f_{22}\Psi^{*}=-(p+W)\Psi'+G^{-1}(W\Psi^{*})=-p\Psi'-W\Psi'+G^{-1}WG\Psi'\Leftrightarrow
$$
$$
G^{-1}\left(f_{22}-W\right)\Psi^{*}=-(p+W)\Psi'=(f_{22}-W)\Psi'\Leftrightarrow
$$
$$
 (f_{22}-W)G\Psi'=G(f_{22}-W)\Psi'\Leftrightarrow\left[W-f_{22},G\right]\Psi'=0.
$$
In the same way as above, for to have
\begin{equation}
\left(D_t G^{-1}\right)\Psi^{*}=\gamma_0\left(D_{t'}G^{-1}\right)\Psi^{*}-\gamma_0 v\left(D_{x'}G^{-1}\right)\Psi^{*},
\end{equation}
it must be equivalent
\begin{equation}
\left[H-f_{11},G\right]\Psi'=0.
\end{equation}
Hence\\
\\
\textbf{Theorem 12.7}\\
Assume that $(E-H)\Psi=(-p-W)\Psi=0$ and $(E'-H')\Psi'=(-p'-W')\Psi'=0$ (is the Lorentz). Also assume that exists operators $f_{11}$ and $f_{22}$, such that
\begin{equation}
H'\Phi=\gamma_0f_{11}\Phi+v\gamma_0f_{22}\Phi
\end{equation}
and 
\begin{equation}
W'\Phi=\gamma_0 v c^{-2}f_{11}\Phi+\gamma_0 f_{22}\Phi,
\end{equation}
for any function $\Phi$ in all Hilbert space $M$. Also $\Psi^{*}=G\Psi'$ and $\Psi'=S\Psi$. Then
\begin{equation}
\left(D_{t'}G^{-1}\right)\Psi^{*}=\gamma_0\left(D_tG^{-1}\right)\Psi^{*}+v\gamma_0\left(D_xG^{-1}\right)\Psi^{*}
\end{equation}
and
\begin{equation}
\left(D_{x'}G^{-1}\right)\Psi^{*}=\gamma_0 \left(D_{x}G^{-1}\right)\Psi^{*}+\gamma_0 v c^{-2}\left(D_{t}G^{-1}\right)\Psi^{*},
\end{equation}
iff
\begin{equation}
\left[H-f_{11},G\right]\Psi'=0
\end{equation}
and
\begin{equation}
\left[W-f_{22},G\right]\Psi'=0,
\end{equation}
respectively.

\[
\]

An example of such operators $f_{11}, f_{22}$, are if we set $f_{11}=\left(D_{t}G^{-1}\right)$ and $f_{22}=\left(D_xG^{-1}\right)$, where $G$ is almost any operator, then from (157),(158) we have in $G(S(L))$:
\begin{equation}
H'\Psi^{*}=\gamma_0\left(D_tG^{-1}\right)\Psi^{*}+\gamma_0 v\left(D_xG^{-1}\right)\Psi^{*},
\end{equation}
\begin{equation}
W'\Psi^{*}=\gamma_0 v c^{-2}\left(D_tG^{-1}\right)\Psi^{*}+\gamma_0\left(D_xG^{-1}\right)\Psi^{*}.
\end{equation}
Hence $H'=D_{t'}G^{-1}$ and $W'=D_{x'}G^{-1}$, in $G(S(L))$ and
\begin{equation}
\left[(D_{t}G^{-1})-H,G\right]\Psi'=0\textrm{ and }\left[(D_{x}G^{-1})-W,G\right]\Psi'=0.
\end{equation}
\\
\textbf{Proposition.}\\
Assume the notation of Theorems 12.5 and 13 above. Then
\begin{equation}
\left(D_{t'}G^{-1}\right)G\Psi'=-G^{-1}(E'-H')G\Psi'
\end{equation} 
and
\begin{equation}
\left(D_{t'}G\right)\Psi'=\left(E'-H'\right)G\Psi'.
\end{equation}
Hence
\begin{equation}
\left(D_{t'}G\right)\Psi'=-G\left(D_{t'}G^{-1}\right)G\Psi'.
\end{equation}
Also
\begin{equation}
\left(D_{x'}G\right)\Psi'=-(p'+W')G\Psi'
\end{equation}
and
\begin{equation}
\left(D_{x'}G^{-1}\right)G\Psi'=G^{-1}(p'+W')G\Psi'.
\end{equation}
Hence
\begin{equation}
\left(D_{x'}G\right)\Psi'=-G \left(D_{x'}G^{-1}\right)G\Psi'.
\end{equation}
\\
\textbf{Theorem 12.8}\\
Assume that $(E-H)\Psi=(-p-W)\Psi=0 : (eq)$ and $(E'-H')\Psi'=(-p'-W')\Psi'=0$ (is the Lorentz). 
Also when $\Psi,\Phi$ are solutions of $(eq)$ and\\
\textbf{1.}\\The inner product of the functions $h,g$ is
\begin{equation}
\left\langle h,g\right\rangle:=\int^{+\infty}_{-\infty}\overline{h(x,t)}g(x,t)dx.
\end{equation}
\textbf{2.} 
\begin{equation}
\partial_t\left(\int^{+\infty}_{-\infty}\overline{\Psi(x,t)}\Phi(x,t) dx\right)=0.
\end{equation} 
\textbf{3.}\\
Exist operators $f_{11}$ and $f_{22}$, such that
\begin{equation}
H'\Phi=\gamma_0f_{11}\Phi+v\gamma_0f_{22}\Phi
\end{equation}
and 
\begin{equation}
W'\Phi=\gamma_0 v c^{-2}f_{11}\Phi+\gamma_0 f_{22}\Phi,
\end{equation}
for any function $\Phi$ in all Hilbert space $M$.\\
Then if $L$ is the set of all solutions $\Psi$ of $(eq)$ (i.e. the section of the two equations) and $E=i\hbar\partial_{t}$, $p=-i\hbar\partial_{x}$, we have:\\
\textbf{a)}
\begin{equation}
E\Psi'=f_{11}\Psi'\textrm{ and }-p\Psi'=f_{22}\Psi',
\end{equation}
where $\Psi'(x',t')=S\Psi(x,t)$.\\
\textbf{b)}\\
The following three conditions are equivalent to each other:\\
\textbf{I.}
\begin{equation}
\left(D_t S^{-1}\right)\Psi'=\gamma_0\left(D_{t'}S^{-1}\right)\Psi'-\gamma_0 v\left(D_{x'}S^{-1}\right)\Psi',
\end{equation}
\begin{equation}
\left(D_x S^{-1}\right)\Psi'=-\gamma_0v c^{-2}\left(D_{t'}S^{-1}\right)\Psi'+\gamma_0\left(D_{x'}S^{-1}\right)\Psi'.
\end{equation}
\textbf{II.}
\begin{equation}
\left[H-f_{11},S\right]\Psi=0\textrm{ and }\left[W-f_{22},S\right]\Psi=0.
\end{equation}
\textbf{III.}\\
If we define the operators $B,\Gamma$, as
\begin{equation}
B=(H-f_{11})S^{-1}\textrm{ and }\Gamma=(W-f_{22})S^{-1},
\end{equation}
then 
\begin{equation}
\left[B,S\right]\Psi=0\textrm{ and }\left[\Gamma,S\right]\Psi=0\textrm{, (resp.)}
\end{equation}
\textbf{c)}\\
The system is a relativistic quantum mechanical system and 
\begin{equation}
\partial_t\left\langle S\right\rangle=0
\end{equation}
and (see Theorem 11):
\begin{equation}
\partial_t\{\left\langle S\Phi,\Psi\right\rangle-\left\langle\Phi,S\Psi\right\rangle\}=0
\end{equation}
and we say that $S$ is ''almost self adjoint''. However $S$ is a physical conserved quantity.

\newpage

\[
\]

\centerline{\bf References}

[1]: Wolfgang Rindler. ''Introduction to Special Relativity''. Oxford University Press. Inc. New York. (1991).

[2]: Edward Prugovecki. ''Quantum Mechanics in Hilbert Space''. Academic Press. New York, London. (1971).

[3]: Nikos Bagis. ''Some Results on the Theory of Infinite Series and Divisor Sums''. arXiv:0912.4815

[4]: Murray R. Spiegel. 'Schaum's Outline of Theory and Problems of Fourier Analysis with Applications to Boundary Value Problems'. McGraw-Hill, New York. 1974.

[5]: Frigyes Riesz and Bela Sz.-Nagy. 'Functional Analysis'. Dover Publications, Inc. New York. 1990.

[6]: Jerome A. Goldstein. 'Semigroups of Linear Operators and Applications'. Oxford University Press-New York. Clarendon Press-Oxford. 1985.

\end{document}